\documentclass{article}
\usepackage{epsfig}
\usepackage{calc}
\usepackage{amsfonts}
\usepackage{amsmath}
\usepackage{amssymb}
\usepackage{theorem}

\theorembodyfont{\rmfamily}
\newtheorem{definition}{Definition}[section]
\newtheorem{notation}[definition]{Notation}
\newtheorem{definition-notation}[definition]{Definition and notation}
\newtheorem{lemma}[definition]{Lemma}
\newtheorem{proposition}[definition]{Proposition}
\newtheorem{theorem}[definition]{Theorem}
\newtheorem{corollary}[definition]{Corollary}
\newtheorem{algorithm}[definition]{Algorithm}

\newcommand{\worklabel}[1]{\label{#1}}

\newcommand{\fig}[1]{\begin{figure}\begin{center}#1\end{center}\end{figure}}

\begin{document}

\title{
	\begin{minipage}{8cm} 
	An algorithm to generate exactly once every tiling with lozenges of a
	domain\footnote{
	This work has been partially supported by the University of
	Santiago, Chile; many thanks to Prof. \'Eric Gol\`es.}
	\end{minipage}
}

\vskip2cm

\author{
	S\'ebastien Desreux%
	\footnote{
		This work would not have been possible without the kind
		explanations, help and re-readings of \'Eric R\'emila,
		Michel Morvan and Jean-Christophe Novelli. 
		The author would also like to thank the 
		reviewers for their thorough evaluation, helpful comments 
		and pertinent suggestions.
	}
	\footnote{
		\textsc{Liafa}, Universit\'e Paris 7, case 7014,
		2 place Jussieu, 75\,251 Paris Cedex 05,\newline
		\texttt{Sebastien.Desreux@liafa.jussieu.fr}
	}
}

\maketitle

%%%%%%%%%%%%%%%%%%%%%%%%%%%%%%%%%%%%%%%%%%%%%%%%%%%%%%%%%%%%%%%%%%%%%%%
%%%%%%%%%%%%%%%%%%%%%%%%%%%%%%%%%%%%%%%%%%%%%%%%%%%%%%%%%%%%%%%%%%%%%%%
%%%%%%%%%%%%%%%%%%%%%%%%%%%%%%            %%%%%%%%%%%%%%%%%%%%%%%%%%%%%
%%%%%%%%%%%%%%%%%%%%%%%%%%%%%%   Résumé   %%%%%%%%%%%%%%%%%%%%%%%%%%%%%
%%%%%%%%%%%%%%%%%%%%%%%%%%%%%%            %%%%%%%%%%%%%%%%%%%%%%%%%%%%%
%%%%%%%%%%%%%%%%%%%%%%%%%%%%%%%%%%%%%%%%%%%%%%%%%%%%%%%%%%%%%%%%%%%%%%%
%%%%%%%%%%%%%%%%%%%%%%%%%%%%%%%%%%%%%%%%%%%%%%%%%%%%%%%%%%%%%%%%%%%%%%%

\begin{abstract}

We first show that the tilings of a domain ${\cal D}$ form a lattice (using
the same kind of arguments as in \cite{Remila}) which we then undertake to
decompose and generate without any redundance. To this end, we study
extensively the relatively simple case of hexagons and their deformations. We
show that general domains can be broken up into hexagon-like parts. Finally we
give an algorithm to generate exactly once every element in the lattice of the
tilings of a general domain ${\cal D}$.


\end{abstract}

%%%%%

\section{Introduction}

Tilings is an age-old topic for specialists and amateurs alike.  
Lozenge tilings in particular have intrigued generations of curious people 
because they can easily be seen as piles of cubes.

In the past few years, mathematics, theoretical physics and computer
science have started shedding a new light.  We will build on these results
to provide an answer to a very natural question: given a tileable domain
and sufficiently many lozenges tiles, how can one generate all the tilings
of the domain without repeating twice the same tiling?

Conway and Lagarias introduced in \cite{Conway} a new, powerful tool to
study tilings-related topics: tiling groups, which give a necessary
condition for a domain to be tileable, and provide and important
bijection in the case of lozenges. Thurston went a step further in
\cite{Thurston} and showed by a constructive algorithm that one can quickly
decide whether a domain is tileable; to this end he uses height functions
(see Section~\ref{basic:tools}).

Thurston also hinted that the set of the tilings of a domain, partially
ordered with the height functions, should have the structure of a lattice.  
This was proved in \cite{Remila} by R\'emila in 1999 in the case of
dominoes. We adapt his proof to the case of lozenges in
Section~\ref{flips:lattices}.

We then proceed to the really new material: after a rather extensive study of
the case of hexagonal and hexagonal-like domains (see
Section~\ref{hexagons:ph}), we use a geometrical point of view (justified by a
bijection between Conway and Lagarias' lozenge group and $\mathbb{Z}^3$) and
the identification of meaningful hexagonal-like sub-domains (see
Section~\ref{domains:fracture:lines:seeds}) to exhibit in the general case a
maximal chain of intervals in the lattice of the tilings (see
Section~\ref{lattices:intervals}). This chain is a natural extension to
Thurston's minimal and maximal tilings. We finally introduce new minimal
tilings which allow us to achieve our goal: (uniquely) generating all the
elements of the lattice formed by the tilings (see
Section~\ref{section:algo:main}).

We believe that the tools introduced in this paper, notably the seeds
and their ranges, should prove fruitful in tackling related problems
(see Section~\ref{conclusion:perspectives}).

%%%%%%%%%%%%%%%%%%%%%%%%%%%%%%%%%%%%%%%%%%%%%%%%%%%%%%%%%%%%%%%%%%%%%%%
%%%%%%%%%%%%%%%%%%%%%%%%%%%%%%%%%%%%%%%%%%%%%%%%%%%%%%%%%%%%%%%%%%%%%%%
%%%%%%%%%%%%%%%%%%%%%%%%%%%                  %%%%%%%%%%%%%%%%%%%%%%%%%%
%%%%%%%%%%%%%%%%%%%%%%%%%%%   Présentation   %%%%%%%%%%%%%%%%%%%%%%%%%%
%%%%%%%%%%%%%%%%%%%%%%%%%%%                  %%%%%%%%%%%%%%%%%%%%%%%%%%
%%%%%%%%%%%%%%%%%%%%%%%%%%%%%%%%%%%%%%%%%%%%%%%%%%%%%%%%%%%%%%%%%%%%%%%
%%%%%%%%%%%%%%%%%%%%%%%%%%%%%%%%%%%%%%%%%%%%%%%%%%%%%%%%%%%%%%%%%%%%%%%

\section{Basic tools and definitions}
\worklabel{basic:tools}

In this section, we present the definitions of classical objects which will be
used in this paper, in an attempt to make it reasonably self-contained. All
the new objects are defined later in the paper, at the moment when they are
needed.

%!!!!!!!!!!!!!!!!!!!!!!!!!!!                 !!!!!!!!!!!!!!!!!!!!!!!!!!!
%!!!!!!!!!!!!!!!!!!!!!!!!!!!   Définitions   !!!!!!!!!!!!!!!!!!!!!!!!!!!
%!!!!!!!!!!!!!!!!!!!!!!!!!!!                 !!!!!!!!!!!!!!!!!!!!!!!!!!!

\subsection{Tiling with lozenges}

First of all, let us define what we mean by a tiling.  One needs two regular
triangles (see Figure~\ref{triangles:lozenges}~(a)). The whole plane can be
covered with a repetition of these figures (see 
Figure~\ref{triangular:grid}~(a)), which gives rise to the \emph{triangular grid}. A
\emph{domain} is a finite union of triangles in the grid. It is \emph{simply
connected} if it is connected and its complement in the plane is connected. A
\emph{polygon} is a simply connected domain.

\fig{
	\begin{tabular}{p{3.5cm}@{\hskip1cm}p{6cm}}

	\mbox{}\hfill
	\epsfig{file=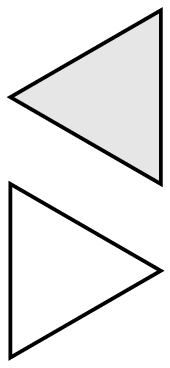, height=1.5cm}
	\hfill\mbox{}
	&
	\mbox{}\hfill
	\epsfig{file=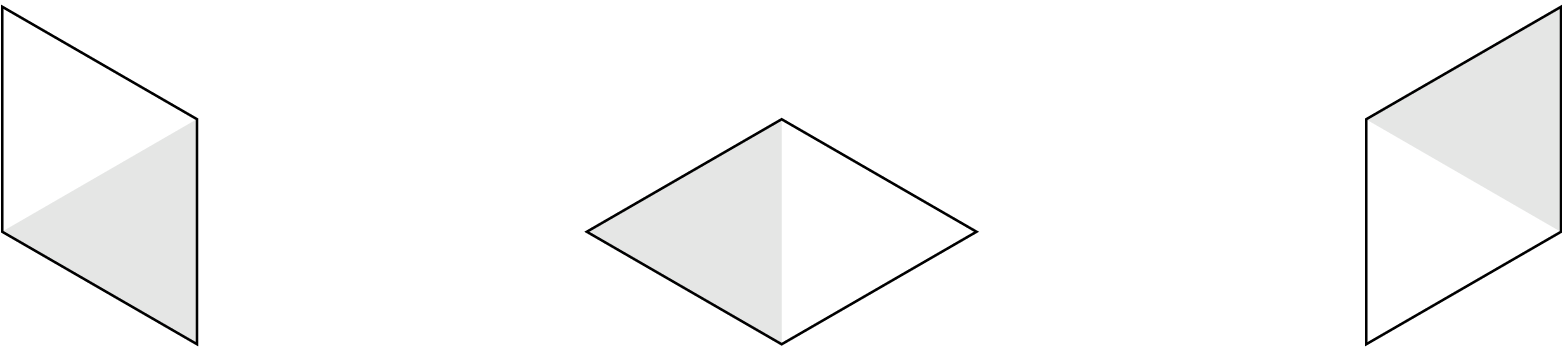, width=5.5cm, clip=}
	\hfill\mbox{}
	\\[1mm]

	\mbox{}\hfill
	(a) 	\begin{minipage}[t]{2.9cm}
		The two admissible regular triangles
		\end{minipage}
	\hfill\mbox{}
	&
	\mbox{}\hfill
	(b) 	\begin{minipage}[t]{4cm}
		There are only three admissible lozenges
		\end{minipage}
	\hfill\mbox{}
	\\

	\end{tabular}
	\caption{Triangles and lozenges \worklabel{triangles:lozenges}}
}

\fig{
	\begin{tabular}{p{5cm}@{\hskip1cm}p{5cm}}

	\mbox{}\hfill
	\epsfig{file=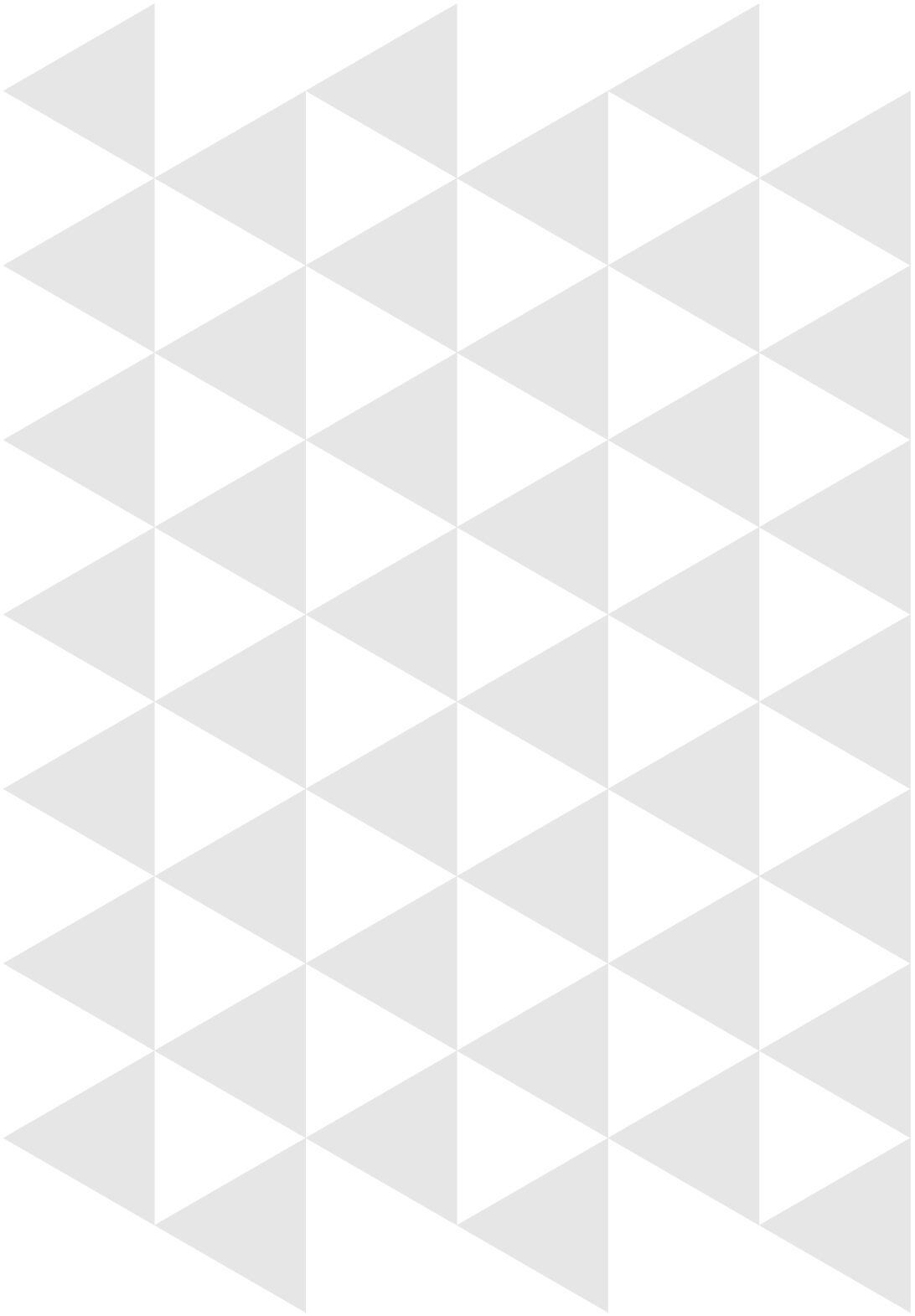, 
		width=2cm, clip=}
	\hfill\mbox{}
	&
	\mbox{}\hfill
	\epsfig{file=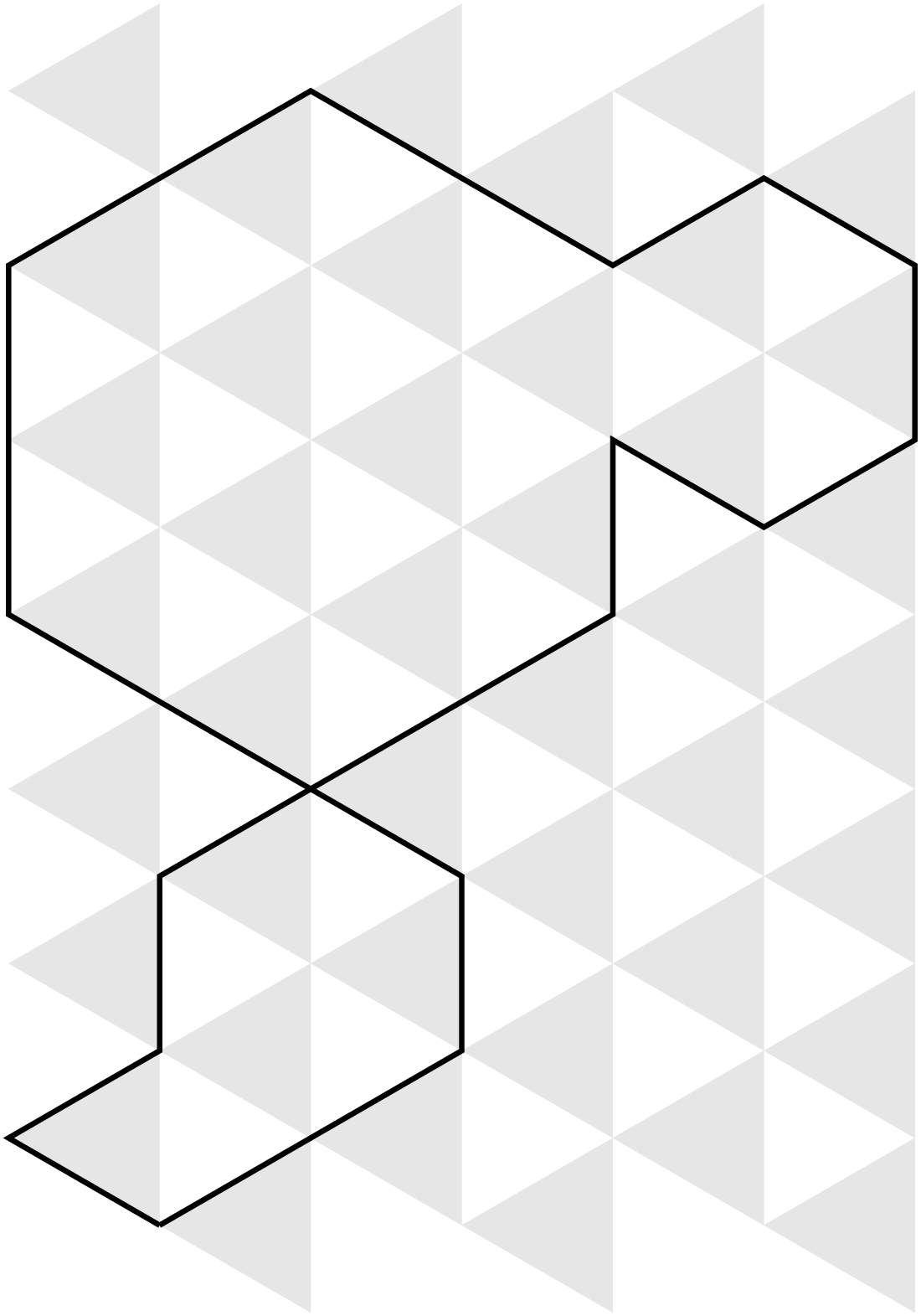, width=2cm, clip=}
	\begin{picture}(0,0)
	\put(-70,50){\vector(1,0){25}}
	\put(-81,47){$\cal D$}
	\put(17,45){\vector(-2,1){20}}
	\put(19,41){$\cal P$}
	\end{picture}
	\hfill\mbox{}
	\\[1mm]

	\mbox{}\hfill
	(a) 	\begin{minipage}[t]{4cm}
		The plane can be covered with regular triangles
		\end{minipage}
	\hfill\mbox{}
	&
	\mbox{}\hfill
	(b) 	\begin{minipage}[t]{4cm}
		A general domain $\cal D$ limited by a closed path $\cal P$
		\end{minipage}
	\hfill\mbox{}
	\\

	\end{tabular}
	\caption{A closed path in the triangular grid 
		\worklabel{triangular:grid}}
}

We now define three tiles by gluing together two regular triangles; let us
call them \emph{lozenges} (see Figure~\ref{triangles:lozenges}~(b)). A
\emph{tiling} of $\cal D$ is a set of lozenges included in $\cal D$ with
pairwise disjoint interiors such that the union of the lozenges is $\cal D$
itself. A domain is \emph{tileable} if it admits at least one tiling. Its
\emph{boundary} or \emph{contour} is the set of its edges that belong to
exactly one of its triangles.

Any finite-length closed path ${\cal P}$ whose edges belong to triangles
in the triangular grid can be viewed as the boundary path of a domain
${\cal D}$ (see Figure~\ref{triangular:grid}~(b)).

\fig{
	\begin{tabular}{p{4.5cm}@{\hskip1cm}p{4.5cm}}

	\mbox{}\hfill
	\epsfig{file=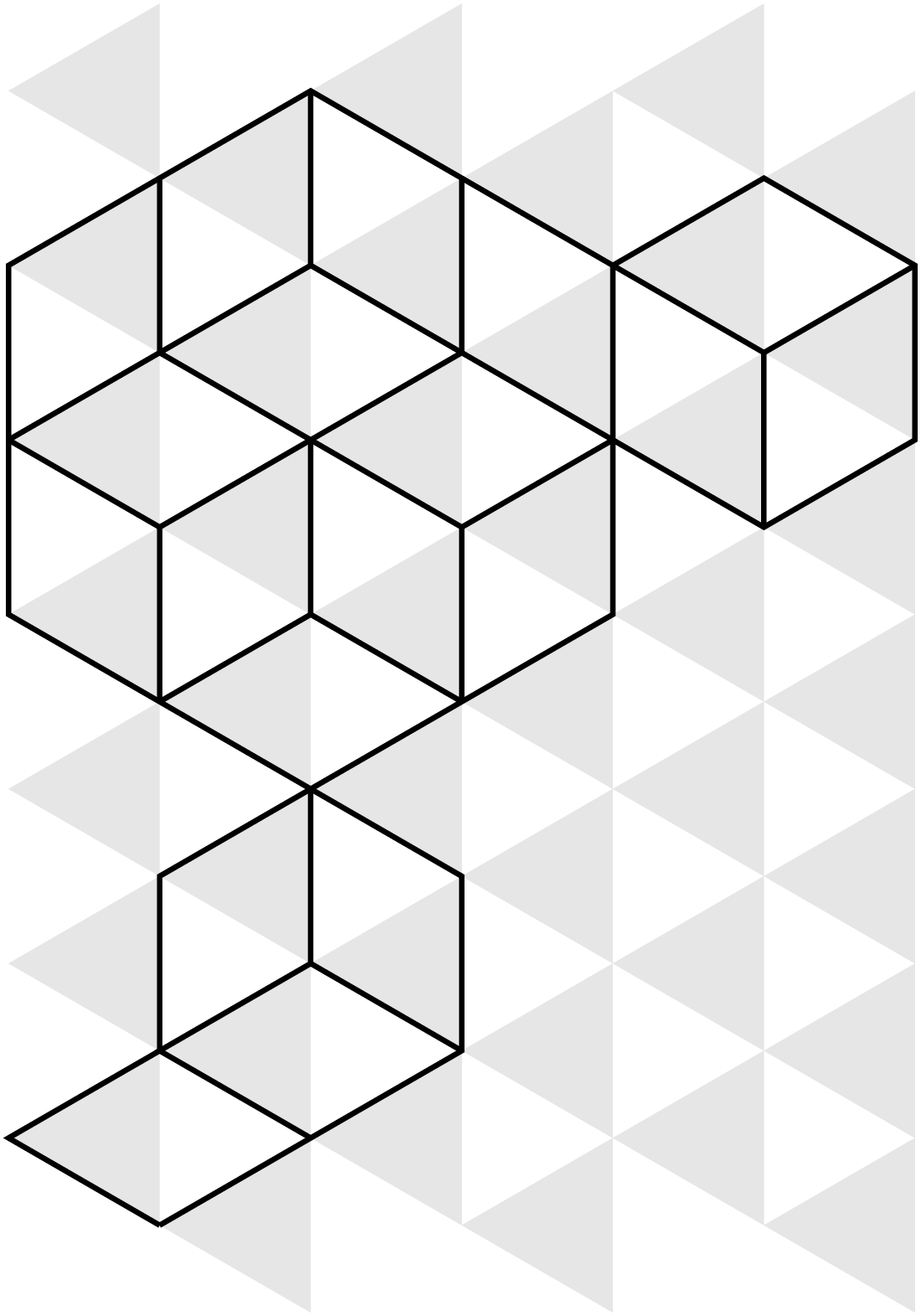, width=2cm, clip=}
	\hfill\mbox{}
	&
	\mbox{}\hfill
	\epsfig{file=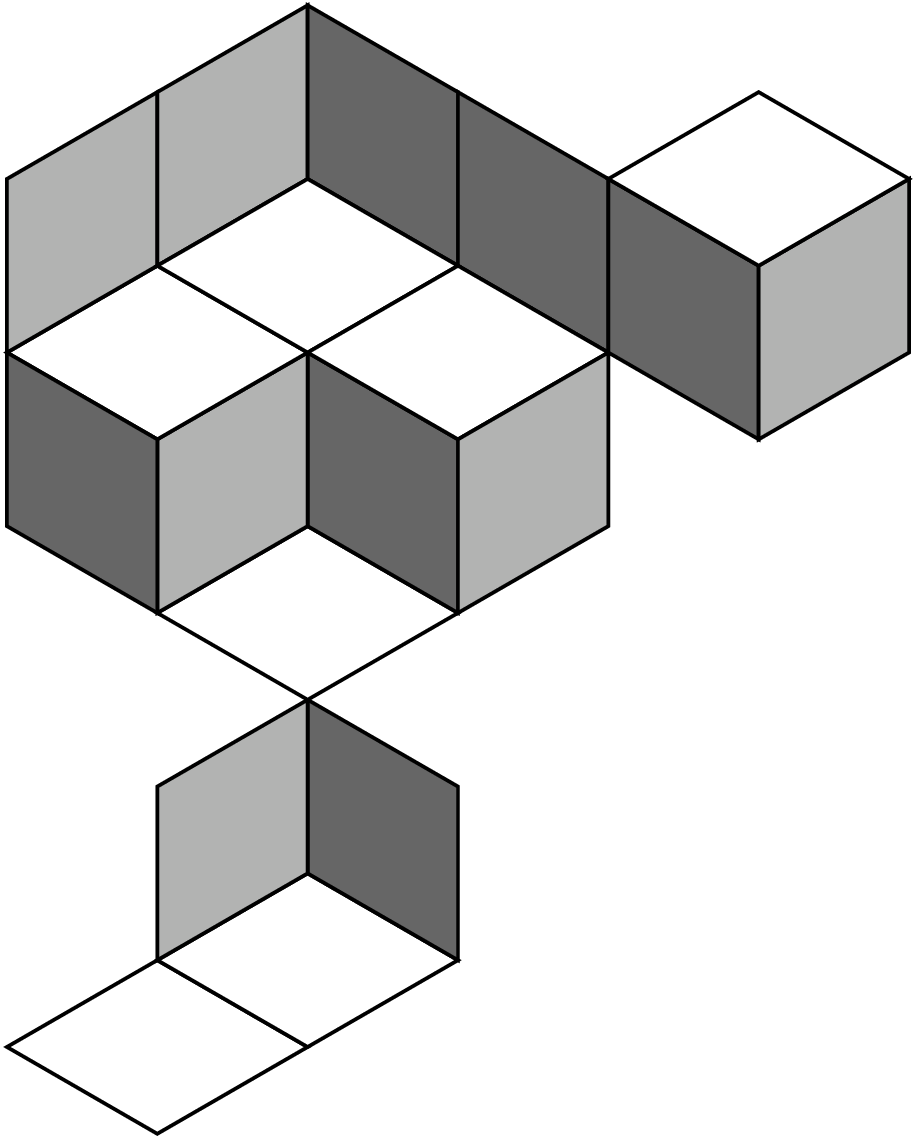, width=2cm, clip=}
	\hfill\mbox{}
	\\

	\mbox{}\hfill (a) \hfill\mbox{}
	&
	\mbox{}\hfill (b) \hfill\mbox{}
	\\

	\end{tabular}
	\caption{A tiling of $\cal D$\worklabel{first:tiling}}
}

%!!!!!!!!!!!!!!!!!!!!!!!!                       !!!!!!!!!!!!!!!!!!!!!!!!
%!!!!!!!!!!!!!!!!!!!!!!!!   Groupes de pavage   !!!!!!!!!!!!!!!!!!!!!!!!
%!!!!!!!!!!!!!!!!!!!!!!!!                       !!!!!!!!!!!!!!!!!!!!!!!!

\subsection{Tiling groups}

Figure~\ref{first:tiling}~(b) suggests a connection between tilings by
lozenges and certain piles of cubes, which we now attempt to clarify. Our tool
here is Conway and Lagarias' tiling groups, a formal description of which can
be found in \cite{Conway,Thurston}, so we will restrict ourselves to an
intuitive (but rigourous) approach. The idea is to simplify contour words so
that those of tileable domains are equivalent to the empty word. 

First, label the sides of triangles with letters; in our case, the set
$\{a,b,c\}$ is sufficient (see Figure~\ref{contour:words}~(a)). Our purpose is
to use words to partially describe tilings. Let $T$ be a tiling of a polygon
and let $P$ be a directed path (in the triangular grid) whose edges belong to
lozenges in~$T$; that is, $P$ never cuts the interior of a lozenge.  Such a
path will be called $T$-\emph{valid} or \emph{valid for~$T$}. Given a starting
point, $P$ is completely encoded by a word~$x$ on the alphabet
$\{a,b,c,a^{-1}, b^{-1},c^{-1}\}$.

We want $P$ to be rather straightforward, so that following an edge and then
following it in the opposite direction should not change~$x$. We thus impose
$aa^{-1}=bb^{-1}=cc^{-1}=\varepsilon$ (the empty word) and make all possible
simplifications in~$x$, so that~$P$ is now encoded by a word in the free group
$F(\{a,b,c\})$. Although this is not a pre-requisite for the rest of this
paper, the reader unfamiliar with free groups can look up page~257 in
\cite{Fraleigh} for instance.

If $P'$ is another $T$-valid path, encoded by a word~$y$ in $F(\{a,b,c\})$,
with the same starting and ending points as~$P$, then~$P$ and~$P'$ define a
polygon~$\cal D$ which is tiled by lozenges, and therefore tileable. To
reflect this fact, we set the contour word $x\,y^{-1}$ of $\cal D$ to be the
empty word,~$\varepsilon$. In particular, the contour words of each of our
three lozenges should be set to~$\varepsilon$, that is (see
Figure~\ref{contour:words}~(b))  $bcb^{-1}c^{-1}=\varepsilon$,
$aba^{-1}b^{-1}=\varepsilon$ and $cac^{-1}a^{-1}=\varepsilon$ by reading
counterclockwise (starting from a different vertex only changes the word by a
circular permutation). Note that this process is equivalent to removing loops
in~$P$, or stating that two $T$-valid paths having the same starting and
ending points should be described by the same word, as we will see in
Corollary~\ref{coro:same:word}.

\fig{
	\begin{tabular}{p{2cm}@{\hskip2cm}p{7.6cm}}

	\mbox{}\hfill
	\epsfig{file=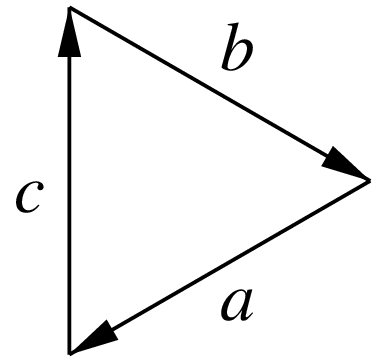, height=1.4cm}
	\hfill\mbox{}
	&
	\mbox{}\hfill
	\epsfig{file=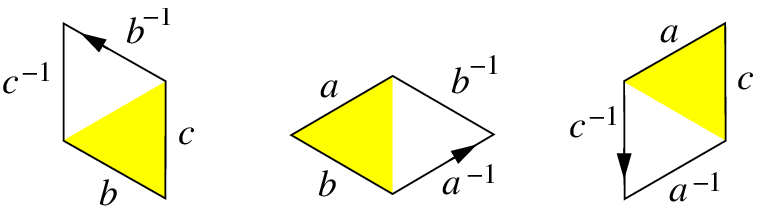, width=7cm, clip=}
	\hfill\mbox{}
	\\[1mm]

	\mbox{}\hfill (a) \hfill\mbox{}
	&
	\mbox{}\hfill (b) \hfill\mbox{}
	\\

	\end{tabular}
	\caption{Labelling the sides \worklabel{contour:words}}
}

The words we consider now belong to the group

$$L = \langle a,b,c \:|\: 	bcb^{-1}c^{-1} 
				= aba^{-1}b^{-1}
				= cac^{-1}a^{-1}
				= \varepsilon \rangle$$

\noindent which is Conway and Lagarias' \emph{lozenge group}. 

\begin{definition}
	Let $T$ be a tiling of a polygon and let $P=(v_0,v_1,\ldots,v_p)$
be a valid directed path in $T$. Each edge $(v_i,v_{i+1})$ can be labelled
by an element of the alphabet $\{a,b,c,a^{-1},b^{-1},c^{-1}\}$; the
\emph{label} of $P$ is the word $w$ obtained by the concatenation of the
labels of its edges in the order $(v_0,v_1)$,\dots,$(v_{p-1},v_p)$. The
\emph{free label} $\alpha(w)$ of $P$ is the representative element of $w$ in
the free group $F(\{a,b,c\})$. The \emph{$L$-label} $\ell(w)$ of $P$ is
the representative element of $w$ in the lozenge group $L$.

	For the sake of simplicity, the free label and the $L$-label of
$P$ will also be written $\alpha(P)$ and $\ell(P)$.
	\worklabel{def:labels}
	\end{definition}

The following lemma gives a necessary condition for a domain to be
tileable:

\begin{proposition}[\cite{Conway}]
	If a polygon is tileable, then the $L$-label of its boundary path 
is trivial.
	\worklabel{necessary:condition}
	\end{proposition}

\noindent \emph{Proof} \quad We make the proof by induction on the surface.
(The reader familiar with De Bruijn's worms will readily find a direct proof.)
It is enough to prove the result for elementary cycles. Besides, the result
holds for single lozenges.

Let $P$ be a valid path in a tiling $T$ of the polygon $\cal P$, such that
$P$ is distinct from the boundary path $B$ of $\cal P$. $P$ cuts $\cal P$
into two tileable polygons ${\cal P}_1$ and ${\cal P}_2$. There exist $w$
and $w'$ such that the boundary paths of ${\cal P}_1$ and ${\cal P}_2$ are
$w\alpha(P)$ and $\alpha(P)^{-1}w'$, both of which are trivial by the
induction hypothesis. Therefore $\ell(B) = \ell(w \cdot w') =
\ell(w\alpha(P) \cdot \alpha(P)^{-1}w') = \ell(w\alpha(P)) \cdot
\ell(\alpha(P)^{-1}w') = \varepsilon$.\mbox{}\hfill$\square$

\begin{corollary}
	Let $T$ be a tiling of a polygon. The $L$-label of any closed
valid path in $T$ is trivial.
	\end{corollary}

\noindent \textit{Proof} \quad Indeed, a valid path delimits a tileable
polygon of which it is the boundary path.\mbox{}\hfill$\square$

\begin{corollary}
	\worklabel{coro:same:word}
	Let $T$ be a tiling of a polygon. Two valid paths in $T$ having 
the same starting and ending points have the same $L$-label.
	\end{corollary}

\noindent \textit{Proof} \quad Indeed, if $w$ and $w'$ are the labels of these
paths, $w^{-1} \cdot w'$ is the label of a valid closed path in
$T$.\mbox{}\hfill$\square$

\medskip

The relations that appear in the definition of~$L$ can be rewritten as
$ab=ba$, $ac=ca$ and $bc=cb$, so that $L$ has three generators that
commute with each other and is therefore isomorphic to~$\mathbb{Z}^3$.
How can we interpret this nice result?

\smallskip

Let $T$ be a tiling of a polygon~$\cal P$ and let~$v$ be a vertex of~${\cal
P}$. A valid path in~$T$ can be associated with a unique word in~$L$, which in
turn corresponds to a unique path in the 1-skeleton of a cubical
tesselation~$\cal T$ of space. There is therefore a one-to-one correspondence
between vertices of $\cal D$ (resp. segments in~$T$) and vertices of
$\mathbb{Z}^3$ (resp. segments in~$\cal T$). Following an edge in~$T$ is
analogous to following one of $\mathbb{Z}^3$'s generators. This bijection
allows us to lift each edge of~$T$ to an edge in~$\cal T$, each lozenge to a
square in~$\mathbb{Z}^3$, so that~$T$ is equivalent to a collection of squares
in~$\cal T$. This squares may or may not define the visible parts of cubes,
depending on~$T$. The squares can be projected along the direction $(1,1,1)$
in $\mathbb{Z}^3$ to give back~$T$ (see Figures~\ref{projection}
and~\ref{first:tiling}~(b)).

\fig{
	\epsfig{file=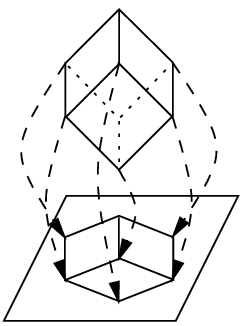, width=3cm}
	\caption{Projection of a cube onto the plane $x+y+z=0$
		\worklabel{projection}}
}

%!!!!!!!!!!!!!!!!!!!!!!                          !!!!!!!!!!!!!!!!!!!!!!
%!!!!!!!!!!!!!!!!!!!!!!   Fonctions de hauteur   !!!!!!!!!!!!!!!!!!!!!!
%!!!!!!!!!!!!!!!!!!!!!!                          !!!!!!!!!!!!!!!!!!!!!!

\subsection{Height functions}
\worklabel{subsection:height:functions}

The contour word of a polygon may well have a trivial image in the lozenge
group even though the polygon is \emph{not} tileable. In other words, the
condition stated in Proposition~\ref{necessary:condition} is not sufficient.
Consider for instance Figure~\ref{integer:values}~(a): the polygon is clearly
non-tileable, but its contour word $bccab^{-1}c^{-1}c^{-1}a^{-1}$ has a
trivial image in~$L$ because letters commute.  More generally, any closed path
in $\mathbb{Z}^3$ corresponds to a word in which the number of $a$'s (resp.
$b$, $c$) is equal to the number of $a^{-1}$'s (resp. $b^{-1}$, $c^{-1}$), so
that the image of this word in~$L$ is trivial; it would be quite surprising if
the projection onto the plane of any closed path in $\mathbb{Z}^3$ was
tileable. More complicated examples can be exhibited, for instance using
Fournier's obstructions (see~\cite{Fournier3}). Instead of examining each
case, let us develop from this simple example a general idea.

\fig{
	\begin{tabular}{p{5cm}@{\hskip1cm}p{4cm}}

	\mbox{}\hfill
	\epsfig{file=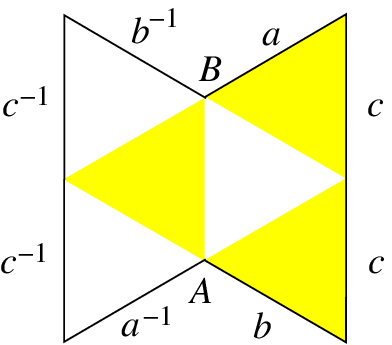, width=3.4cm}
	\hfill\mbox{}
	&
	\mbox{}\hfill
	\epsfig{file=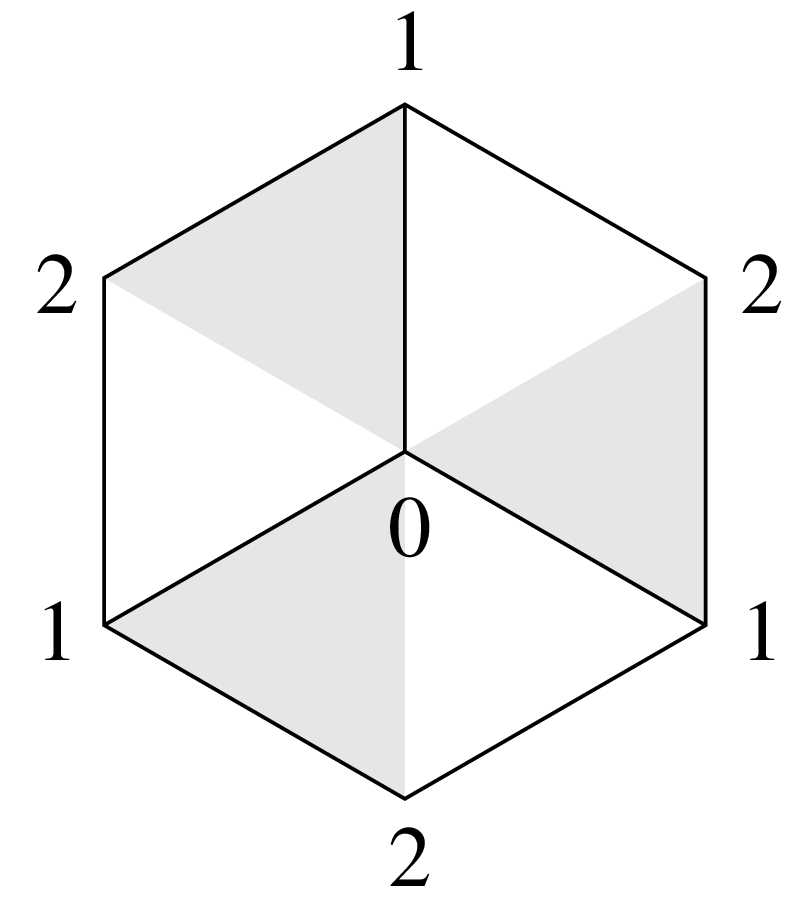, width=3cm}
	\hfill\mbox{}
	\\[1mm]

	\mbox{}\hfill
	(a) 	\begin{minipage}[t]{4cm}
		Lemma~\ref{necessary:condition} does not give a sufficient 
		condition
		\end{minipage}
	\hfill\mbox{}
	&
	\mbox{}\hfill
	(b) 	\begin{minipage}[t]{3.3cm}
		Height functions take integer values
		\end{minipage}
	\hfill\mbox{}
	\\

	\end{tabular}
	\caption{The need for height functions \worklabel{integer:values}}
}

Let $a$, $b$ and $c$ correspond to $(1,0,0)$, $(0,1,0)$ and $(0,0,1)$: we can
now see Figure~\ref{integer:values}~(a) as a closed path in $\mathbb{Z}^3$ (it
looks like a cyclohexane molecule in chair conformation). If the domain was
tileable, one could view a tiling as a collection of squares and cubes in
$\mathbb{Z}^3$. Since the distance between $A$ and $B$ is only one edge in the
triangular grid, these points would be linked by either the side or the
diagonal of a square in $\mathbb{Z}^3$, so that their distance in
$\mathbb{Z}^3$ would be~1 or~$\sqrt{2}$; but if~$A$ has coordinates $(0,0,0)$
then $B$ has coordinates $(1,1,2)$ since starting from~$A$ we follow $b$ once,
$c$ twice and $a$ once.

The distance (in $\mathbb{Z}^3$) between $A$ and $B$ can easily be made
greater: remark that the figure looks like a butterfly, and make the wings
bigger without moving $A$ or $B$. This adds 3 to the distance between $A$ and
$B$ at each step.

We see that the distance between two points seems important. Two points at
distance~1 in the triangular grid should not be distant ones in $\mathbb{Z}^3$
if the domain is tileable. What really matters in our butterfly example though
is the distance between $A$ and $B$ along the $(1,1,1)$ axis: following an
edge in a tiling is like following one of $\mathbb{Z}^3$'s generators, each of
which yields the same height increase, say, $i$, along the $(1,1,1)$ axis so
that every point in $\mathbb{Z}^3$, once orthogonally projected onto the
$(1,1,1)$ axis, is at a distance to the origin that is a multiple of~$i$.  To
each point $p$ of $\mathbb{Z}^3$ one can therefore give a value (multiple of
$i$) which corresponds to the distance between a fixed origin point and the
orthogonal projection of~$p$ onto the $(1,1,1)$ axis.  For convenience, we
will forget the geometrical interpretation, place the origin at an arbitrary
vertex and set~$i=1$ (see Figure~\ref{integer:values}~(b)).

%\fig{
%	\begin{tabular}{p{4cm}@{\hskip2cm}p{4cm}}
%
%	\mbox{}\hfill
%	\epsfig{file=dessins/helix.eps, width=2.5cm}
%	\hfill\mbox{}
%	&
%	\mbox{}\hfill
%	\epsfig{file=dessins/chemins_diriges.ps, width=1.5cm}
%	\hfill\mbox{}
%	\\[1mm]
%
%	\end{tabular}
%	\caption{An orientation of the triangular grid 
%		\worklabel{orientation}
%	}
%}

\medskip

$\!$We now proceed to properly define height functions, mainly following
Thurs\-ton (see~\cite{Thurston}) while preserving the algebraic point of view.

\begin{definition} 
	We call \emph{evaluation function} the morphism $\varphi$ from the
lozenge group $(L,\cdot)$ with generators $a$, $b$, $c$, to the group
$(\mathbb{Z},+)$ of integers such that $\varphi(a) = \varphi(b) =
\varphi(c) = 1$.
	\end{definition}

Note that this implies $\varphi(\varepsilon) = 0$ and $\varphi(w\cdot w') =
\varphi(w) + \varphi(w')$ for any $w$, $w'$ in $L$, whence $\varphi(a^{-1}) =
\varphi(b^{-1}) = \varphi(c^{-1}) = -1$.

\begin{definition}
	\worklabel{def:height:function}
	Let $T$ be a tiling of a polygon $\cal P$ and let $v$ be a vertex
on the boundary path of ${\cal P}$. The height function induced by $T$ and
$v$ is the function that maps each vertex $x$ of ${\cal P}$ to the image
by the evaluation function of any valid path from $v$ to $x$.
	\end{definition}

The correctness of this definition stems from Corollary~\ref{coro:same:word}.

\begin{lemma}
	On the boundary path of $\cal P$, the heights of the vertices do 
not depend on~$T$.

	\worklabel{lemma:boundary:height}
	\end{lemma}

\noindent\textit{Proof} \quad Indeed, if $x$ is such a vertex, there exists a
path from $v$ to $x$ that lies on the boundary path of $\cal P$, so the result
follow by induction.\mbox{}\hfill$\square$

\medskip

Note that only the vertices on the boundary path have fixed heights; those of
inner vertices \emph{do} depend on~$T$ (see for instance
Figure~\ref{fig:basic:tools}~(a)). 

Since changing the reference vertex $v$ only changes the height function by a
constant, we will often refer to a height function without mentionning~$v$.  
Thus we will write ``the height function associated with the tiling''.

%!!!!!!!!!!!!!!!!!!!!!                            !!!!!!!!!!!!!!!!!!!!!
%!!!!!!!!!!!!!!!!!!!!!   Algo de pavage maximal   !!!!!!!!!!!!!!!!!!!!!
%!!!!!!!!!!!!!!!!!!!!!                            !!!!!!!!!!!!!!!!!!!!!

\subsection{Thurston's algorithm}
\worklabel{Thurston:algo}

We have seen (see Section~\ref{subsection:height:functions}) that
Proposition~\ref{necessary:condition} does not give a sufficient condition for
the tileability of a polygon. Thurston's height functions provides a
constructive algorithm, outlined in \cite{Thurston}, to determine whether a
polygon can be tiled, and exhibit a tiling if it can be done.

We will build the minimal tiling of the polygon, in a sense that will be made
clear in Section~\ref{flips:lattices}. There is a natural (partial) order on
the height functions as the tiling changes and this is interpreted as an order
on the tilings. 

The following algorithm, defined by Thurston in~\cite{Thurston}, is based on
this simple result (see Proposition~\ref{prop:local:extremum} and
Definition~\ref{def:flip} in Section~\ref{flips:lattices} to know more about
flips):

\begin{lemma}
	Let $v$ be a vertex of a polygon $\cal P$ such that the height
function associated with the minimal tiling of~$\cal P$ is maximal on~$v$.  
This vertex cannot belong to the interior of $\cal P$, otherwise it could be
flipped down: therefore $v$ lies on the boundary path of $\cal P$.

\end{lemma}

\begin{algorithm}
	\worklabel{algo:Thurston}
\end{algorithm}

\begin{itemize}

\item \textbf{Input}: A polygon $\cal P$.

\item \textbf{Output}: The minimal tiling of $\cal P$ if the polygon is 
tileable, untileability otherwise.

\item \textbf{Initialization}: Initialize the list $L$ to~$\emptyset$.

\item \textbf{Step 1}: If $\cal P$ is a single point, return $L$.

\item \textbf{Step 2}: Compute the height function on the boundary path of
$\cal P$. If a vertex is given two different heights, return untileability.

\item \textbf{Step 3}: Let $v$ be a vertex of $\cal P$ of maximal height.  
There exist $u$ and $w$ on the boundary path $B$ of $\cal P$ such that $(u,v)$
and $(v,w)$ are segments of~$B$. Place a lozenge $\ell$ on $\cal P$ so that
$u$, $v$ and $w$ belong to it. Add $\ell$ and its position to $L$. Update
$\cal P$ to $\cal P \setminus \ell$. Go back to step 1.

\end{itemize}

\noindent For the proof of the algorithm, the reader can refer to
\cite{Thurston}. Its complexity is linear in the number of triangles in~$\cal
P$. From this algorithm one can easily deduce an algorithm to build the
maximal tiling.

%!!!!!!!!!!!!!!!!!!!!!!!!!!!!              !!!!!!!!!!!!!!!!!!!!!!!!!!!!
%!!!!!!!!!!!!!!!!!!!!!!!!!!!!   Treillis   !!!!!!!!!!!!!!!!!!!!!!!!!!!!
%!!!!!!!!!!!!!!!!!!!!!!!!!!!!              !!!!!!!!!!!!!!!!!!!!!!!!!!!!

\subsection{Lattices}

We will see in Section~\ref{lattice:structure} that the height functions
associated with the tilings of a polygon define a partial order on this
tilings; this order has the interesting property of being a lattice, which we
now define. See Figure~\ref{fig:lattice:example} for an example of a graphical
representation of a lattice.

\begin{definition}[Lattice]
	A set $S$ partially ordered by a relation $\preccurlyeq$ is a
\emph{lattice} if for any $a$ and $b$ in $S$ there exist $i$ and $s$ in $S$
such that

\begin{itemize}

\item $i \preccurlyeq a \preccurlyeq s$ and $i \preccurlyeq b 
\preccurlyeq s$;

\item if $x \preccurlyeq a$ and $x \preccurlyeq b$ then $x \preccurlyeq i$ 
for any $x \in S$;

\item if $a \preccurlyeq y$ and $b \preccurlyeq y$ then $s \preccurlyeq y$ 
for any $y \in S$;

\end{itemize}

\noindent Such elements $i$ and $s$ are called the \emph{infimum} and the
\emph{supremum} of $a$ and $b$, noted $\text{\rm inf}(a,b)$ and $\text{\rm
sup}(a,b)$ (or $a \wedge b$ and $a \vee b$ for short). If $a$ and $b$ happen
to be comparable (\textit{e.g.} $a \preccurlyeq b$) then $i$ and $s$ are
called the \emph{minimum} and the \emph{maximum} of $a$ and $b$, noted
$\text{\rm min}(a,b)$ and $\text{\rm max}(a,b)$ (equal respectively to $a$ and
$b$ if $a \preccurlyeq b$).

	\worklabel{lattice:definition}
	\end{definition}

If $S$ is a finite lattice (such as the tiling lattices considered in this
paper), then it admits one maximal and one minimal element.

A lattice $(S,\preccurlyeq)$ is \emph{distributive} if for any $a$, $b$ and
$c$ in $S$ one has

\[
\left\{\hskip-1.5mm
\begin{array}{r@{\ =\ }l}

a \wedge (b \vee c) & (a \wedge b) \vee (a \wedge c) \\
a \vee (b \wedge c) & (a \vee b) \wedge (a \vee c). \\

\end{array}
\right.
\]

\begin{definition}[Interval]
	An \emph{interval} $[a;b]$ (with $a \preccurlyeq b$) in a lattice
$({\cal L},\preccurlyeq)$ is the set of all $x\in\cal L$ such that $a
\preccurlyeq x \preccurlyeq b$.

	\worklabel{def:interval}
	\end{definition}

Note that $a$ and $b$ need be comparable elements.

\begin{proposition}
	Let $({\cal L},\preccurlyeq)$ be a lattice and $\cal I$ an interval in
$\cal L$. The order $\preccurlyeq$ defines on $\cal I$ a structure of lattice.

	\worklabel{prop:interval:sub:lattice}
	\end{proposition}

%!!!!!!!!!!!!!!!!!!!!!!                           !!!!!!!!!!!!!!!!!!!!!!
%!!!!!!!!!!!!!!!!!!!!!!   Ferrers et partitions   !!!!!!!!!!!!!!!!!!!!!!
%!!!!!!!!!!!!!!!!!!!!!!                           !!!!!!!!!!!!!!!!!!!!!!

\subsection{Partitions, Ferrers diagrams and plane partitions}

These definitions will be needed in Section~\ref{section:algo:hexagons}. See 
also Figure~\ref{fig:Ferrers:pp}.

\begin{definition}
	A \emph{partition} of an integer $n$ is a non-increasing list of
positive integers with sum equal to $n$.

	\worklabel{def:partition}
	\end{definition}

\begin{definition}
	A \emph{Ferrers diagram} is a geometrical representation of a
partition: given a partition $(\lambda_1, \ldots , \lambda_p)$, the
corresponding Ferrers diagram is a collection of $p$ consecutive
left-justified rows respectively containing $\lambda_1$, \ldots, $\lambda_p$
consecutive squares.

	\worklabel{def:Ferrers:diagram}
	\end{definition}

\fig{
	\begin{tabular}{p{6.5cm}@{\hskip1cm}p{3.5cm}}

	\mbox{}\hfill
	\raisebox{4mm}{
		\epsfig{file=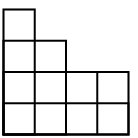, width=1.3cm}
	}
	\hfill\mbox{}
	&
	\mbox{}\hfill
	\epsfig{file=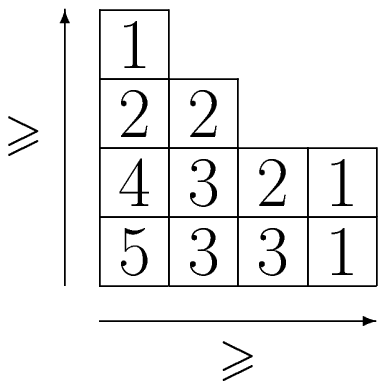, width=1.8cm, clip=}
	\hfill\mbox{}
	\\[1mm]

	\mbox{}\hfill
	(a) 	\begin{minipage}[t]{4.5cm}
		The Ferrers diagram of the partition $(4,4,2,1)$
		\end{minipage}
	\hfill\mbox{}
	&
	\mbox{}\hfill
	(b) 	\begin{minipage}[t]{2.7cm}
		A plane partition of~27.
		\end{minipage}
	\hfill\mbox{}
	\\

	\end{tabular}
	\caption{Ferrers diagrams and plane partitions 
		\worklabel{fig:Ferrers:pp}
	}
}

\begin{definition}
	A \emph{plane partition} of an integer $n$ is a filling of a Ferrers
diagram with integers (of sum $n$), such that the values along the rows and
the columns is non-increasing.

	\worklabel{def:plane:partition}
	\end{definition}

%%%%%%%%%%%%%%%%%%%%%%%%%%%%%%%%%%%%%%%%%%%%%%%%%%%%%%%%%%%%%%%%%%%%%%%%
%%%%%%%%%%%%%%%%%%%%%%%%%%%%%%%%%%%%%%%%%%%%%%%%%%%%%%%%%%%%%%%%%%%%%%%%
%%%%%%%%%%%%%%%%%%%%%%%                           %%%%%%%%%%%%%%%%%%%%%%
%%%%%%%%%%%%%%%%%%%%%%%   Structure de treillis   %%%%%%%%%%%%%%%%%%%%%%
%%%%%%%%%%%%%%%%%%%%%%%                           %%%%%%%%%%%%%%%%%%%%%%
%%%%%%%%%%%%%%%%%%%%%%%%%%%%%%%%%%%%%%%%%%%%%%%%%%%%%%%%%%%%%%%%%%%%%%%%
%%%%%%%%%%%%%%%%%%%%%%%%%%%%%%%%%%%%%%%%%%%%%%%%%%%%%%%%%%%%%%%%%%%%%%%%

\section{Flips and lattices}
\worklabel{flips:lattices}

In this section, we first show that there is a one-to-one correspondence
between tilings and height functions (up to constant) for a given polygon
$\cal P$, which allows us to show that the set of the tilings of $\cal P$
under the partial order defined by the height functions has the structure of a
lattice. We then consider a local operation (called a flip) on a tiling and
prove that it has a straightforward translation in terms of the underlying
height function.

\fig{
	\begin{tabular}{p{4.5cm}@{\hskip1.5cm}p{5.5cm}}

	\mbox{}\hfill
	\epsfig{file=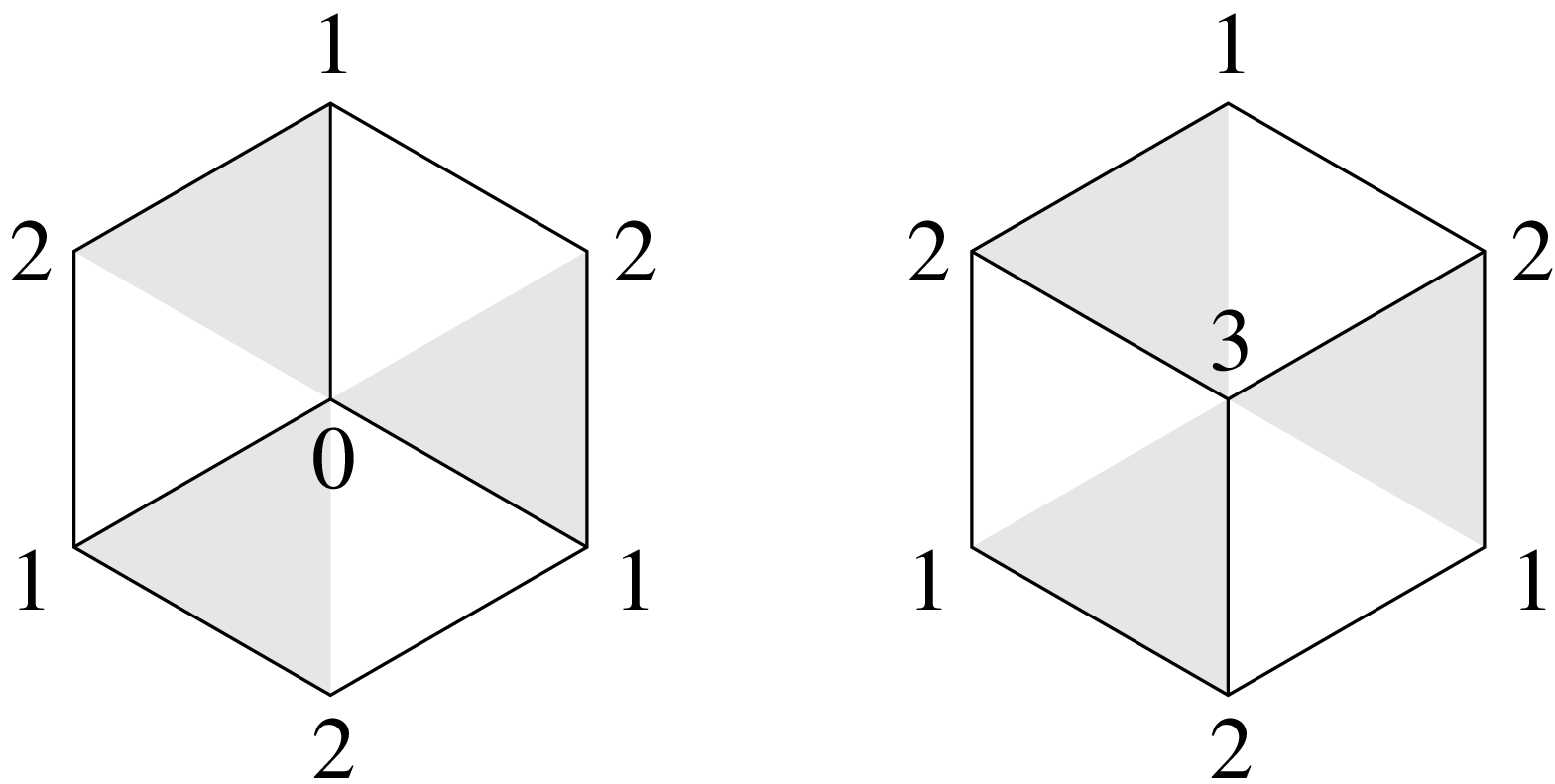, width=4cm}
	\hfill\mbox{}
	&
	\mbox{}\hfill
	\epsfig{file=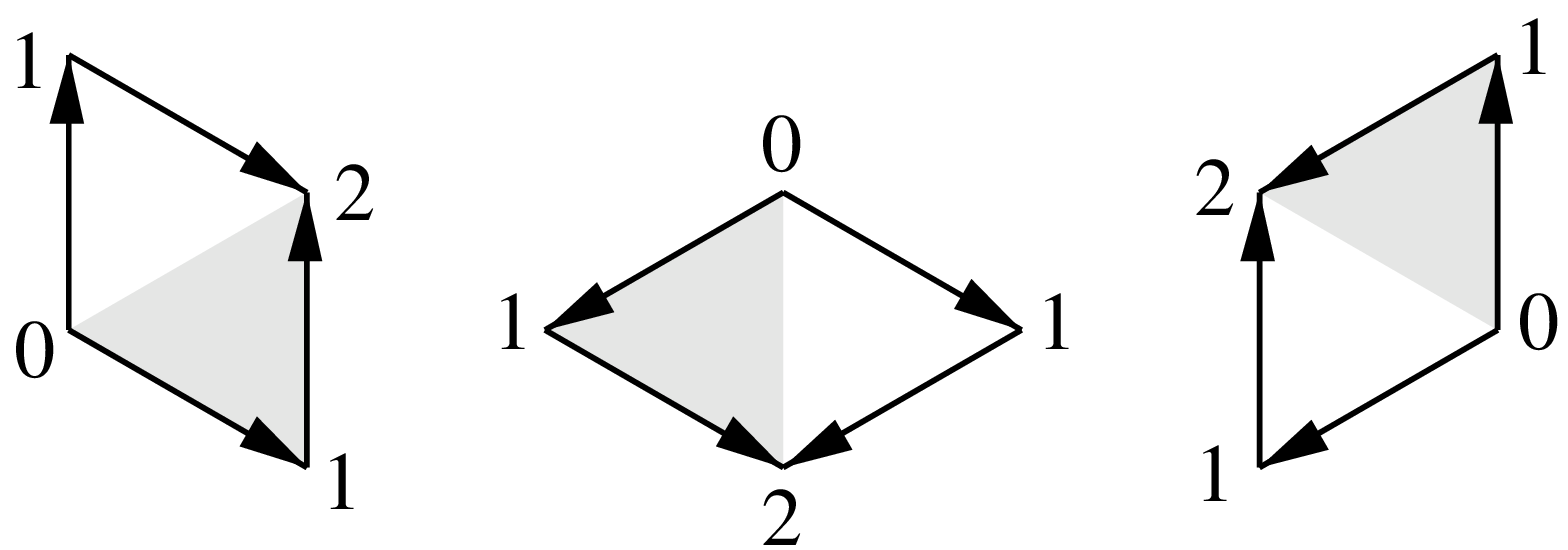, width=5cm}
	\hfill\mbox{}
	\\[1mm]

	\mbox{}\hfill
	(a) 	\begin{minipage}[t]{2.9cm}
		A hexagon can be tiled in two ways
		\end{minipage}
	\hfill\mbox{}
	&
	\mbox{}\hfill
	(b) 	\begin{minipage}[t]{4cm}
		Height functions for the three lozenges
		\end{minipage}
	\hfill\mbox{}
	\\

	\end{tabular}
	\caption{Our building blocks \worklabel{fig:basic:tools}}
}

A basic remark will be the leading idea for most of this section: a hexagon
can be tiled in exactly two ways with lozenges (see
Figure~\ref{fig:basic:tools}~(a)).  Switching between these two tilings is
called \emph{flipping}.  This operation allows us to deduce new tilings from a
known one. The really challenging idea is that this allows us in fact to
navigate between \emph{all} the tilings of~$\cal D$. All the results and
proofs in this section are adapted to the case of lozenges from \cite{Remila},
which deals with dominoes.

%!!!!!!!!!!!!!!!!!                                     !!!!!!!!!!!!!!!!!
%!!!!!!!!!!!!!!!!!   Fonctions de hauteur et pavages   !!!!!!!!!!!!!!!!!
%!!!!!!!!!!!!!!!!!                                     !!!!!!!!!!!!!!!!!

\subsection{From height functions to tilings}
\worklabel{height:functions}

We investigate in this section the link between height functions and tilings.
An edge will be said \emph{positively directed} if it is labelled by $a$, $b$
or~$c$. Recall that an edge in the triangular grid is valid for a tiling~$T$
if and only if it does not cut the interior of a lozenge in~$T$.

We have seen (see Definition~\ref{def:height:function}) how a height function
$h$ can be defined, starting from a particular tiling $T$ of a polygon $\cal
P$. Along any valid positively directed edge the height changes by~1; along
any invalid positively directed edge, the height changes by~$-2$ (see
Figures~\ref{contour:words}~(b) and~\ref{fig:basic:tools}~(b)). Therefore $h$
is merely an encoding of $T$, and the latter can be reconstructed from the
former.

\begin{lemma}
	There is a one-to-one correspondence between the tilings of a 
polygon~$\cal P$ and the associated height functions (up to constant). In 
other words, if $T$ and $T'$ are two tilings of $\cal P$ and if $h_T(v) = 
h_{T'}(v)$ for any $v \in \cal P$ then $T = T'$.

	\worklabel{lemma:height:tiling}
	\end{lemma}

\noindent\textit{Proof}\quad Assume $\cal P$ and a height function $h_T$
given, the latter corresponding to an unknown tiling $T$ of $\cal P$. To
rebuild $T$ from $h_T$, it suffices to draw all the triangles inside $\cal P$
and to erase those edges whose endpoints have a height difference of~2 in
absolute value. These edges lie on the interior of lozenges and completely
characterize them so we are done.\mbox{}\hfill$\square$

\begin{proposition}
	Let $T$ and $T'$ be two tilings of a polygon $\cal P$. For each vertex
$v$ of $\cal P$, $h_T(v) - h_{T'}(v)$ is a multiple of 3. 
	\worklabel{prop:multiple:3}
	\end{proposition}

\noindent \textit{Proof} \quad We make the proof by induction on the vertices.  
The proposition is true on the boundary path of $\cal P$ since all height
functions take the same values there (see Lemma~\ref{lemma:boundary:height}).
Assume that the proposition holds for $v$ and let $v'$ be any of its
neighbours in the triangular grid. Let $\cal T$ be either $T$ or $T'$. If
$(v,v')$ is a valid positively directed edge in $\cal T$ then $h_{\cal T}(v')
- h_{\cal T}(v) = 1$; otherwise, $h_{\cal T}(v') - h_{\cal T}(v) = -2$ and
thus $h_{\cal T}(v') - h_{\cal T}(v)$ takes the same value modulo~3 whether
$(v,v')$ is valid or not. Consequently, if $h_T(v) - h_{T'}(v)$ is a multiple
of 3, then so is $h_T(v') - h_{T'}(v')$.~\mbox{}\hfill$\square$

\medskip

An intuitive way to consider this proposition is to think in terms of cubes:
adding~3 to the value of a height function on a vertex $v$ is equivalent to
adding exactly one cube (see Figure~\ref{fig:basic:tools}~(a)).

%!!!!!!!!!!!!!!!!!!!!!!                           !!!!!!!!!!!!!!!!!!!!!!
%!!!!!!!!!!!!!!!!!!!!!!   Structure de treillis   !!!!!!!!!!!!!!!!!!!!!!
%!!!!!!!!!!!!!!!!!!!!!!                           !!!!!!!!!!!!!!!!!!!!!!

\subsection{The lattice structure}
\worklabel{lattice:structure}

In this section, we prove that the set of the tilings of a polygon can be
endowed with a lattice structure. To this end, we define a partial order
between tilings by using their height functions: 

\begin{definition}
	Let $h_{T}$ and $h_{T'}$ be two height functions associated with the
tilings~$T$ and~$T'$ of a polygon~$\cal P$ such that $h_T$ and $h_{T'}$ take
the same values on the boundary path of~$\cal P$. $h_T$ is \emph{less than}
$h_{T'}$ (and we note $h_{T} \leqslant h_{T'}$) if $h_{T}(v) \leqslant
h_{T'}(v)$ for any vertex $v \in \cal P$.

	The order between height functions allow us to endow the set of the 
tilings of~$\cal P$ with a natural order: $T \preccurlyeq T'$ if and only if 
$h_T \leqslant h_{T'}$.

	\worklabel{def:order:tilings}
	\end{definition}

\begin{proposition}
	Let $T$ and $T'$ be two tilings of polygon $P$ and let $h_T$ and
$h_{T'}$ be their height functions. The functions $h_{\text{\rm min}} =
\text{\rm min}\:(h_{T},h_{T'})$ and $h_{\text{\rm max}} = \text{\rm
max}\:(h_{T},h_{T'})$ are themselves height functions.

\end{proposition}

\noindent \textit{Proof} \quad We prove the result only for $h_{\text{min}}$
since the other case is analogous. We make the proof by induction. Since $h_T$
and $h_{T'}$ take the same values on the boundary path $\cal B$ of $\cal P$
(see Lemma~\ref{lemma:boundary:height}), $h_{\text{min}}(v) = h_T(v) =
h_{T'}(v)$ for every $v \in \cal B$. Let $(v,v')$ be any positively directed
edge between vertices of $\cal P$ in the triangular grid. We claim that
$h_{\text{min}}(v') - h_{\text{min}}(v)$ equals either~1 or~$-2$.

\begin{itemize}

\item If $h_T(v)=h_{T'}(v)$ then $h_{\text{min}}(v') - h_{\text{min}}(v)$
equals either $h_{T}(v')-h_{T}(v)$ or $h_{T'}(v')-h_{T'}(v)$ and must
therefore be equal to either 1 or $-2$.

\item We can now assume without loss of generality that $h_T(v) < h_{T'}(v)$,
so that $h_{\text{min}}(v) = h_T(v)$. $h_T(v')$ can only be $h_T(v)+1$ or
$h_T(v) - 2$, and $h_{T'}(v)$ is at least $h_T(v) + 3$ by
Proposition~\ref{prop:multiple:3}. Moreover $h_{T'}(v')$ can only be
$h_{T'}(v) + 1$ or $h_{T'}(v)-2$, and therefore $h_{T'}(v')$ is at least
$h_T(v)+1$. Thus $h_{\text{min}}(v') = h_T(v')$, from which we derive
$h_{\text{min}}(v') - h_{\text{min}}(v) = h_{T}(v')-h_{T}(v)$ which must be
either 1 or $-2$.

\end{itemize}

We have shown that $h_{\text{min}}$ increases by either 1 or $-2$ along any
positively directed edge in the triangular grid.  Consequently, the set of
heights (for $h_{\text{min}}$) modulo~3 of the vertices of any triangle in the
triangular grid must exactly be $\{0,1,2\}$. Moreover the height difference
along an edge does not depend on which of the two neighbouring triangles was
chosen. Erase all edges whose endpoints have a difference of heights of $-2$.
What is left is a tiling of $\cal D$ with lozenges; that is, $h_{\text{min}}$
is a height function.~\mbox{}\hfill$\square$

\begin{corollary}
	The order $\preccurlyeq$ induces a structure of distributive lattice
on the set of the tilings of $\cal D$.
	\worklabel{coro:lattice:structure}
	\end{corollary}

\noindent \textit{Proof} \quad If $T_1$ and $T_2$ are two tilings of a domain
${\cal D}$ then the height functions $h_{\text{\rm min}}(T_1,T_2)$ and
$h_{\text{\rm max}}(T_1,T_2)$ are clearly the infimum and supremum of the
height functions $h_{T_1}$ and $h_{T_2}$ (see
Definition~\ref{lattice:definition}), and since height functions encode
tilings (see Proposition \ref{prop:multiple:3}) we have defined the infimum
and supremum of $T_1$ and~$T_2$.

To prove distributivity, it suffices to consider each vertex of $\cal D$,
which brings us back to checking the relation for integers and the usual min
and max functions.~\mbox{}\hfill$\square$

\medskip

\begin{figure}
\begin{center}
\epsfig{file=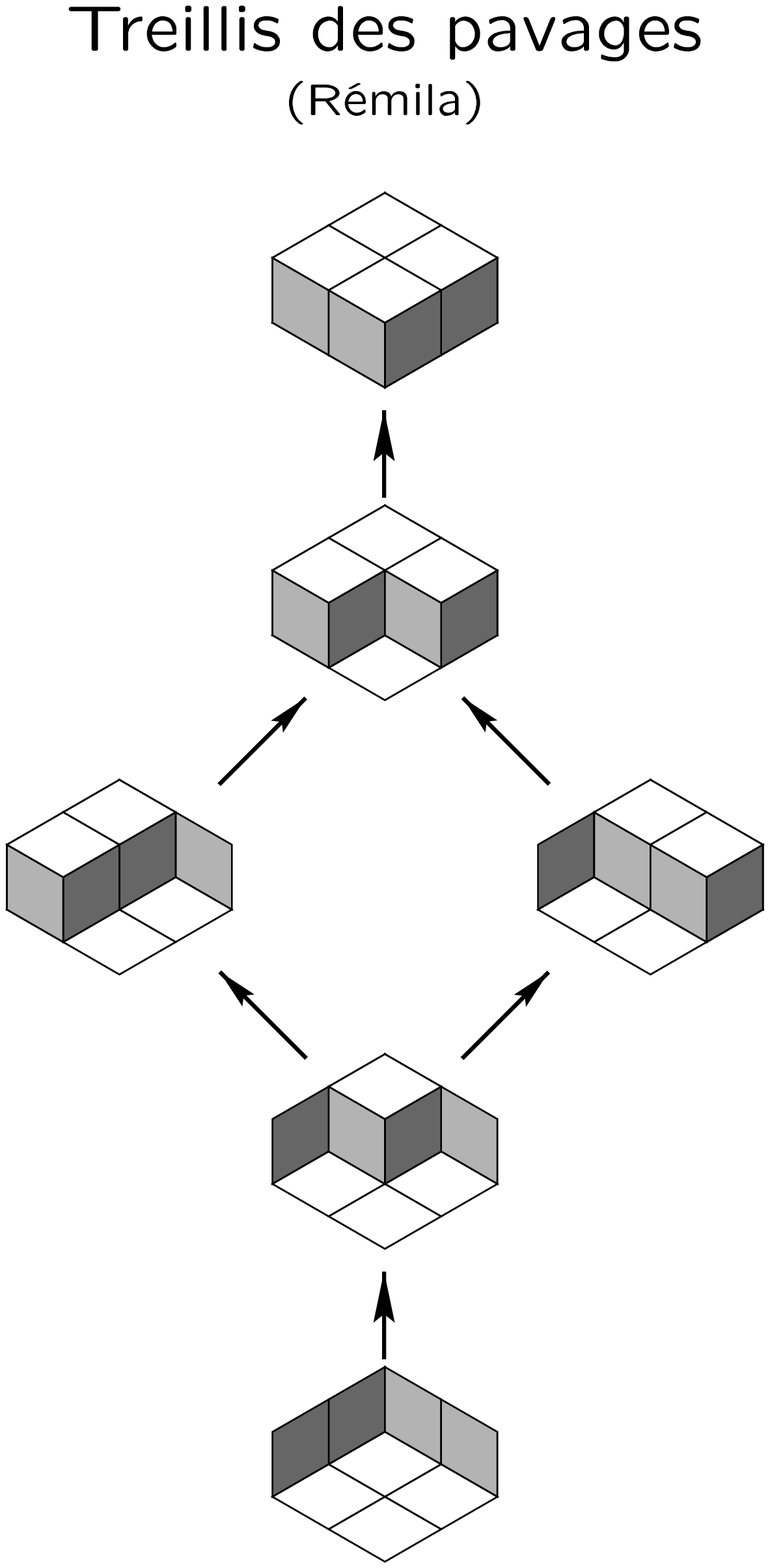, width=2cm, clip}
\caption{An example of lattice
	\worklabel{fig:lattice:example}
	}
\end{center}
\end{figure}

Figure~\ref{fig:lattice:example} provides an example of the graphical
representation of a lattice of tilings. Intuitively, one obtains the infimum
of two tilings $T_1$ and $T_2$ by selecting only the cubes which appear in
both, and their supremum by selecting the cubes which appear in either one. In
other words, $\text{\rm inf}(T_1,T_2)$ is encoded by $h_{\text{\rm
min}}(T_1,T_2)$ and $\text{\rm sup}(T_1,T_2)$ is encoded by $h_{\text{\rm
max}}(T_1,T_2)$.

%!!!!!!!!!!!!!!!!!!!!!!!!!!!!!!           !!!!!!!!!!!!!!!!!!!!!!!!!!!!!!
%!!!!!!!!!!!!!!!!!!!!!!!!!!!!!!   Flips   !!!!!!!!!!!!!!!!!!!!!!!!!!!!!!
%!!!!!!!!!!!!!!!!!!!!!!!!!!!!!!           !!!!!!!!!!!!!!!!!!!!!!!!!!!!!!

\subsection{Flips}
\worklabel{Flips}

In this section, we introduce an elementary operation classically called a
flip; we prove that the set of the tilings of a polygon is connected by flips.

\begin{definition}
	Let $T$ be a tiling of a polygon~$\cal P$. A \emph{local maximum}
(resp. \emph{minimum}\/) of the height function $h_T$ is a vertex $v$ in the
interior of $\cal P$ such that $h_T(v) \geqslant h_T(v')$ (resp. $h_T(v)  
\leqslant h_T(v')$) for any $v'$ neighbour of $v$.

	\worklabel{def:local:max:min}
	\end{definition}

\begin{proposition}
	Let $T$ be a tiling of a polygon $P$. A vertex $v$ in the interior
of~$\cal P$ is a local extremum of the height function $h_T$ if and only if it
is the center of a hexagon tiled with three lozenges.

	\worklabel{prop:local:extremum}
	\end{proposition}

%\begin{figure}
%\begin{center}
%\epsfig{file=dessins/hauteurs_des_hexagones.ps, width=9cm}
%\caption{Extrema of a height function\worklabel{up:down:flips}}
%\end{center}
%\end{figure}

\noindent \textit{Proof} \quad Let $v'$ and $v''$ be two neighbours of~$v$ so
that $(v, v', v'')$ is a triangle. If $(v',v)$ and $(v,v'')$ were both valid
edges in~$T$, then $h(v)$ would be less than $h(v')$ and more than $h(v'')$
(or the converse), so it would not be an extremum. Moreover $(v',v)$ and
$(v,v'')$ cannot both be invalid edges, therefore exactly one of them is.
Consequently, the hexagon around $v$ is tiled with exactly three lozenges.

The converse part of the proof is obvious (see
Figure~\ref{fig:basic:tools}~(a)).~\mbox{}\hfill$\square$

\begin{definition}[Flip]
	A \emph{flip} is the operation by which one switches from one tiling of 
a hexagon to the other. An \emph{up-flip} increases the height while a 
\emph{down-flip} decreases~it.
	\worklabel{def:flip}
	\end{definition}

We now prove that flips allow us to reach any tiling of~$\cal P$ from any
other tiling of~$\cal P$. To this end, we need one more definition:

\begin{definition}

If $T$ and $T'$ are two tilings of a polygon $\cal P$, then the distance 
between them is

\[ \Delta(T,T') = \sum_{v\in\cal D} |h_T(v) - h_{T'}(v)| \]

	\worklabel{def:distance:tilings}
	\end{definition}

\noindent This function is indeed a distance since it is symmetrical,
satisfies separation ($\Delta(T,T') = 0$ if and only if $T=T'$ by
Proposition~\ref{prop:multiple:3}) and the triangular inequality (since
$|\cdot|$ does).

\begin{proposition}
	Let $T$ and $T'$ be two tilings of a polygon $\cal P$. If $T
\preccurlyeq T'$ then there exists a sequence $(T_0=T , T_1, \ldots , T_n=T')$
of tilings of $\cal D$ such that $T_{p+1}$ is deduced from $T_p$ by a single
up-flip, $0 \leqslant p \leqslant n-1$.

	\worklabel{prop:atomic}
	\end{proposition}

\noindent \textit{Proof} \quad Assume $T \prec T'$, let $m = \text{min}
\{h_T(v) \:|\: v \in {\cal P} \text{ and } h_T(v) < h_{T'}(v)\}$ and let $v$
be a vertex of $\cal P$ such that $h_T(v) = m$. The vertex $v$ cannot belong
to the boundary path of $\cal P$ by Lemma~\ref{lemma:boundary:height} and
$h_{T'}(v)$ must be at least $h_T(v) + 3$ by
Proposition~\ref{prop:multiple:3}. Let $v'$ be a neighbour of~$v$ such that
the edge $(v,v')$ is valid for $h_T$. $h_{T'}(v')$ is at least $h_{T'}(v)-2
\geqslant h_T(v)+1 > h_T(v)$. Because $m$ is minimal, the edge $(v,v')$ must
be directed from $v$ to $v'$. Since this is true for any edge anchored in $v$
that is valid in $h_T$, we conclude that $v$ is a local minimum of $h_T$.

By an up-flip on $v$, we send $T$ to a tiling $T_1$ such that $T \prec T_1
\preccurlyeq T'$. Besides, $h_{T_1}$ differs from $h_T$ only on the
vertex~$v$, so that $\Delta (T,T_1) = 3$. Using induction, we thus build a
sequence of ever greater tilings. Since the distance between two tilings must
be a multiple of~3 (see Proposition~\ref{prop:multiple:3}), there comes a
tiling $T_n$ such that $\Delta (T_n,T')=0$ and $T_n \preccurlyeq T'$, and
since height functions encode tilings (see Lemma~\ref{lemma:height:tiling}),
$T_n=T'$.~\hfill$\square$

\begin{theorem}
	The set of the tilings of~$\cal P$ is connex: any tiling~$T$ of~$\cal
P$ can be reached from any other tiling~$T'$ of~$\cal P$ by using only flips.

	\worklabel{th:connex}
	\end{theorem}

\noindent \textit{Proof} \quad Indeed, any tiling of~$\cal P$ can be reached
from the minimal one by using only up-flips (see
Proposition~\ref{prop:atomic}) and the minimal tiling of~$\cal P$ can be
reached from any other tiling by using only down-flips.~\mbox{}\hfill$\square$

%%%%%%%%%%%%%%%%%%%%%%%%%%%%%%%%%%%%%%%%%%%%%%%%%%%%%%%%%%%%%%%%%%%%%%%
%%%%%%%%%%%%%%%%%%%%%%%%%%%%%%%%%%%%%%%%%%%%%%%%%%%%%%%%%%%%%%%%%%%%%%%
%%%%%%%%%%%%%%%%%%%%%%%                          %%%%%%%%%%%%%%%%%%%%%%
%%%%%%%%%%%%%%%%%%%%%%%   Les pseudo-hexagones   %%%%%%%%%%%%%%%%%%%%%%
%%%%%%%%%%%%%%%%%%%%%%%                          %%%%%%%%%%%%%%%%%%%%%%
%%%%%%%%%%%%%%%%%%%%%%%%%%%%%%%%%%%%%%%%%%%%%%%%%%%%%%%%%%%%%%%%%%%%%%%
%%%%%%%%%%%%%%%%%%%%%%%%%%%%%%%%%%%%%%%%%%%%%%%%%%%%%%%%%%%%%%%%%%%%%%%

\section{Hexagons and pseudo-hexagons}
\worklabel{hexagons:ph}

The ultimate goal of this paper is to provide an algorithm for uniquely
generating all the elements of the lattice of the tilings of a polygon~$P$. As
we will see in Section~\ref{domains:fracture:lines:seeds}, some hexagonal-like
subsets of~$P$ play a important role. We thus start with the rather simple
case of $\cal P$ being a hexagon or a pseudohexagon (see definitions below).

\subsection{Piles of cubes}

It is well known in tilings lore that the tilings of a hexagon of side~$n$ are
bijectively related to some piles of cubes. Using Conway and Lagarias' lozenge
group, we have explicited in Section~\ref{subsection:height:functions}
(following~\cite{Thurston}) the algebraic translation of the geometrical
intuition. The bijection between the lozenge group and $\mathbb{Z}^3$ implies
that it is equivalent to consider lozenges in a tiling or 2-cells (squares) in
the cubical tesselation of~$\mathbb{Z}^3$. Until now, we have been interested
in manipulating tilings and watching the translation in the figure. For
instance, rearranging the three tiles of a hexagon of side~1 is equivalent to
raising the height of exactly one point by~3: this is an up-flip. Now we look
at the bijection the other way: we manipulate squares in $\mathbb{Z}^3$ to
produce new informations on tilings. For instance, it is easily seen from
Figure~\ref{fig:lattice:example} that while using only flips suffices to
generate all the tilings of a polygon, it also generates multiple times the
same tiling.

Not all collections of squares in $\mathbb{Z}^3$ correspond, once projected
onto the plane, to tilings in the triangular grid: see for instance
Figure~\ref{fig:compact:ph}~(a). So we need conditions on the squares. In
fact, we will build our reasoning not on squares but on cubes, because the
elementary operation on a tiling, the flip, has a natural, visual and
intuitive translation in terms of cubes (see Figure~\ref{squares:cubes}). 
Moreover, fracture lines (see Section~\ref{fracture:lines}) can be used to get 
rid of squares that do not correspond to cubes.

\fig{
	\begin{tabular}{p{2cm}@{\hskip5mm}p{2cm}@{\hskip2mm}
		p{3cm}@{\hskip2mm}p{2cm}}

	\mbox{}\hfill
	\epsfig{file=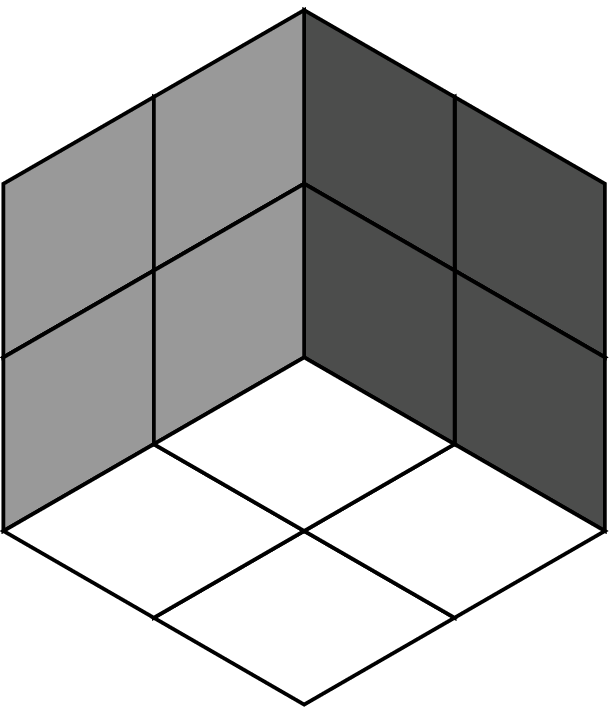, width=1.55cm}
	\hfill\mbox{}
	&
	\mbox{}\hfill
	\epsfig{file=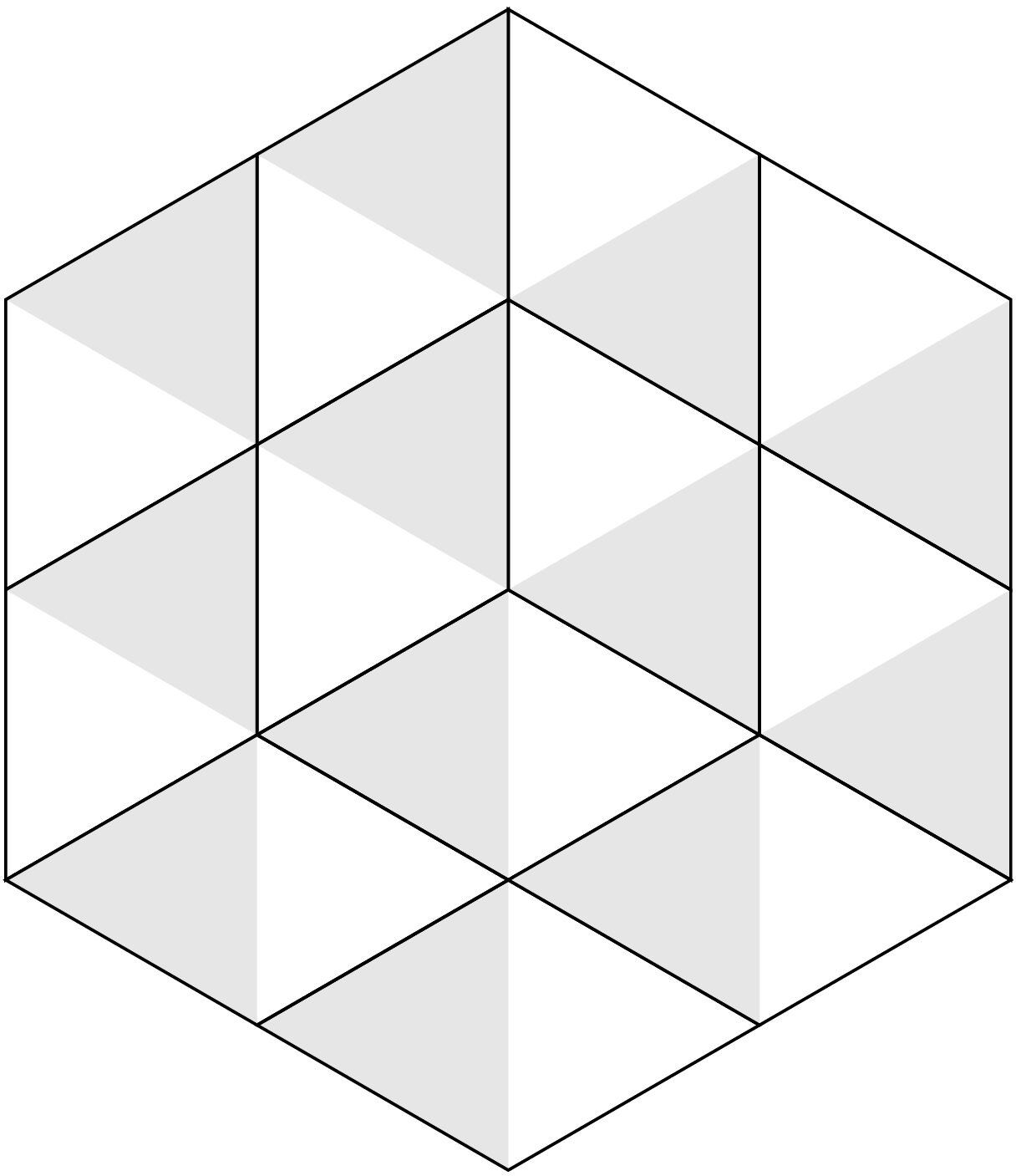, width=1.55cm}
	\hfill\mbox{}
	&
	\mbox{}\hfill
	\epsfig{file=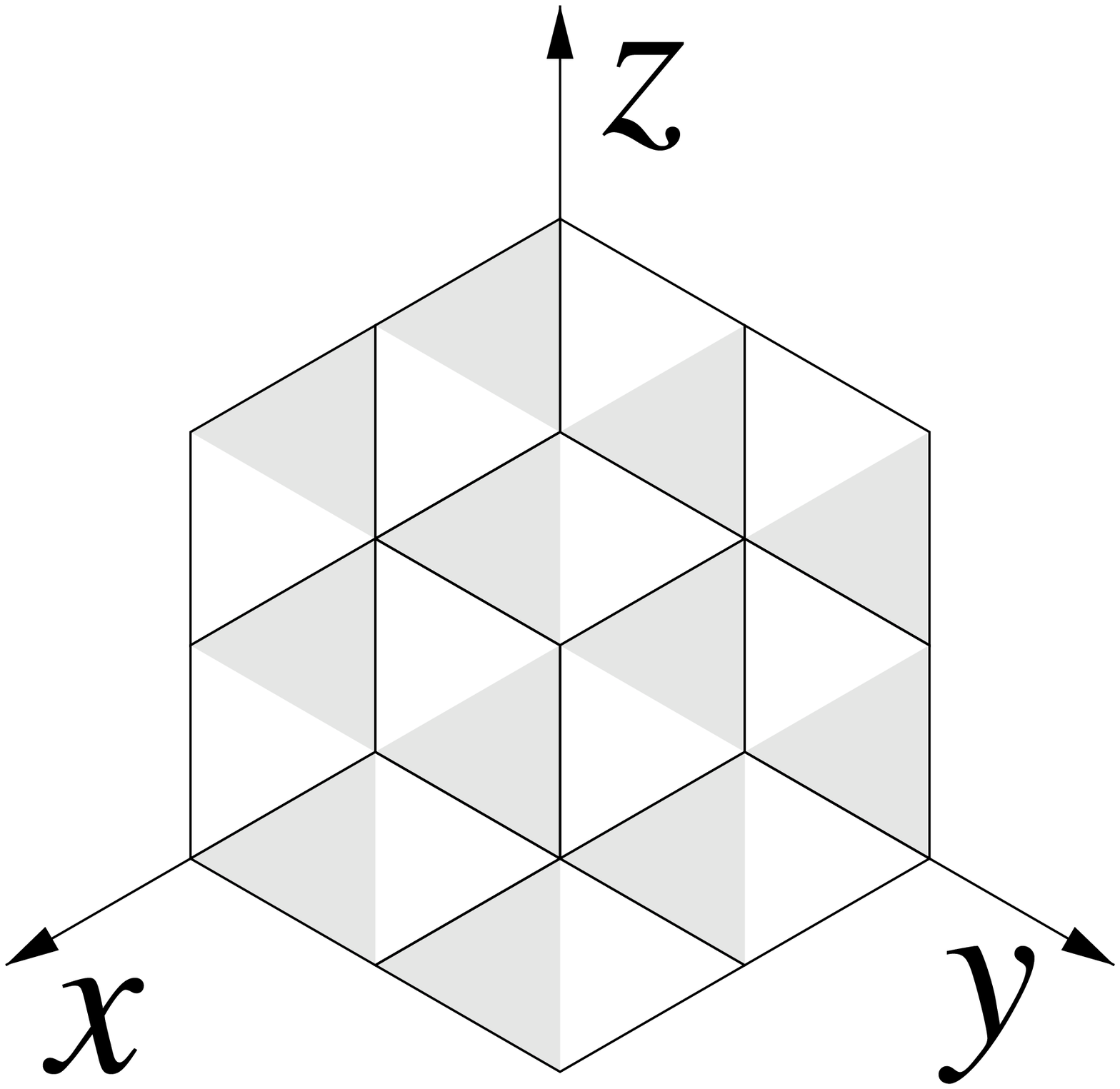, width=2.5cm}
	\hfill\mbox{}
	&
	\mbox{}\hfill
	\epsfig{file=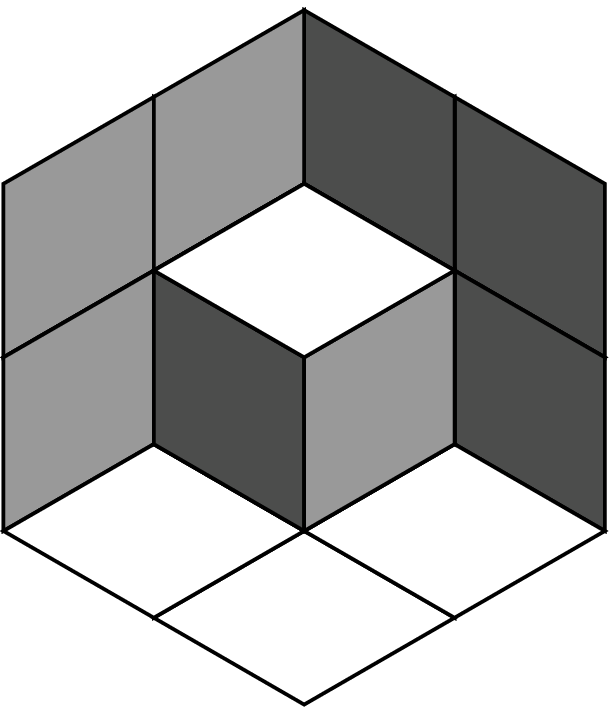, width=1.55cm}
	\hfill\mbox{}

	\end{tabular}
	\caption{An up-flip is equivalent to adding a cube 
		\worklabel{squares:cubes}}
}

The minimal tiling of a polygon has no local maximum in the interior of the
polygon, otherwise this maximum could be flipped down. In other words, there
is no cube associated to such a tiling. It merely provides squares on which
one can put cubes. In the case of a hexagon, the association of the squares
can be seen as base planes in $(\mathbb{Z}_+)^3$ so that any cube can be
encoded by the (integer) coordinates of its lowest corner.

Cubes can be piled, but if we want to preserve a tiling after projection, we
must pile them in a way that corresponds to flips. Such piles will be called
compact. 

\begin{definition}[Compact pile]
	A pile of cubes is \emph{compact} if for any cube of the pile at
$(i_0,j_0,k_0)$ there are cubes at $(i,j_0,k_0)$, $(i_0,j,k_0)$ and
$(i_0,j_0,k)$ with $i$, $j$ and $k$ ranging in $[0,i_0]$, $[0,j_0]$ and
$[0,k_0]$.
	\worklabel{def:compact:pile}
	\end{definition}

The next lemma easily stems from this definition:

\begin{lemma}
	A pile of cubes is compact if and only if any section of the pile by a
plane orthogonal to one axis ($Ox$, $Oy$ or $Oz$) yields a Ferrers diagram.
	\worklabel{lemma:compact:Ferrers}
	\end{lemma}

It will not be enough to consider only ``perfect'' hexagons of size $n \times
n \times n$. The hexagon-like sub-polygons that appear in tilings require a
somewhat more general approach (see Section~\ref{section:seeds}):

\fig{
	\begin{tabular}{p{5cm}@{\hskip7mm}p{6cm}}

	\mbox{}\hfill
	\epsfig{file=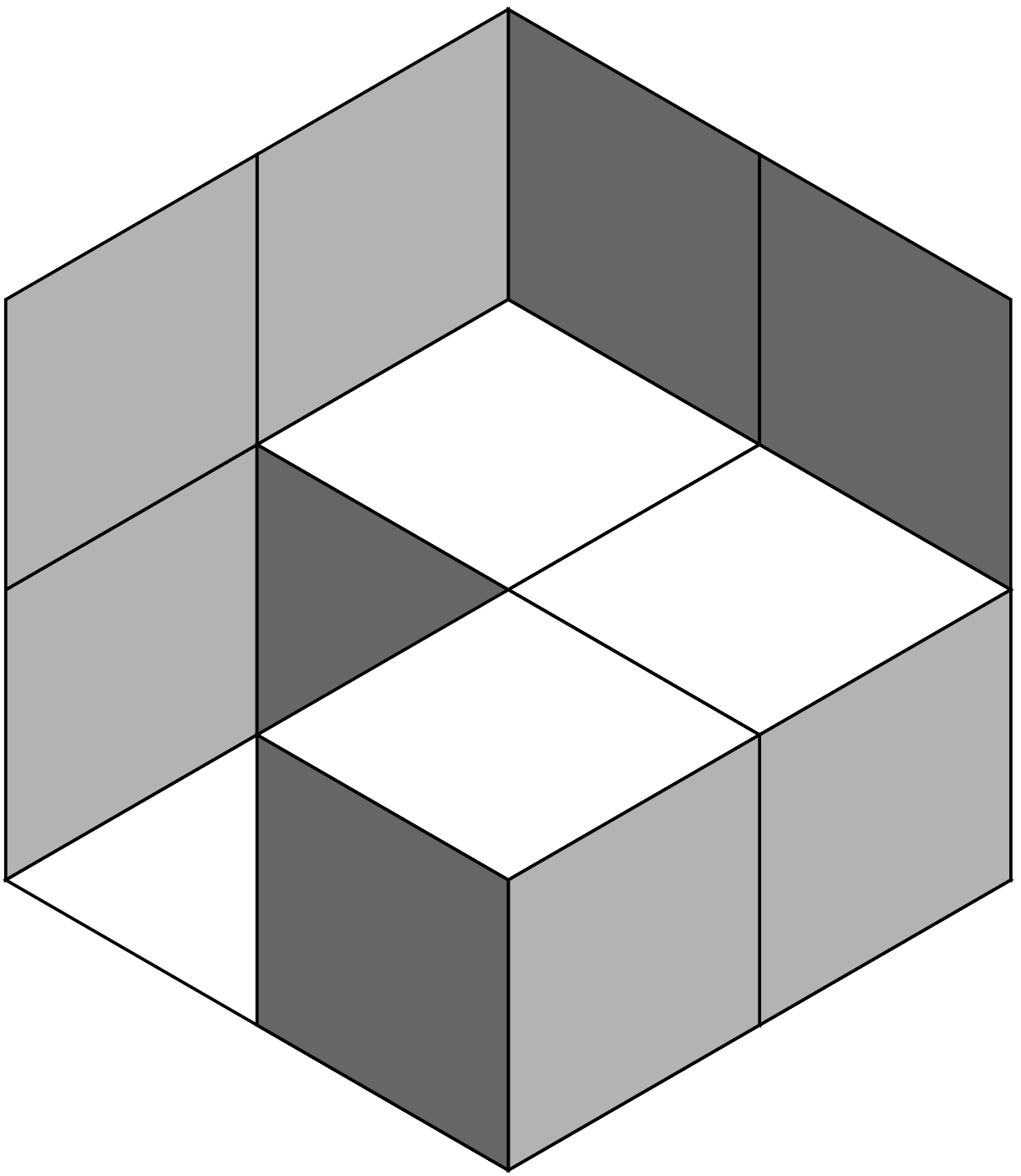, width=1cm}
	\hfill\mbox{}
	&
	\mbox{}\hfill
	\epsfig{file=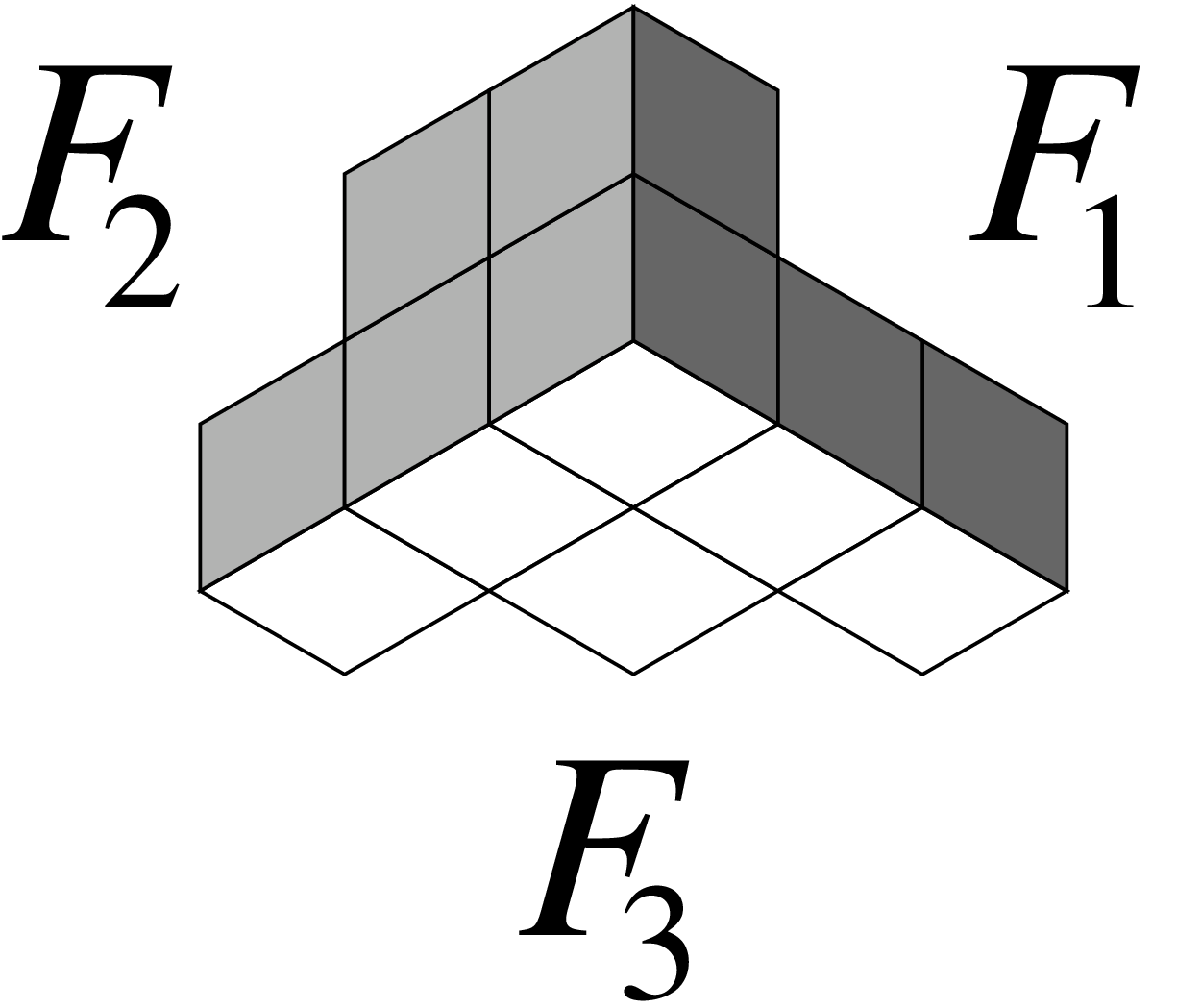, width=1.7cm}
	\hfill\mbox{}
	\\[1mm]

	\mbox{}\hfill
	(a) 	\begin{minipage}[t]{4.2cm}
		Neither a tiling nor a \emph{compact} pile of cubes
		\end{minipage}
	\hfill\mbox{}
	&
	\mbox{}\hfill
	(b) 	\begin{minipage}[t]{5.0cm}
		A pseudo-hexagon delimited by three Ferrers diagrams
		\end{minipage}
	\hfill\mbox{}
	\\

	\end{tabular}
	\caption{Compact piles and pseudo-hexagons
		\worklabel{fig:compact:ph}
	}
}

\begin{definition}[Pseudo-hexagon]
	$\!$A \emph{pseudo-hexagon} is the domain obtained by starting from a
compact pile of cubes and performing all the possible down-flips.

	\worklabel{def:ph}
	\end{definition}

Note that by Lemma~\ref{lemma:compact:Ferrers} a pseudo-hexagon can also be
seen as a figure delimited by three Ferrers diagrams fitting neatly (see
Figure~\ref{fig:compact:ph}~(b)). It is also tileable by construction.

\begin{proposition}
	A pile of cubes in a peudo-hexagon corresponds to a tiling of the
pseudo-hexagon if and only if it is compact.

	\worklabel{prop:tilings:cp}
	\end{proposition}

\noindent \textit{Proof} \quad A pile of cubes corresponds to a tiling of a
domain ${\cal D}$ if and only if it can be generated (starting from the
minimal tiling) by using only up-flips (see
Corollary~\ref{coro:lattice:structure}). An up-flip adds a cube that lies on
squares that belong either to the minimal tiling or to an already added cube.
The result follows by induction on the number of cubes and
Lemma~\ref{lemma:compact:Ferrers}.~\mbox{}\hfill$\square$

%!!!!!!!!!!!!!!!!!!!!!!!                         !!!!!!!!!!!!!!!!!!!!!!!
%!!!!!!!!!!!!!!!!!!!!!!!   Algo pour hexagones   !!!!!!!!!!!!!!!!!!!!!!!
%!!!!!!!!!!!!!!!!!!!!!!!                         !!!!!!!!!!!!!!!!!!!!!!!

\subsection{An algorithm to generate the tilings of a\\ pseudo-hexagon}
\worklabel{section:algo:hexagons}

Proposition~\ref{prop:tilings:cp}
means that generating the tilings of a pseudo-hexagon is
equivalent to generating compact arrangements of cubes, which is quite easier.
Indeed, the latter are related one-to-one with plane partitions (see
Definition~\ref{def:plane:partition}, Figure~\ref{fig:partitions}~(a)
and~\cite{Elser}). 

\begin{definition-notation}
	The parts of the partition~$p$ associated with a Ferrers diagram~$F$
are noted~$p[j]$, and~$F[j]$ denotes a collection of~$p[j]$ linearly adjacent
squares of~$\mathbb{Z}^2$. A partition on~$F[j]$ is a non-increasing sequence
of integers placed on this squares. Two partitions~$p$ and~$q$ are
\emph{comparable} and \mbox{$p \leqslant q$} if $p[j] \leqslant q[j]$ for
all~$j$.

A plane partition~$P$ on a Ferrers diagram~$F$ is a non-increasing sequence of
partitions on the~$F[j]$. The partition on~$F[j]$ is noted~$P[j]$ and called
the \emph{$j$-slice} of~$P$.  Two plane partitions~$P$ and~$Q$ are
\emph{comparable} and \mbox{$P \leqslant Q$} if $P[j] \leqslant Q[j]$ for
all~$j$.

The \emph{weight} of a partition is its sum; the \emph{weight} of a plane 
partition is its~sum.

	\worklabel{def:order:PP}
	\end{definition-notation}

\fig{
	\begin{tabular}{p{5cm}@{\hskip1cm}p{5cm}}

	\mbox{}\hfill
	\epsfig{file=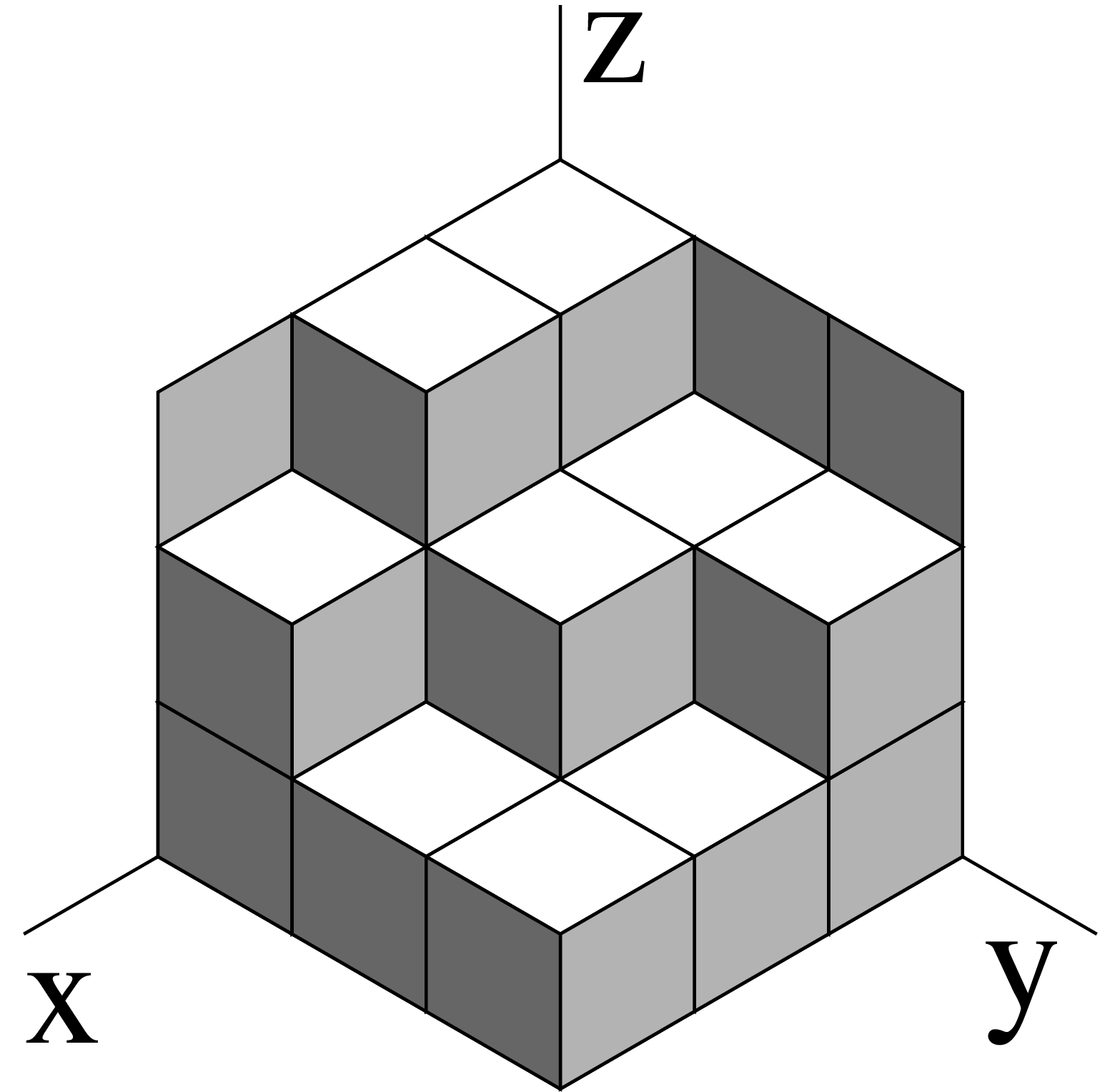, width=2cm}
	\raisebox{0.2cm}{
		\epsfig{file=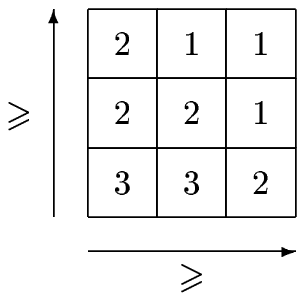, width=2cm}
	}
	\hfill\mbox{}
	&
	\mbox{}\hfill
	\epsfig{file=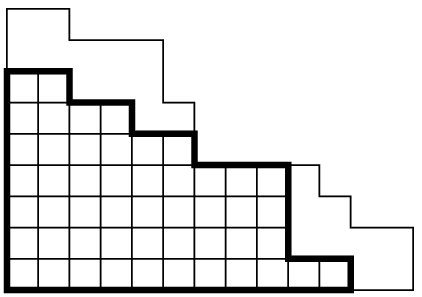, width=2.5cm}
	\hfill\mbox{}
	\\[1mm]

	\mbox{}\hfill
	(a) 	\begin{minipage}[t]{4cm}
		Tilings are related one-to-one to plane partitions
		\end{minipage}
	\hfill\mbox{}
	&
	\mbox{}\hfill
	(b) 	\begin{minipage}[t]{4cm}
		A limited partition
		\end{minipage}
	\hfill\mbox{}
	\\

	\end{tabular}
	\caption{Plane partitions and limited partitions 
		\worklabel{fig:partitions}
	}
}

The plane partitions we are dealing with correspond to the tilings of a
pseudo-hexagon; they are thus limited, in the sense that the corresponding
pile must be embedded in the pile that defines the pseudo-hexagon.

\begin{definition}
	A $(s,P)$ \emph{limited plane partition} is a plane partition~$A$ of
the integer~$s$ such that the plane partitions~$A$ and~$P$ are comparable and
\mbox{$A \leqslant P$}.

	\worklabel{def:limited:plane:partitions}
	\end{definition}

We will generate all $(s,P)$ limited plane partitions recursively (see
Algorithm~\ref{algo:tilings}), adding one $j$-slice at a time. To this end, we
need to generate the partitions that are less than a given partition.

\begin{definition}
	A $(s,F)$ \emph{limited partition} is a partition
$a$ of the integer~$s$ such that~$a \leqslant p$ where~$p$ is the partition 
associated with the Ferrers diagram~$F$.

	\worklabel{def:limited:partitions}
	\end{definition}

Figure~\ref{fig:partitions}~(b) shows of a graphical representation of a
limited partition.

\medskip

The following recursive algorithm uniquely generates all $(s,F)$ limited
partition by adding one by one the parts of the final partitions. The average
step is thus to add a row with~$k$ squares to a yet unfinished Ferrers
diagram; the value of this part must at least be permitted by the ``geometry''
of~$F$ and be less than the last added part. The value must also allow the
updated~$s$ to fit in the remaining ``space''; this is the meaning
of~$r_{q+1}$ below.

\begin{algorithm}
\worklabel{algo:limited:standard:partitions}
\end{algorithm}

\begin{itemize}

\item \textbf{Input}: Integers $s, p_1, \ldots, p_n$ such that $s \leqslant
p_1 + \cdots + p_n$.

\item \textbf{Output}: List $\cal L$ of all $(s,(p_1, \ldots , p_n))$ limited
partitions.

\item \textbf{Compute} $\text{loop}(s, \emptyset, (p_1, \ldots, p_n))$ where

$\text{loop}(t, (a_1, \ldots, a_q), (p_{q+1}, \ldots, p_n)) =$

\begin{itemize}

\item if $t=0$, ${\cal L} \leftarrow (a_1, \ldots , a_q , \underbrace{0,
\ldots , 0}_{(n-q)\ \text{times}}\!\!)$;

\item else if $n = q+1$, ${\cal L} \leftarrow (a_1, \ldots , a_q, t)$;

\item else compute 

$\text{loop}(t-k, (a_1, \ldots , a_q, k), (\text{min}(p_{q+2},k) , \ldots ,
\text{min}(p_n,k)))$

for $k = \text{min}(p_{q+1},t,a_q) \text{ downto } r_{q+1}$

where $r_{q+1} = \text{min}\{ 1 \leqslant j \leqslant 
\text{min}(p_{q+1},t,a_q) \: |$

\mbox{}\hskip3cm $j + \text{min}(p_{q+2},j) + \cdots + 
\text{min}(p_n,j) \geqslant t \}$.

\end{itemize}

\end{itemize}

\noindent \textit{Proof of the algorithm} \quad A more basic way to write the
algorithm would be to add any number of squares ($k = 1 ..  
\text{min}(p_{q+1},t,a_q)$) at each step: the min simply ensures that the
resulting collection of squares is indeed a Ferrers diagram (the new part must
be less than its counterpart in the limiting partition ($p_{q+1}$), less than
the number of remaining squares ($t$) and less than the last added part
($a_q$)) and this procedure obviously generates all the desired partitions.

But it is also wasteful, in that many such combinations would not satisfy $t
\leqslant p_{q+1} + \cdots + p_n$\,, so we restrict the range of~$k$ until its
minimal value guarantees, by construction, that the following steps will
result in a partition limited by~$p$. Since the values we eliminate would not 
ultimately yield such a partition, we still generate all the desired 
partitions.

Finally, our algorithm executes a depth-first search on a dynamically
generated tree $\cal T$ of height~$n$ and whose $p{}^{\text{th}}$ level
features all the possible combinations for the~$p$ first rows (the other rows
being empty), so that no partition is obtained twice.~\mbox{}\hfill$\square$

\medskip

\noindent \textit{Complexity} \quad Since we perform a depth-first search
on~$\cal T$, the execution space of the algorithm is bounded by the height
of~$\cal T$ ($n$, the number of parts in the limiting partition) times the
number of leaves, which is the number of~$(s,(p_1, \ldots , p_n))$ limited
partitions. We believe no closed formula is known for the latter, but
generating series techniques (see for instance \cite{Flajolet} and
\cite{Flaj}) yield lemma~\ref{lemma:count:lp}:

\begin{definition}
	The \emph{multidegree} of a monomial ${u_1}^{i_1} \cdots {u_n}^{i_n}$
is the sequence $(i_1, \ldots , i_n)$. Its \emph{cumulated multidegree} is the
sequence $(i_1 + \cdots + i_n, i_2 + \cdots + i_n, \ldots, i_n)$.

	\worklabel{def:multidegree}
	\end{definition}

\begin{lemma}
	The number of $(s,p=(p_1, \ldots, p_n))$ limited partitions is the
number of monomials whose cumulated multidegree is less than the conjugate
partition of $p$ in $$P(u_1, \ldots , u_n) = [z^s] \prod_{k=1}^n \dfrac{1}{1 -
u_k\,z^k}$$

	\worklabel{lemma:count:lp}
	\end{lemma}

\noindent \emph{Proof (Sketch)} \quad The exponent of $u_k$ in a monomial
${u_1}^{i_1} \cdots {u_p}^{i_n}$ in $P(u_1, \ldots , u_n)$ counts the number
of rows of length~$k$ in the corresponding partition (see
Figure~\ref{fig:counting}). This monomial encodes a partition limited by~$p$
if and only if it satisfies $i_r + i_{r+1} + \cdots + i_n \leqslant m(r) +
m(r+1) + \cdots + m(n)$ for $1 \leqslant r \leqslant n$ where $m(k)$ is the
number of occurences of~$k$ in the conjugate partition
of~$p$.~\mbox{}\hfill$\square$

\begin{figure}[ht]
	\begin{center}
	\epsfig{file=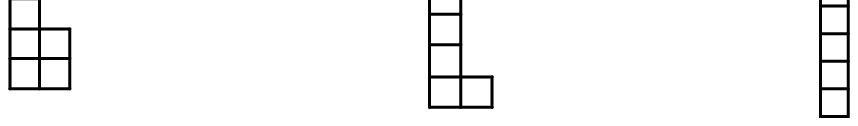, width=8.5cm}
	\caption{Generating series and Ferrers diagrams
		\worklabel{fig:counting}}
	\end{center}
	\end{figure}

Finally, the execution time of the algorithm is proportionnal to its execution
space, the multiplying factor being the time needed to compute~$r_{q+1}$ at
each step, which is at most~$n \times p_1$ (exhaustive
search).~\mbox{}\hfill$\square$

\medskip

Note that since truncated generating series can be computed by standard
programs such as Maple, lemma~\ref{lemma:count:lp} provides another
(sub-optimal, but handy) way to generate limited partitions.

\medskip

%\fig{
%	\begin{tabular}{p{5cm}@{\hskip5mm}p{5.5cm}}
%
%	\mbox{}\hfill
%	\raisebox{1.5cm}{
%		\epsfig{file=dessins/before_after_collapse.ps, width=4.7cm}
%	}
%	\hfill\mbox{}
%	&
%	\mbox{}\hfill
%	\epsfig{file=dessins/arbre.ps, width=5cm, clip=} 
%	\hfill\mbox{}
%	\\[1mm]
%
%	\mbox{}\hfill
%	(a) 	\begin{minipage}[t]{3.5cm}
%		Collapsing stacks into $j$-slices
%		\end{minipage}
%	\hfill\mbox{}
%	&
%	\mbox{}\hfill
%	(b) 	\begin{minipage}[t]{4.5cm}
%		Algorithm 2 for 5 cubes in a box of side 2
%		\end{minipage}
%	\hfill\mbox{}
%	\\
%
%	\end{tabular}
%	\caption{Plane partitions and cubes \worklabel{fig:partitions:cubes}}
%}

We now undertake to generate all the tilings of a pseudo-hexagon, or rather
all plane partitions limited by a given plane partition~$P=(P_1 , \ldots ,
P_n)$ where the~$P_j$'s are the $j$-slices of~$P$. It is enough to generate
all $(s,P)$ limited plane partitions and to let~$s$ range from~0 to the
weight~$w(P)$ of~$P$.

\begin{algorithm}
\worklabel{algo:tilings}
\end{algorithm}

\begin{itemize}

\item \textbf{Input}: An integer $s$ and a plane partition~$P = (P_1 , \ldots
, P_n)$ such that $s \leqslant w(P)$.

\item \textbf{Output}: The list $\cal L$ of all $(s,P)$ limited plane
partitions.

\item \textbf{Compute} $\text{loop}(s,\emptyset,(P_1,\ldots,P_n))$ where

$\text{loop}(t,(A_1, \ldots ,A_q) , (P_{q+1} , \ldots, P_n))=$

	\begin{itemize}

	\item if $t=0$, ${\cal L} \leftarrow (A_1, \ldots, A_q,
		\underbrace{0,\ldots,0}_{(n-q)\text{ times}})$;

	\item else if $n = q+1$,

		\begin{itemize}
		\item generate all $(t, \text{min}(A_q,P_{q+1}))$ limited 
			partitions with 
			Algo\-rithm \ref{algo:limited:standard:partitions};
		\item for each such partition~$p$, ${\cal L} \leftarrow
			(A_1, \ldots, A_q,p)$;
		\end{itemize}

	\item else for all partition~$k$ such that

		$$\begin{cases}
		w(k) \leqslant t \\
		k \leqslant \text{min}(A_q,P_{q+1}) \\
		k + w(\text{min}(k,P_{q+1})) + \cdots + w(\text{min}(k,P_n))
			\geqslant t \\
		\end{cases}$$

		compute 

		$\text{loop}(t-w(k),(A_1,\ldots,A_q,k), 
		(\text{min}(k,P_{q+2}), \ldots, \text{min}(k,P_n)))$;

	\end{itemize}

\end{itemize}

\noindent \textit{Proof and complexity of the algorithm} \quad The analysis is
completely analogous to Algorithm \ref{algo:limited:standard:partitions}'s.  
The execution space is at most~$n$ times the number of~$(s,P)$ limited plane
partitions (for which we believe no closed formula is known, but as above a
characterization with generating series techniques can be obtained) and the
execution time is at most the execution space times~$n$ times the number of
$(u,P_1)$ limited partitions for~$u=0..w(P_1)$ (exhaustive
search).~\mbox{}\hfill$\square$

\medskip

Note that Algorithm \ref{algo:limited:standard:partitions} is a special case
of Algorithm~\ref{algo:tilings}, in the same manner that a partition is a
special case of plane partition; we presented both for the sake of clarity.

%%%%%%%%%%%%%%%%%%%%%%%%%%%%%%%%%%%%%%%%%%%%%%%%%%%%%%%%%%%%%%%%%%%%%%%%
%%%%%%%%%%%%%%%%%%%%%%%%%%%%%%%%%%%%%%%%%%%%%%%%%%%%%%%%%%%%%%%%%%%%%%%%
%%%%%%%%%%%%%%%%%%%%                                 %%%%%%%%%%%%%%%%%%%
%%%%%%%%%%%%%%%%%%%%   Encore quelques définitions   %%%%%%%%%%%%%%%%%%%
%%%%%%%%%%%%%%%%%%%%                                 %%%%%%%%%%%%%%%%%%%
%%%%%%%%%%%%%%%%%%%%%%%%%%%%%%%%%%%%%%%%%%%%%%%%%%%%%%%%%%%%%%%%%%%%%%%%
%%%%%%%%%%%%%%%%%%%%%%%%%%%%%%%%%%%%%%%%%%%%%%%%%%%%%%%%%%%%%%%%%%%%%%%%

\section{Domains, fracture lines and seeds}
	\worklabel{domains:fracture:lines:seeds}

Up to now, we have considered polygons, which are hole-free (simply connected)
domains on which one knows how to define a consistent height function, an
essential ingredient of Section~\ref{flips:lattices}. The subsequent results
in this paper hold provided the results in Section~\ref{flips:lattices} do,
but this does not necessarily mean one has to deal only with simply connected
domains (see Section~\ref{domains:with:holes}). It would be quite exciting to
examine how works like \cite{Fournier1} or \cite{Fournier2} extend to the
results of this paper.

For the rest of this paper, we consider domains for which the results in
Section~\ref{flips:lattices} hold. To reflect this fact, we will not use the
word ``polygon''.

\medskip

A domain ${\cal D}$ is defined by its contour path, which is simply any
finite-length closed path $\cal P$ in the triangular grid.  \emph{Fracture
lines} (see Section~\ref{fracture:lines}) will allow us to break the domain
into several parts (called \emph{fertile zones}) whose tilings can be
generated independently. The lattice of the tilings of ${\cal D}$ can then be
obtained as the product of the lattices of the tilings of the fertile zones.
Finally, \emph{seeds} (see Section~\ref{section:seeds}) will allow us to
decompose fertile zones into pseudo-hexagons, for which much is already known.

%!!!!!!!!!!!!!!!!!!!!!!!                        !!!!!!!!!!!!!!!!!!!!!!!!
%!!!!!!!!!!!!!!!!!!!!!!!    Domaines à trous    !!!!!!!!!!!!!!!!!!!!!!!!
%!!!!!!!!!!!!!!!!!!!!!!!                        !!!!!!!!!!!!!!!!!!!!!!!!

\subsection{Domains with holes}
\worklabel{domains:with:holes}

A domain defined by a finite-length closed path ${\cal P}$ in the triangular
grid need not be hole-free: it suffices that ${\cal P}$ cross itself in a
non-trivial way (see Figure~\ref{fig:holes}~(a) for an example). Such a domain
may still be tileable (see Figure~\ref{fig:holes}~(b)). The definitions and
results in this paper apply to such domains provided the results of
Section~\ref{flips:lattices} hold.

\fig{
	\begin{tabular}{p{4cm}@{\hskip1cm}p{5.5cm}}

	\mbox{}\hfill
	\epsfig{file=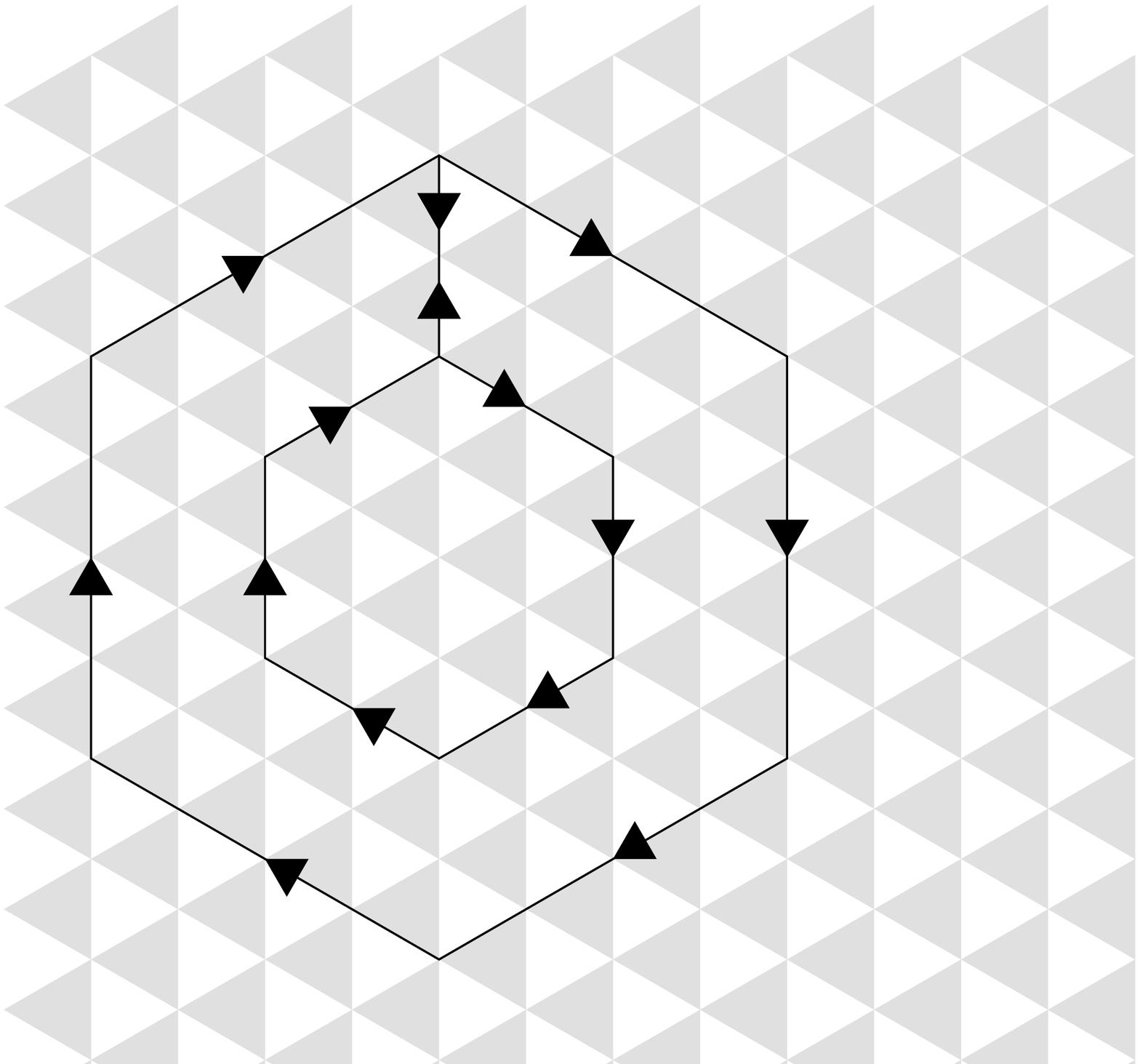, width=1.5cm, clip=}
	\hfill\mbox{}
	&
	\mbox{}\hfill
	\epsfig{file=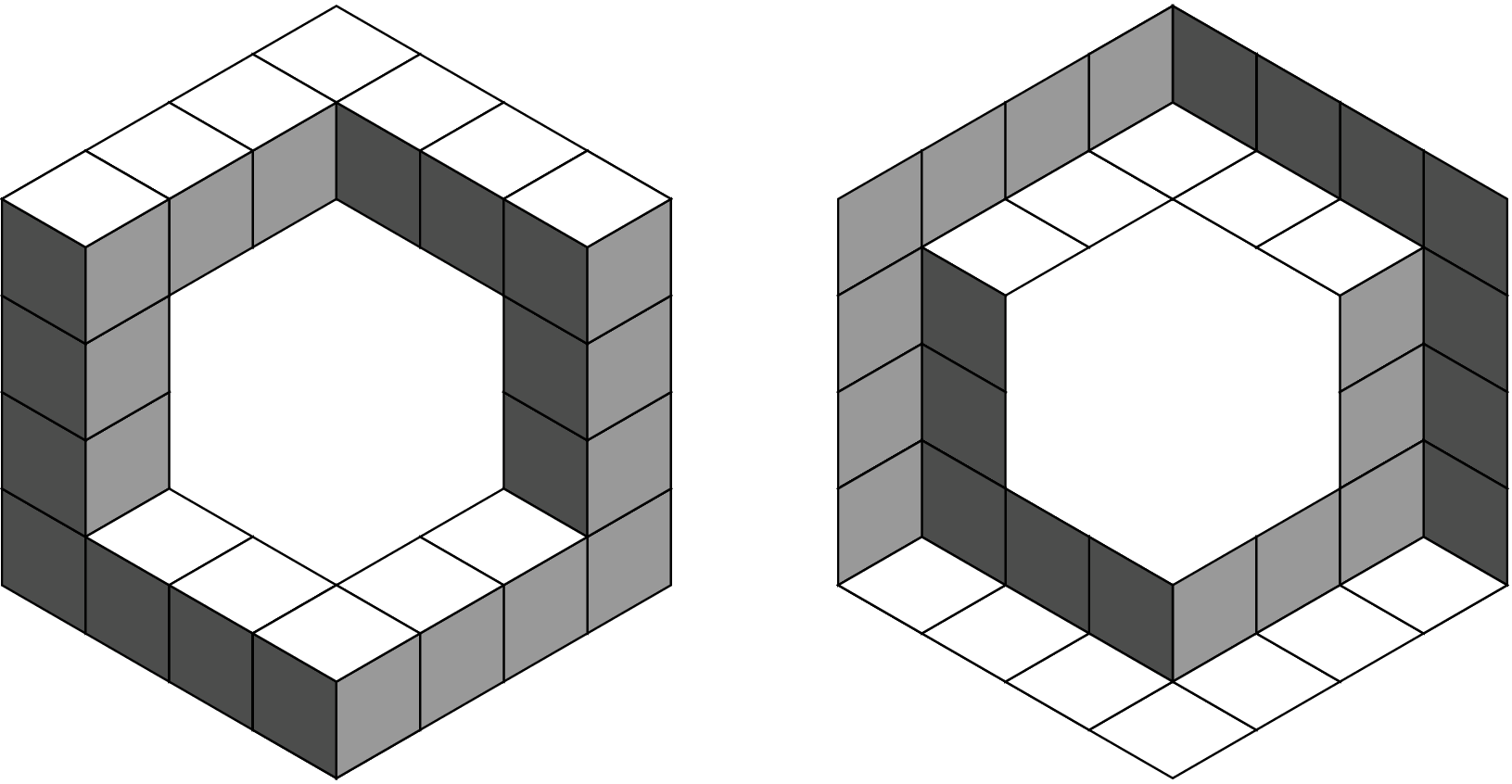, width=3cm, clip=}
	\hfill\mbox{}
	\\[1mm]

	\mbox{}\hfill
	(a) 	\begin{minipage}[t]{3.2cm}
		A path delimiting a domain with a hole
		\end{minipage}
	\hfill\mbox{}
	&
	\mbox{}\hfill
	(b) 	\begin{minipage}[t]{4cm}
		This domain can be tiled even though it has a hole
		\end{minipage}
	\hfill\mbox{}
	\\

	\end{tabular}
	\caption{Domains with holes \worklabel{fig:holes}}
}

%!!!!!!!!!!!!!!!!!!!!!!!                        !!!!!!!!!!!!!!!!!!!!!!!
%!!!!!!!!!!!!!!!!!!!!!!!   Lignes de fracture   !!!!!!!!!!!!!!!!!!!!!!!
%!!!!!!!!!!!!!!!!!!!!!!!                        !!!!!!!!!!!!!!!!!!!!!!!

\subsection{Fracture lines, fertile zones and the fracture algorithm}
\worklabel{fracture:lines}

When considering a general (tileable) domain ${\cal D}$, it may happen that a
given tile is ``fixed'', in the sense that it cannot be changed by a flip, so
we'd like to remove such a tile and concentrate on ``interesting'' zones (see
Section~\ref{example}). The tool for this is the fracture line, made of
vertices whose height does not vary between the minimal and maximal tilings.  
It allows us to ``remove'' fixed lozenges but there's more: it also allows us
to define disjoint sub-domains which can be tiled independently (see
Section~\ref{example}), so that the lattice of the tilings of $\cal D$ can be
obtained by computing a product of smaller, simpler lattices (see
Section~\ref{subsection:application:FT}). This section attempts to clarify all
of the above, step by step, until a complete algorithm can be written, proved
and analyzed.

Fracture lines have been investigated in the case of dominoes, notably in
\cite{Fournier3}. In all this section, $\cal D$ is assumed to be a tileable
domain such that the results of Section~\ref{flips:lattices} hold.

\begin{definition}
	Let $T$ be any tiling of $\cal D$ and let $h$ be the associated height
function. A path $P$ in $\cal D$ is \emph{$T$-valid} if it is a single
vertex of $\cal D$ or if $|h(v') - h(v)|=1$ for $v$, $v'$ two consecutive
vertices of~$P$ (\textit{i.e.} it follows the edges of lozenges in~$T$).

	A path $P$ in $\cal D$ is a $T$-valid cycle if it is a closed,
repetition-free path in $\cal D$ and a $T$-valid path.

	A subset $F$ of $\cal D$ is a \emph{sub-domain} of $\cal D$ if there
exists a tiling $T$ of $\cal D$ and a $T$-valid cycle $C$ such that $C$ is the
contour path of~$F$.
	\worklabel{def:sub:domain}
	\worklabel{def:valid:cycle}
	\end{definition}

\begin{definition}
	Let ${\cal D}$ be any tileable domain and $v$ one of its vertices.  
We denote by $\Delta h(v)$ the difference between the height of $v$ in the
maximal tiling of ${\cal D}$ and its height in the minimal tiling. If $\Delta
h(v)=0$ then $v$ is said to be \emph{solid}.
	\worklabel{def:delta:height}
	\end{definition}

\begin{definition}[Fracture line]
	A \emph{fracture line} in ${\cal D}$ is a $T$-valid cycle in~$\cal D$
(for some tiling~$T$ of~$\cal D$) that is composed only of solid points. Every
domain ${\cal D}$ has at least one such fracture line, namely its contour
path, which we will call the \emph{trivial fracture line} of ${\cal D}$.
	\worklabel{def:fracture:line}
	\end{definition}

We now undertake to prove Theorem~\ref{theo:fracture:zones}. We first prove
(see Proposition~\ref{prop:trivial:fracture:lines}) that a solid point in the
interior of~$\cal D$ cannot be isolated: it belongs at least to a solid
lozenge.

\begin{lemma}
	Let $F$ be a sub-domain of $\cal D$, delimited by a valid cycle in a
tiling $T$ of $\cal D$. If a triangle $t$ belongs to $F$, then so does the
lozenge $\ell$ to which $t$ belongs in $T$.
	\worklabel{lemma:triangle:lozenge}
	\end{lemma}

\noindent \textit{Proof} \quad Let $a$, $b$, $c$ be the vertices of $t$, and
let $d$ be the fourth vertex of $\ell$. There are exactly two couples made of
vertices in $\{a,b,c\}$ which have a height difference equal to 1 in $T$. Let
$a$ and $b$ be the two vertices which have a height difference of 2 (in
absolute value), so that $[a;b]$ is not a valid path in $T$. Since $F$ is
delimited by a $T$-valid cycle, $[a;b]$ cannot be part of the
contour path of $F$. Therefore $d$ belongs to $F$, whence the triangle $abd$
belongs to $F$, and since $\ell$ is the union of the triangles $abc$ and
$abd$, we are done.~\mbox{}\hfill$\square$

\begin{proposition}
	Any sub-domain of a tileable domain is itself tileable.

	\worklabel{prop:tileable:valid:path}
	\end{proposition}

\noindent \textit{Proof} \quad Let $F$, a sub-domain of a tileable domain
$\cal D$, be delimited by a valid cycle $L$ in a tiling $T$ of~$\cal D$.

If $F$ is a single point, we are done (using no tile). Otherwise, let $s$ be a
segment of $L$; it belongs to 2 triangles in the triangular grid, one of which
must belong to $F$. By Lemma~\ref{lemma:triangle:lozenge}, the lozenge $\ell$
to which this triangle belongs in $T$ belongs to $F$. If $F=\ell$, we are
done.

Otherwise, let $F_1$ be the subset of $\cal D$ obtained by removing $\ell$
from $F$. It is \emph{a priori} composed of several disjoint subsets of $\cal
D$, but each of these subsets is delimited by a $T$-valid cycle so that the
situation is equivalent to a single domain $F_1$, which is then a sub-domain
of $\cal D$.

Let $F_p$ be the domain obtained by recursively removing a lozenge until
nothing can be done. If $F_p$ contained a triangle it would also contain a
lozenge by Lemma~\ref{lemma:triangle:lozenge} and so it would not be a fixed
point. Therefore it contains no triangle, which means that it is a single
point and~$F$ is a collection of lozenges and therefore
tileable.~\mbox{}\hfill$\square$

\begin{corollary}
	A fracture line of $\cal D$ delimits a tileable sub-domain of $\cal D$.

	\worklabel{coro:fracture:line:tileable}
	\end{corollary}

\noindent \textit{Proof} \quad Indeed, a fracture line is a valid cycle so it
delimits a sub-domain.~\mbox{}\hfill$\square$

\begin{lemma}
	Let $s_1, \ldots, s_6$ be six vertices of $\cal D$ that define a
hexagon of side~1. Let $c$ be the center of this hexagon. If $\Delta h (s_k)
\neq 0$ for $k=1..6$, then $\Delta h(c) \neq 0$.

	\worklabel{lemma:hexagon:center}
	\end{lemma}

In other words, if $c$ is a vertex in the interior of $\cal D$ that cannot be
flipped, then at least one of its neighbours cannot be flipped.

\medskip

\noindent \textit{Proof} \quad Let $h_{\text{min}}$ and $h_{\text{max}}$
denote the minimal and maximal height functions. Since $c$ and $s_k$ belong to
the same triangle, $h_{\text{min}}(c) \leqslant h_{\text{min}}(s_k) + 2$ (see
Section~\ref{height:functions}). Moreover $h_{\text{max}}(s_k) \geqslant
h_{\text{min}}(s_k)+3$ since $\Delta h(s_k) \neq 0$, from which we derive
$h_{\text{min}}(c) \leqslant h_{\text{max}}(s_k) - 1$. In the maximal tiling,
no vertex can be a local minimum (see Proposition~\ref{prop:local:extremum}
and Definition~\ref{def:flip}), so $h_{\text{max}}(c)$ cannot be less than
$h_{\text{min}}(s_k) - 1$ for all $k=1..6$, whence $h_{\text{max}}(c) \neq
h_{\text{min}}(c)$.~\mbox{}\hfill$\square$

\begin{proposition}
	Let $\cal D$ be any tileable domain, not limited to a single point,
and let $v$ be a solid vertex of $\cal D$. Then $v$ belongs to a fracture line
that encloses at least a lozenge.

	\worklabel{prop:trivial:fracture:lines}
	\end{proposition}

\noindent \textit{Proof} \quad The result is trivial if $v$ is on the boundary
path of $\cal D$. We now assume that we are not in this case. By
Lemma~\ref{lemma:hexagon:center} we know that there exists a solid vertex $v'$
at distance 1 (in the triangular grid) of $v$. Let $T$ be any tiling of $\cal
D$ and let $h$ be the associated height function. If $|h(v)-h(v')|$ is equal
to 2 then the lozenge whose middle segment is $[v;v']$ is defined by two
fracture points and we are done. 

If $|h(v)-h(v')|=1$ then $(v,v')$ is a valid path in~$T$. Let now $s_1,
\ldots, s_5$ and $v'$ be the (distinct) vertices of the hexagon around $v$. We
prove \textit{ab absurdo} that at least one of the $s_k$ is solid: assume that
none of them is solid, so that $h_{\text{min}}(s_k) \neq h_{\text{max}}(s_k)$
for $k=1..5$.

\begin{itemize}

\item Let us assume that $h(v) = h(v') -1$. We also have $h(v) \leq
h_{\text{min}}(s_k) + 2 \leqslant (h_{\text{max}}(s_k) - 3) + 2 \leqslant
h_{\text{max}}(s_k) - 1$ so that $h(v)$ would be less than the height of any
of its immediate neighbours, \textit{i.e.} it would be a local minimum, and
thus $\Delta h (v)$ could not be~0. 

\item Let us now assume $h(v) = h(v') + 1$. We also have $h(v) \geqslant
h_{\text{max}}(s_k) - 2 \geqslant (h_{\text{min}}(s_k) + 3) - 2 =
h_{\text{min}}(s_k) + 1$ so $h(v)$ would be more than the height of any of its
immediate neighbours and could be down-flipped.

\end{itemize}

Either case brings a contradiction so at least one of the $s_k$ is solid; let
us note it $v''$. If $|h(v) - h(v'')|=2$ then, as above, we are done.
Otherwise, we have shown that $v$ has two immediate neighbours which are solid
and linked to $v$ by an arc in $T$. The same result applies to $v'$ and $v''$
so that by induction, since the number of vertices in $\cal D$ is finite, one
obtains a fracture line of $\cal D$. Since the domain delimited by this line
contains at least a triangle, it contains at least a lozenge by
Lemma~\ref{lemma:triangle:lozenge}.~\mbox{}\hfill$\square$

\medskip

We have now completed our preliminary study of fracture lines. Before proving
Theorem~\ref{theo:fracture:zones}, we need to know a bit more about the
sub-domains containing no fracture lines. There are two cases, according to
whether or not all the vertices of the sub-domain are on its contour path. We
start with the second case.

\begin{definition}
	A sub-domain $\cal D'$ of ${\cal D}$ is \emph{pre-fertile} if:

\begin{itemize}

\item it is connected and tileable;

\item at least one of its vertices is not on its contour path~$P$;

\item $\Delta h (v) \neq 0$ for any such vertex;

\item two vertices of $P$ are at distance~1 in the triangular grid if and only
if they are neighbours in~$P$.

\end{itemize}

	\worklabel{def:pre:fertile:domain}
	\end{definition}

A sub-domain is pre-fertile if all its inner points can be flipped at least
once. The union of two pre-fertile zones need not be a pre-fertile zone.
Consider for instance Figure~\ref{fig:pre:fertile:zones} where ${\cal D}_1$ is
the left-hand side hexagon and ${\cal D}_2$ is the right-hand side one, so
that $A$ belongs to the contour path of both. One sees that $A$ cannot be
involved in any flip, whence $\Delta h(A) = 0$; and yet it is not on the
contour path of ${\cal D}_1 \cup {\cal D}_2$. For this reason we give the
following definition (both ${\cal D}_1$ and ${\cal D}_2$ are examples of such
zones):

\fig{
	\epsfig{file=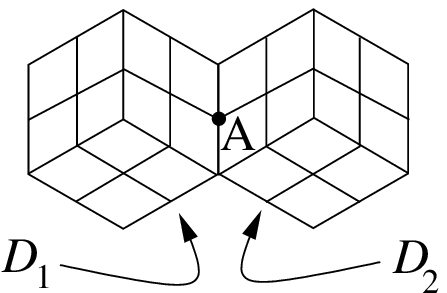, width=2cm}
	\caption{The union of two pre-fertile zones need not be pre-fertile
		\newline
		\worklabel{fig:pre:fertile:zones}}
}

\begin{definition}[Fertile zone]
	Let ${\cal D}_1$ and ${\cal D}_2$ be two pre-fertile sub-domains
of~${\cal D}$. Their union is \emph{fertile} if ${\cal D}_1 \cup {\cal D}_2$
is itself pre-fertile. A sub-domain of ${\cal D}$ is \emph{fertile} if it is
pre-fertile and a maximal element for fertile union.

	\worklabel{def:fertile:union}
	\end{definition}

\begin{proposition}
	Let $F$ be a pre-fertile sub-domain of $\cal D$. It is fertile if and
only if its contour path is a fracture line of $\cal D$.

	\worklabel{prop:fertile:contour}
	\end{proposition}

\noindent \textit{Proof} \quad We first prove that the contour path of a
fertile sub-domain~$F$ is a fracture line in~$\cal D$. Let $a$ be any vertex
of the contour path of $F$. If~$a$ belongs to the contour path of~$\cal D$, it
is solid by Lemma~\ref{lemma:boundary:height}. Otherwise, $a$ is in the
interior of~$\cal D$; let $H_a$ be the hexagon (in the triangular grid) of
side~1 whose center is~$a$. Since $F$ is fertile it is a maximal element for
fertile union so that $F \cup H_a$ is not pre-fertile, which means that at
least one of its inner vertices is solid. Since none of~$F$'s interior is, $a$
must be this vertex. In other words, all the vertices on the contour path
of~$F$ are solid and we are done.

Conversely, we prove that a pre-fertile sub-domain whose contour path is a
fracture line in~$\cal D$ is fertile. Let~$F'$ be any other pre-fertile
sub-domain of~$\cal D$ such that $F \cap F' = \emptyset$ and $F \cup F'$ is
connected. Let~$a$ be a vertex that belongs to both $F$ and $F'$ (it is on the
boundary path of both). If $a$ is not an inner point of $F \cup F'$ for all
such~$a$, then this union is not pre-fertile (see the last part of
Definition~\ref{def:pre:fertile:domain}). Otherwise, $a$ is in the interior
of~$F \cup F'$, and it is also solid since it belongs to the boundary path
of~$F$, whence $F \cup F'$ is not pre-fertile, from which we deduce that~$F$
is indeed a maximal element for fertile union and therefore a fertile
sub-domain.~\mbox{}\hfill$\square$

\begin{proposition}
	Let $F$ be a sub-domain of $\cal D$, delimited by a fracture line~$L$,
such that at least one of its vertices is not on its contour path. Then $F$ is
fertile if and only if it contains no non-trivial fracture line of $\cal D$
other than $L$.

	\worklabel{prop:fracture:lines}
	\end{proposition}

\noindent \textit{Proof} \quad If $F$ is fertile then $\Delta h(v) \neq 0$ for
any of its inner vertices, so that it can contain no fracture line other than
$L$. Conversely, recall that a fracture line can be made of only one point:
therefore if $F$ contains no fracture line other than $L$ then $\Delta h(v)
\neq 0$ for any of its inner vertices. Moreover $F$ is tileable (see
Corollary~\ref{coro:fracture:line:tileable}), connected, and two vertices
of~$L$ that are not neighbours in~$L$ cannot be at distance~1 in the
triangular grid (this would define a non-trivial fracture line) so it is
pre-fertile, and since its contour path is a fracture line of~$\cal D$ it is a
fertile zone by Proposition~\ref{prop:fertile:contour}.~\mbox{}\hfill$\square$

\medskip

We have completed our study of fertile zones; we now turn to sub-domains whose
vertices are all on the contour path.

\begin{lemma}
	Let $L$ be a fracture line of ${\cal D}$ enclosing a domain $F$, not
limited to a single point, such that all the vertices of $F$ belong to $L$.
Then either $F$ is a lozenge or it contains a fracture line which defines a
lozenge.

	\worklabel{lemma:fracture:line:lozenge}
	\end{lemma}

\noindent \textit{Proof} \quad Since $F$ is not limited to a single point, it
contains at least a triangle. This triangle belongs to a lozenge $x$ in a
tiling $T$ of $\cal D$. The fracture line $L$ is a $T$-valid path, therefore
$x$ is embedded in $F$ by lemma~\ref{lemma:triangle:lozenge}. By hypothesis
all its vertices ($v_1$, $v_2$, $v_3$ and $v_4$) are on a fracture line of
$\cal D$, so that $\Delta h(v_k)=0$ for $k=1..4$. Thus the closed path
$v_1$-$v_2$-$v_3$-$v_4$-$v_1$ is a fracture line in $\cal D$ and we are
done.~\mbox{}\hfill$\square$

\medskip

We need one last definition, which corresponds to the sub-domains that are
obtained by recursively looking for fracture lines, and therefore used by
Algorithm~\ref{algo:fracture:lines} (see also Figures~\ref{fig:example}~(b)  
and~\ref{fig:example:zones}~(b)):

\begin{definition}
	A \emph{fracture zone} in a tileable domain $\cal D$ is a sub-domain
of $\cal D$ not reduced to a single point, delimited by a fracture line and
containing no other fracture line.

	\worklabel{def:fracture:zone}
	\end{definition}

\begin{theorem}[Fracture theorem]
	Let $\cal D$ be any tileable domain such that the results of
Section~\ref{flips:lattices} hold. A fracture zone in $\cal D$ is either a
fertile zone or a lozenge.

	\worklabel{theo:fracture:zones}
	\end{theorem}

\noindent \textit{Proof} Let $F$ be a fracture zone of $\cal D$. If all its
vertices belong to its contour path, then it is a lozenge by
Lemma~\ref{lemma:fracture:line:lozenge}. Otherwise $\Delta h(v) \neq 0$ for
any inner vertex~$v$ of~$F$ by Proposition~\ref{prop:trivial:fracture:lines}
since~$F$ is a fracture zone, so that~$F$ is fertile by Proposition
\ref{prop:fracture:lines}.~\mbox{}\hfill$\square$

\medskip

We now have the tools to extract the fertile zones from a domain:

\begin{algorithm}[Fracture algorithm]
\worklabel{algo:fracture:lines}
\end{algorithm}

\begin{itemize}

\item \textbf{Input}: A tileable domain $\cal D$.

\item \textbf{Output}: The fertile zones of $\cal D$.

\item \textbf{Step 1}: Use Thurston's Algorithm~\ref{algo:Thurston} to build
the maximal and minimal tilings of $\cal D$. Store the solid points of $\cal
D$ in a list $S$.

\item \textbf{Step 2}: Use $S$ and the minimal tiling of $\cal D$ to build all
the fracture lines except the trivial one; store them in a list $L$.

\item \textbf{Step 3}: Remove from $L$ all the fracture lines consisting of
exactly 4 vertices.

\item \textbf{Step 4}: Return $L$. 

\end{itemize}

\noindent \emph{Proof of the algorithm} \quad Everything has already been
proved in Theorem~\ref{theo:fracture:zones}.

\medskip

\noindent \emph{Complexity} \quad Steps 1 and 2, which control the execution
time of the algorithm, are both linear in the number of vertices of $\cal D$.

%!!!!!!!!!!!!!!!!!!!!!!                           !!!!!!!!!!!!!!!!!!!!!!
%!!!!!!!!!!!!!!!!!!!!!!   The product algorithm   !!!!!!!!!!!!!!!!!!!!!!
%!!!!!!!!!!!!!!!!!!!!!!                           !!!!!!!!!!!!!!!!!!!!!!

\subsection{An example of the use of the Fracture Theorem}
	\worklabel{example}

We now present an example. Let us consider the finite-length closed path of
Figure~\ref{fig:example}~(a), which delimits a tileable domain (see
Figure~\ref{fig:example}~(b) for its minimal and maximal tilings). Using
fracture lines, one can disconnect the central part, whose tilings give rise to
a trivial lattice (see Figure~\ref{fig:example:zones}~(a)). Thus only the
fertile zones remain (see Figure~\ref{fig:example:zones}~(b)).

\fig{
	\begin{tabular}{p{3cm}@{\hskip2cm}p{6.5cm}}

	\mbox{}\hfill
	\epsfig{file = 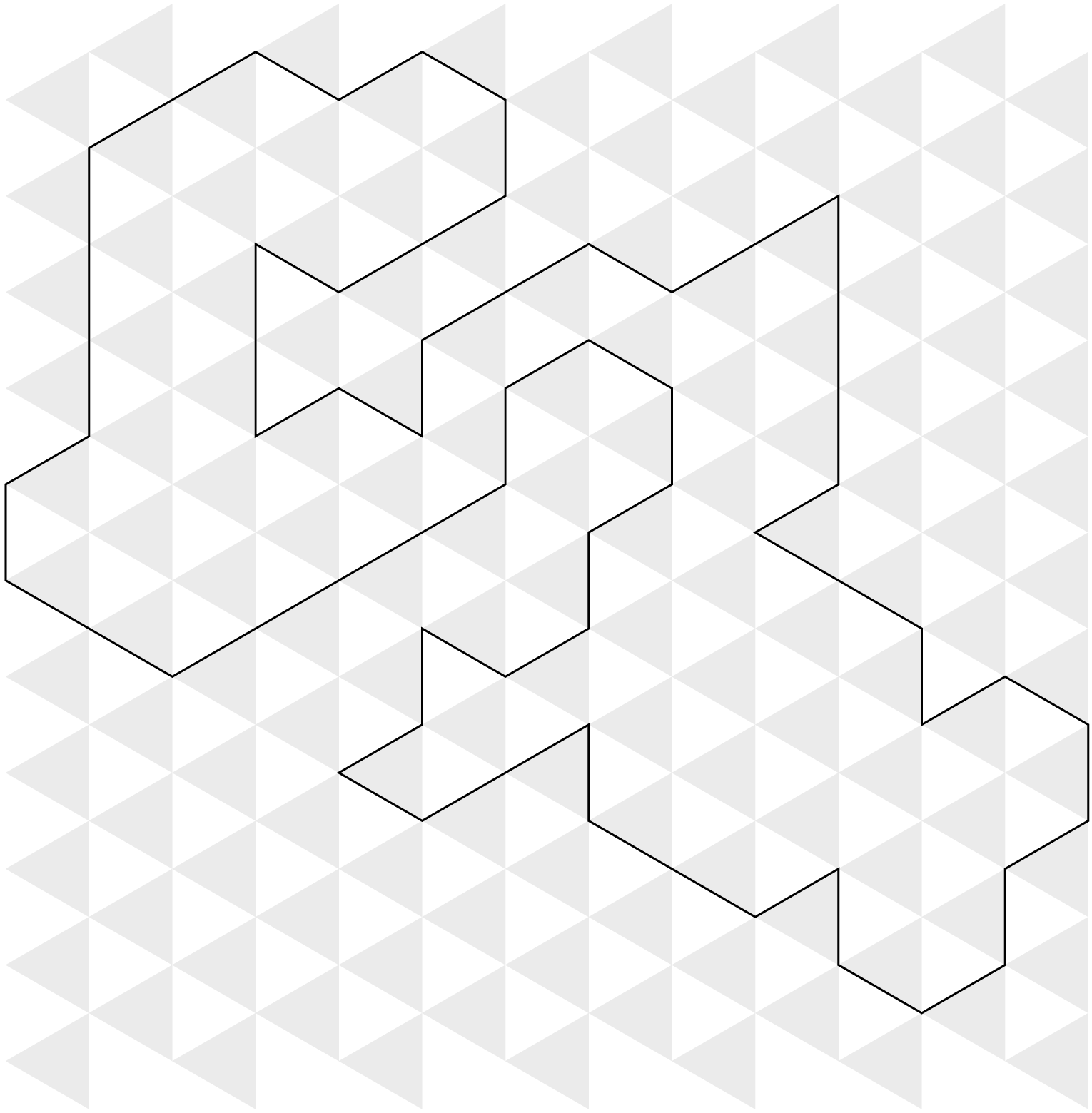, width = 2cm}
	\hfill\mbox{}
	&
	\mbox{}\hfill
	\epsfig{file = 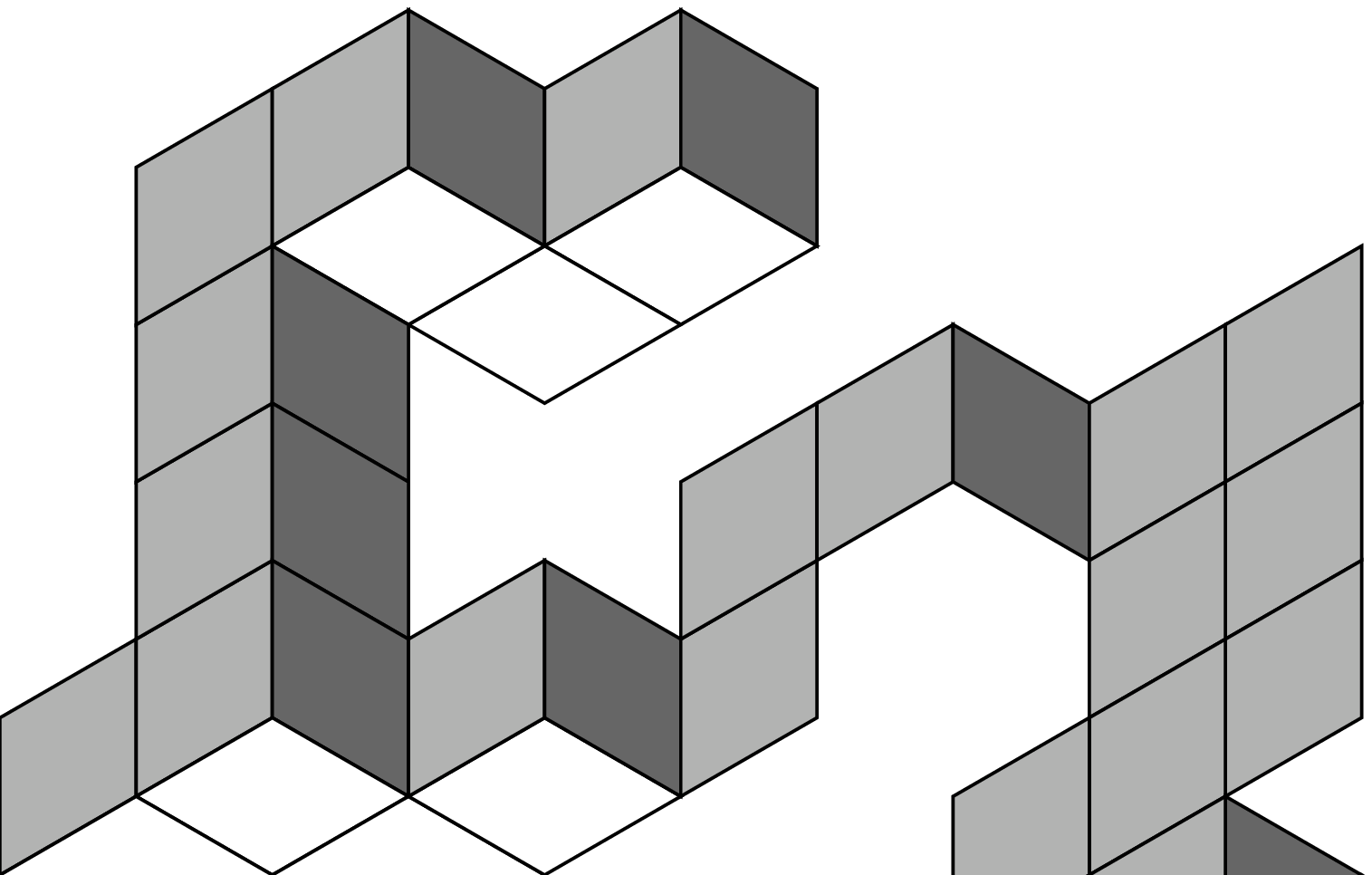, width=2cm}
	\hskip5mm
	\epsfig{file = 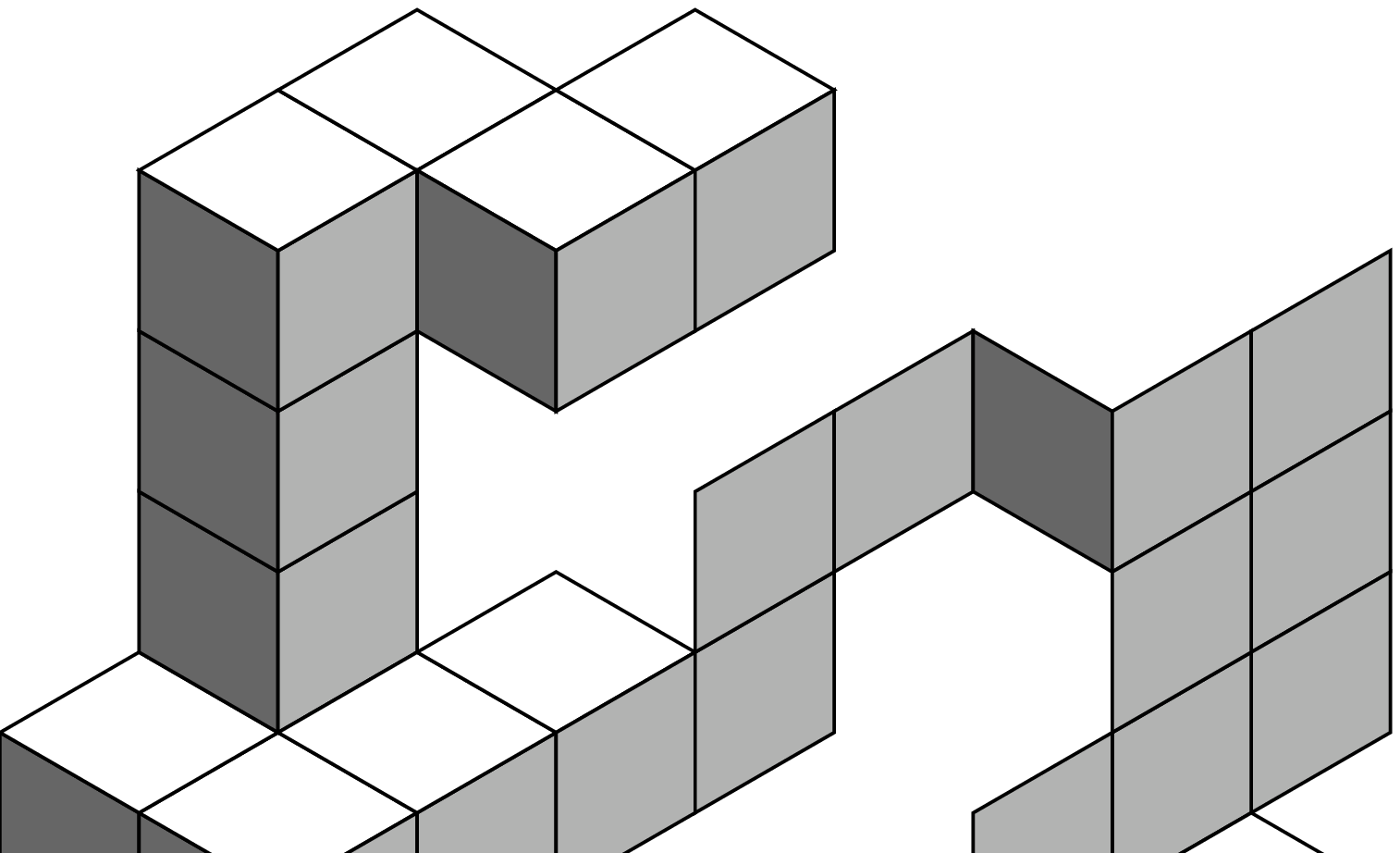, width=2cm}
	\hfill\mbox{}
	\\[1mm]

	\mbox{}\hfill
	(a) 	\begin{minipage}[t]{1.5cm}
		A generic domain
		\end{minipage}
	\hfill\mbox{}
	&
	\mbox{}\hfill
	(b) 	\begin{minipage}[t]{4cm}
		The minimal and maximal tilings of the domain
		\end{minipage}
	\hfill\mbox{}
	\\

	\end{tabular}
	\caption{A domain and its extreme tilings 
		\worklabel{fig:example}
	}
}

\fig{
	\begin{tabular}{p{3.5cm}@{\hskip1.5cm}p{6.5cm}}

	\mbox{}\hfill
	\epsfig{file = 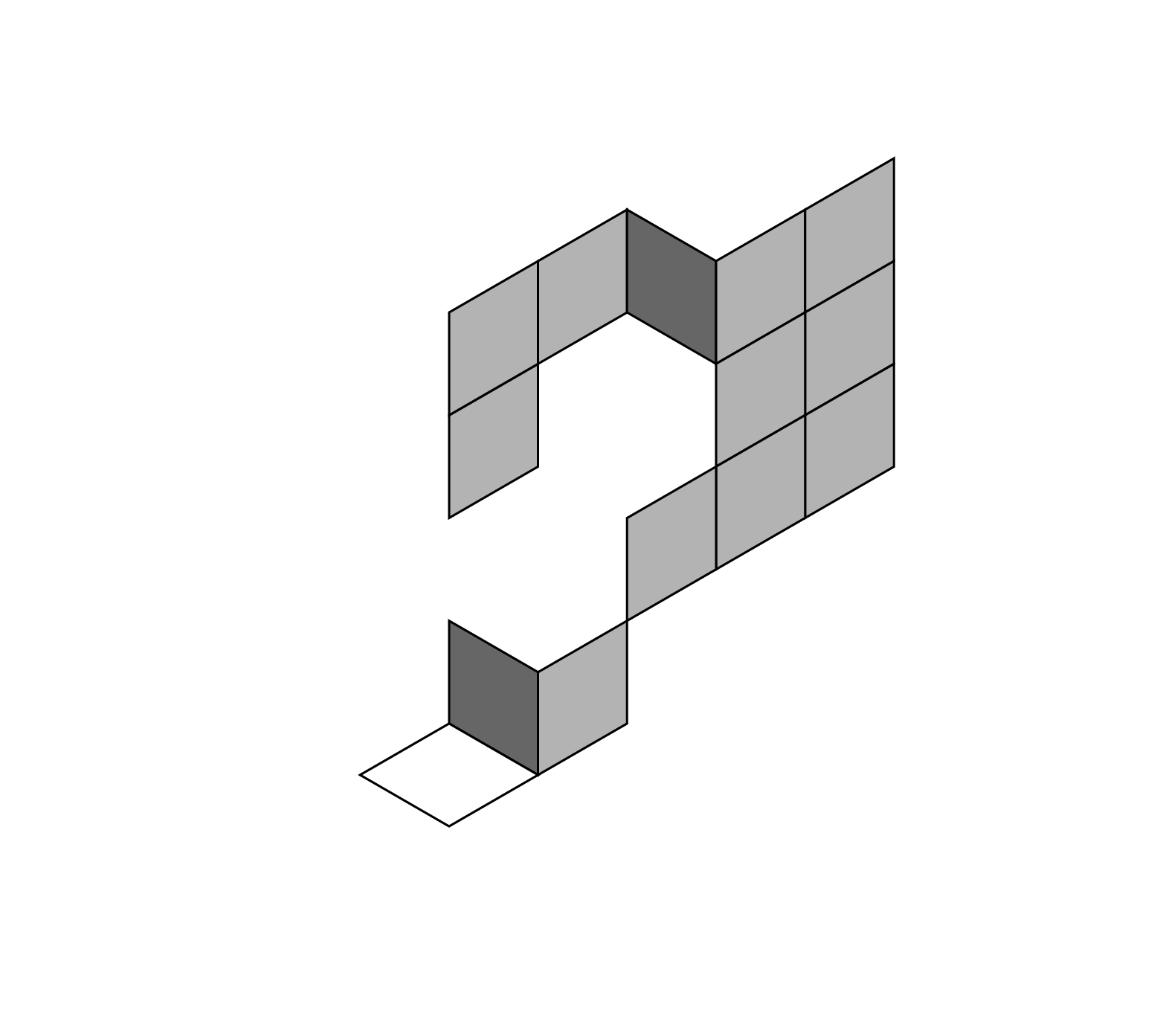, width=2cm}
	\hfill\mbox{}
	&
	\mbox{}\hfill
	\epsfig{file = 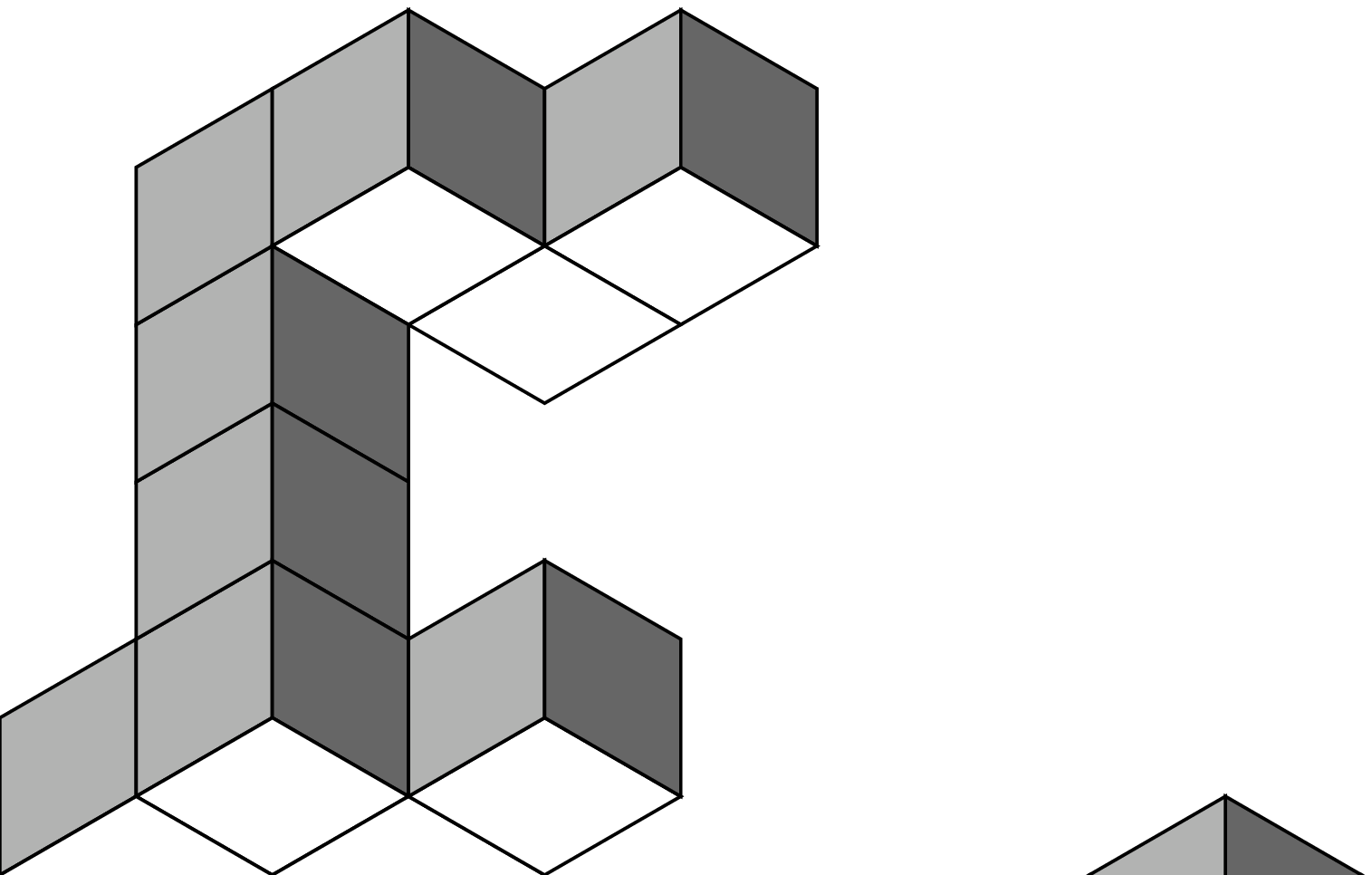, width=2cm}
	\hskip10mm
	\epsfig{file = 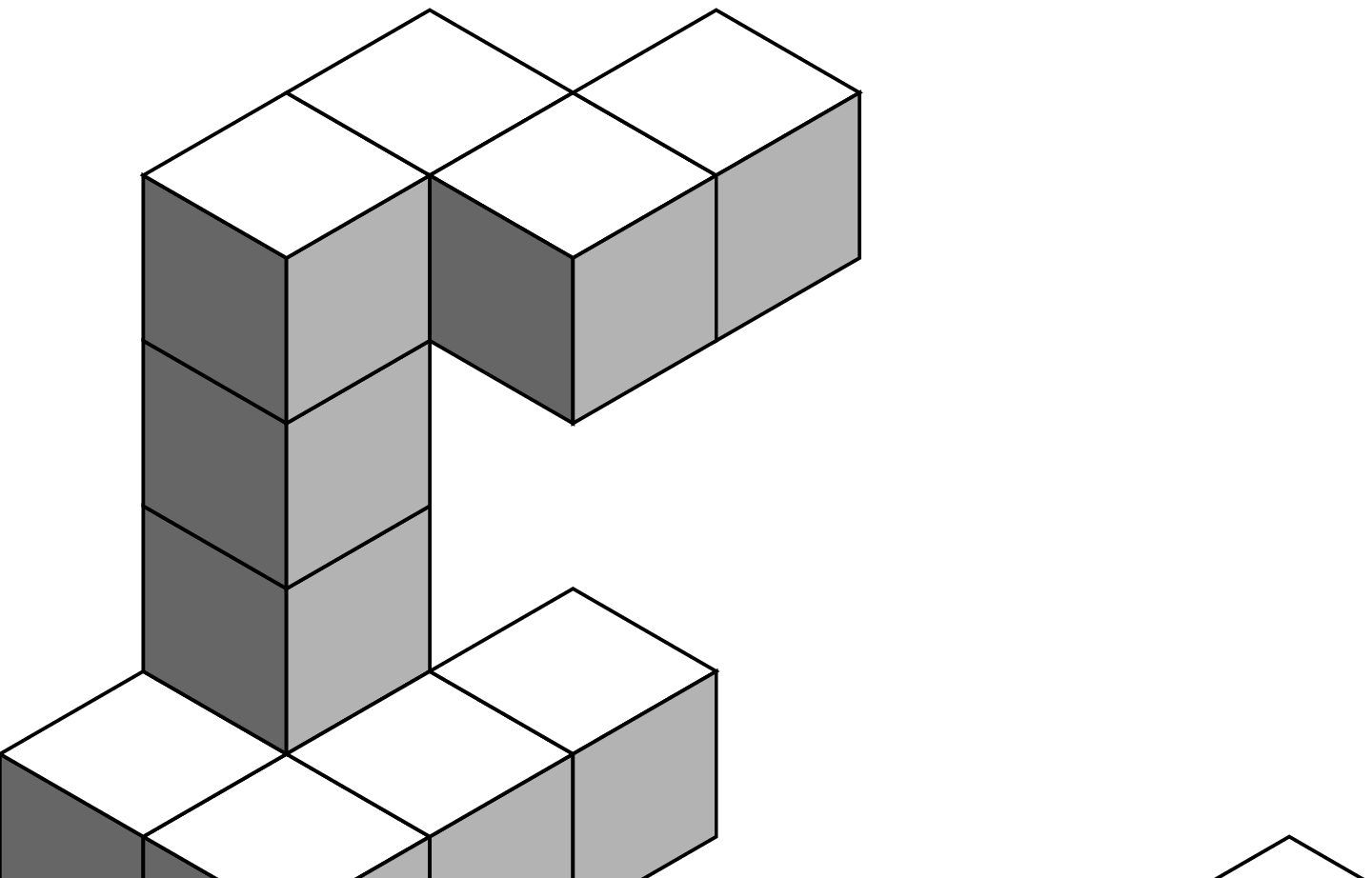, width=2cm}
	\hfill\mbox{}
	\\[1mm]

	\mbox{}\hfill
	(a) 	\begin{minipage}[t]{2.9cm}
		The maximal trivial fracture zone
		\end{minipage}
	\hfill\mbox{}
	&
	\mbox{}\hfill
	(b) 	\begin{minipage}[t]{5cm}
		The fertile zones in the minimal and maximal tilings
		\end{minipage}
	\hfill\mbox{}
	\\

	\end{tabular}
	\caption{The fracture zones
		\worklabel{fig:example:zones}
	}
}

%!!!!!!!!!!!!!!!                                         !!!!!!!!!!!!!!!
%!!!!!!!!!!!!!!!   Application du théorème de fracture   !!!!!!!!!!!!!!!
%!!!!!!!!!!!!!!!                                         !!!!!!!!!!!!!!!

\subsection{An application of the Fracture Theorem}
	\worklabel{subsection:application:FT}

Since fertile zones are delimited by fracture lines, flips done inside one of
them cannot have any effect on the other fertile zones: therefore one can
study independently the tilings of these zones. Since the tilings of $\cal D$
will be obtained by ``gluing'' tilings from these disjoint zones, their
lattice can be obtained by product. Consider Figure~\ref{fig:example:lattice},
in which the domain contains two fertile zones. We will see later (see
Section~\ref{section:algo:main}) an algorithm to generate the lattice of the
tilings of such a zone; on this example each fertile zone is a pseudo-hexagon,
so its tilings can be found using Algorithm~\ref{algo:tilings}. We build the
product lattice of these lattices by connecting every element in the left-hand
side lattice with every element in the right-hand side lattice, which finally
gives the lattice of the tilings of the initial domain.

\fig{
	\begin{tabular}{p{2cm}@{\hskip5mm}
			p{2.5cm}@{\hskip5mm}
			p{2.5cm}@{\hskip5mm}
			p{3.2cm}
			}

	\mbox{}\hfill
	\epsfig{file=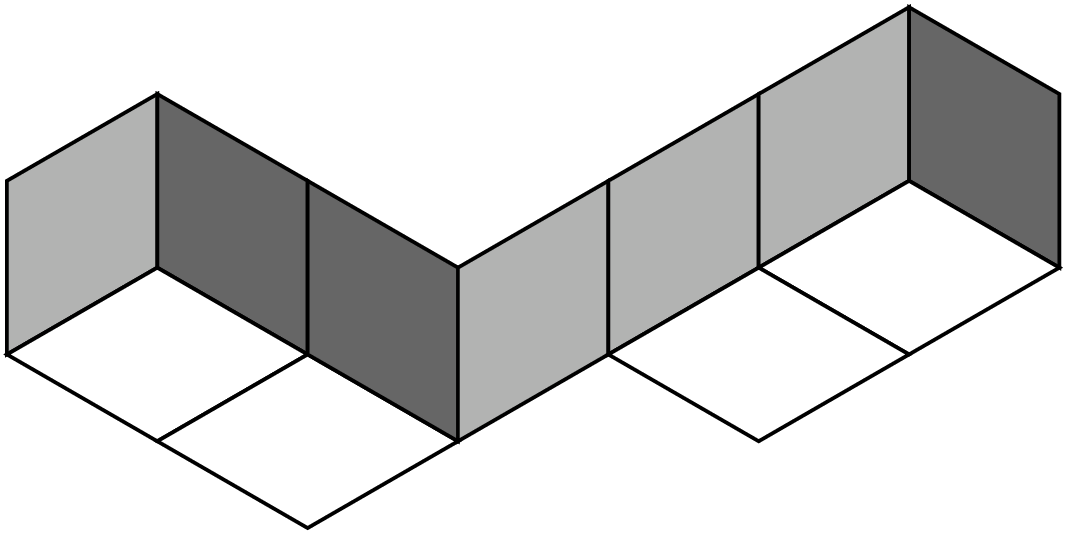, width = 1.2cm}
	\hfill\mbox{}
	&
	\mbox{}\hfill
	\epsfig{file=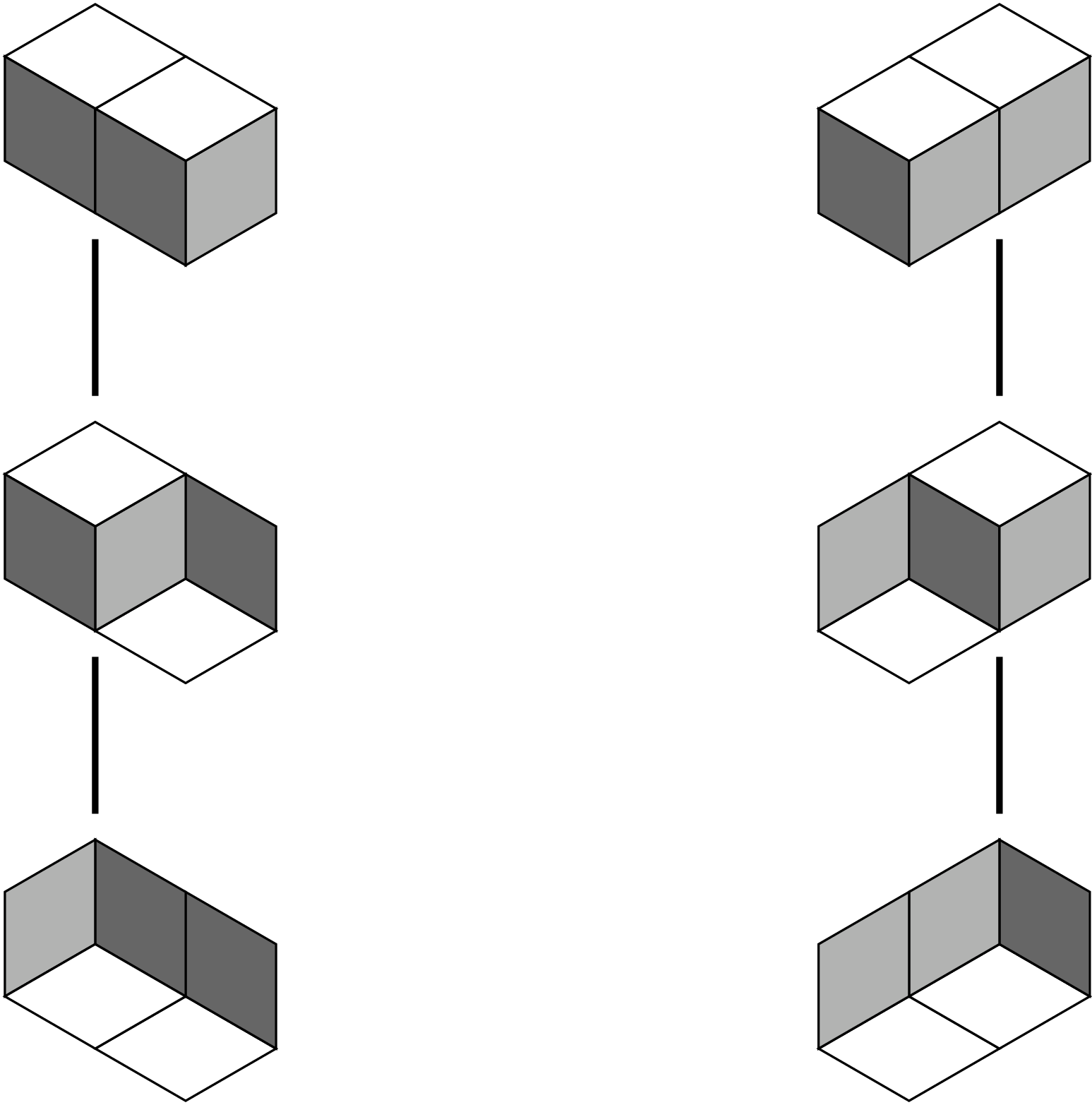, width = 2cm}
	\hfill\mbox{}
	&
	\mbox{}\hfill
	\epsfig{file=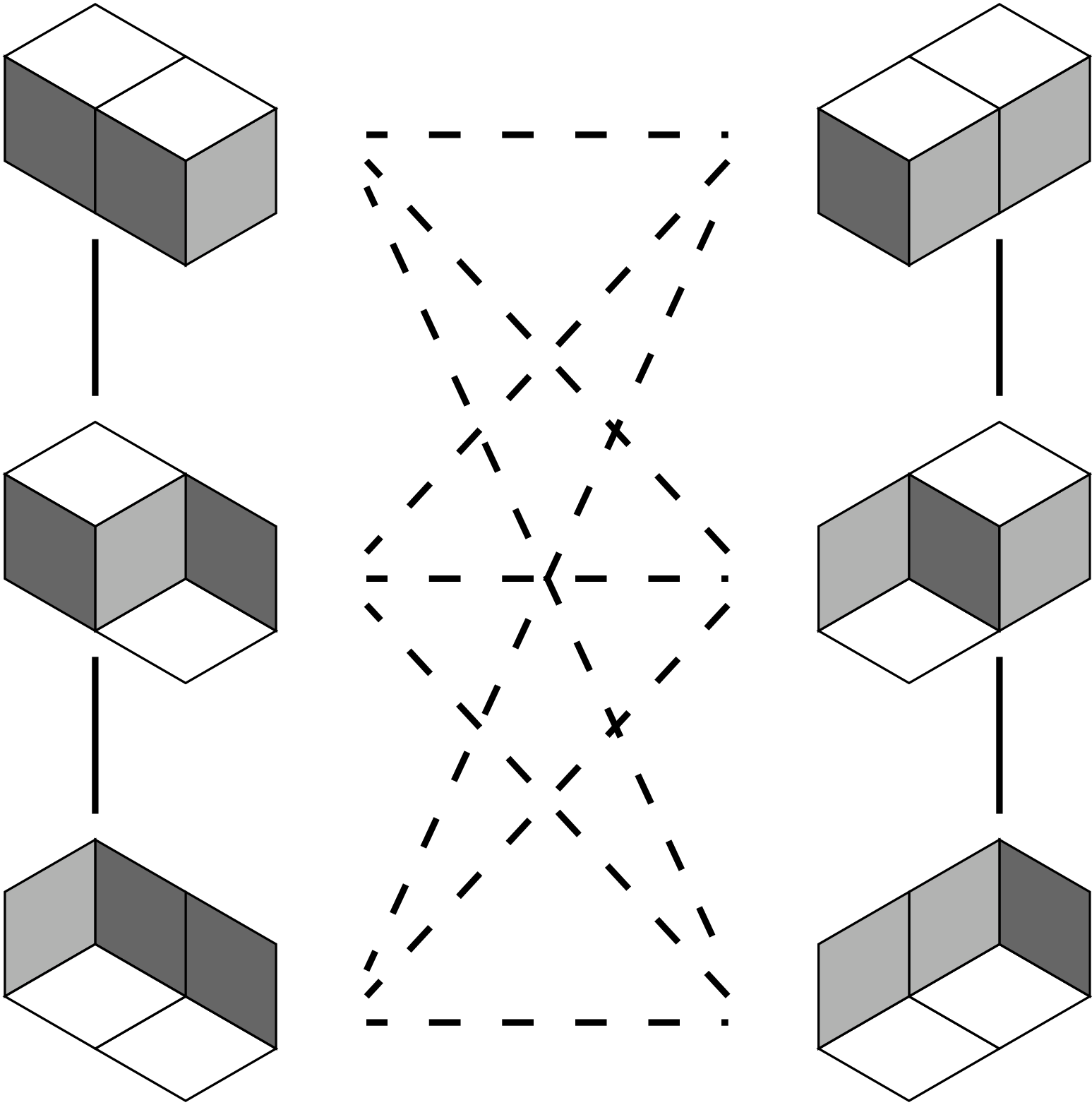, width = 2cm}
	\hfill\mbox{}
	&
	\mbox{}\hfill
	\epsfig{file=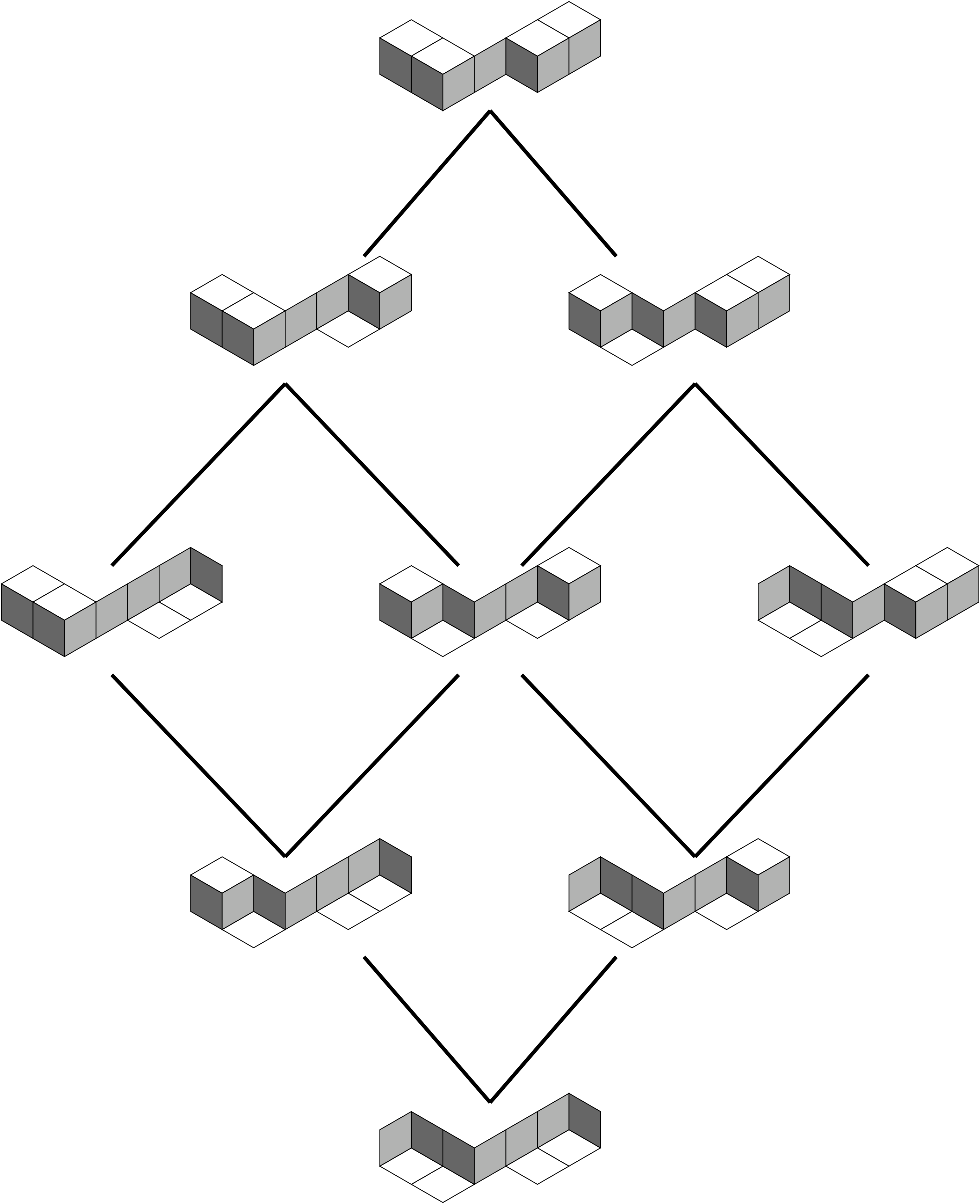, width = 3cm}
	\hfill\mbox{}
	\\[1mm]

	The minimal tiling of a domain
	&
	The lattices of\newline the tilings of the fertile zones
	&
	How to compute a product lattice
	&
	The complete lattice of the tilings of the domain
	\\

	\end{tabular}
	\caption{The main steps in computing the lattice 
		\worklabel{fig:example:lattice}
	}
}

%!!!!!!!!!!!!!!!!!!!!                              !!!!!!!!!!!!!!!!!!!!
%!!!!!!!!!!!!!!!!!!!!   Seeds and their fillings   !!!!!!!!!!!!!!!!!!!!
%!!!!!!!!!!!!!!!!!!!!                              !!!!!!!!!!!!!!!!!!!!

\subsection{Seeds and their fillings}
\worklabel{section:seeds}

We have seen in Section~\ref{fracture:lines} that fertile zones play a major
role in the generation of all the tilings of a generic domain $\cal D$.  We
have also seen in Section~\ref{subsection:application:FT} that we can restrict
our study to these zones and later compute a product lattice. We now show that
fertile zones can be decomposed into collections of pseudo-hexagons (see
Figure~\ref{fig:ranges}~(a) for a start). The correspondence between tilings
and compact piles of cubes allows us to use geometric terms; in other words,
an up-flip can be viewed as adding a cube.

\begin{definition}[Seed]
	A \emph{seed} is the minimal tiling of a hexagon of side 1. 
	\worklabel{def:seed}
	\end{definition}

Note that a pseudo-hexagon (see Definition~\ref{def:ph}) contains exactly one
seed. 

\smallskip

We now investigate the immediate properties of seeds.

\begin{lemma}
	The union of two pseudo-hexagons sharing the same seed is again a
pseudo-hexagon.
	\worklabel{lemma:ph:union}
	\end{lemma}

\noindent \textit{Proof} \quad Using Definition~\ref{def:ph}, it is enough to
show that the union of two compact piles of cubes, aligned on the same set of
axes, is itself a compact pile of cubes. But this immediately follows from
Definition~\ref{def:compact:pile}.~\mbox{}\hfill$\square$

\begin{figure}[ht]
\begin{center}
	\begin{tabular}{p{5.5cm}@{\hskip7mm}p{5.5cm}}

	\mbox{}\hfill
	\epsfig{file=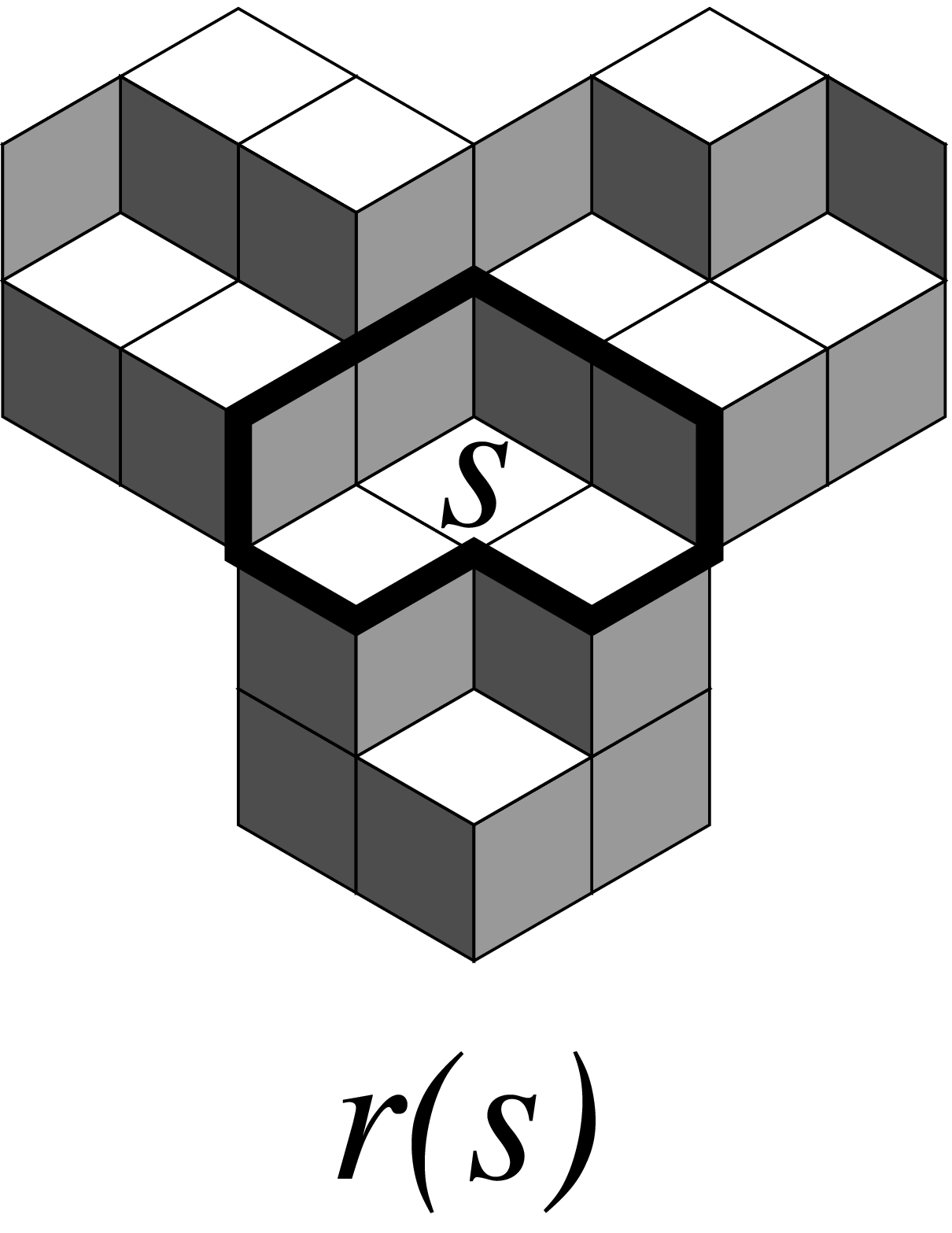, width=1.7cm, clip=}
	\hskip5mm
	\epsfig{file=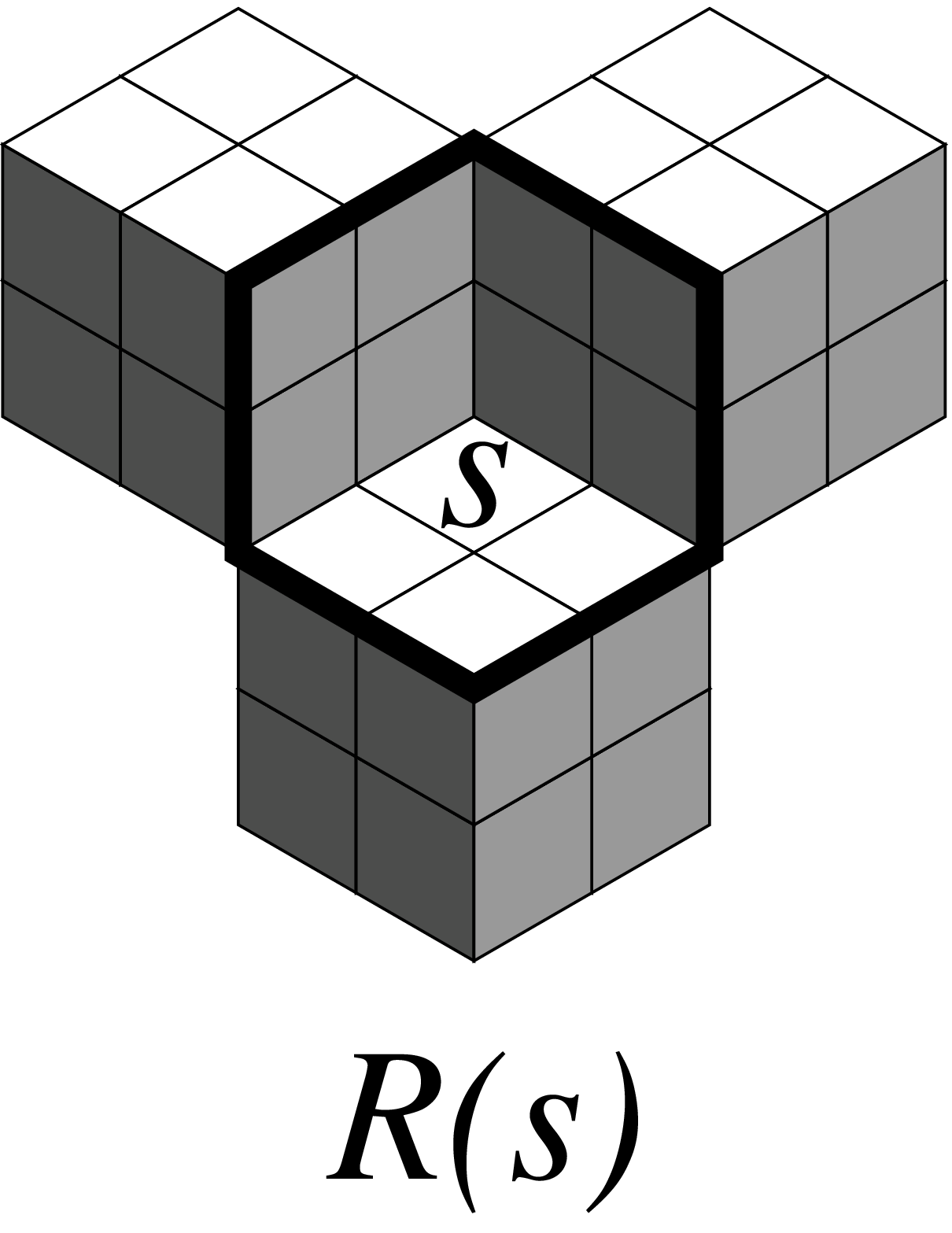, width=1.7cm, clip=}
	\hfill\mbox{}
	&
	\mbox{}\hfill
	\epsfig{file=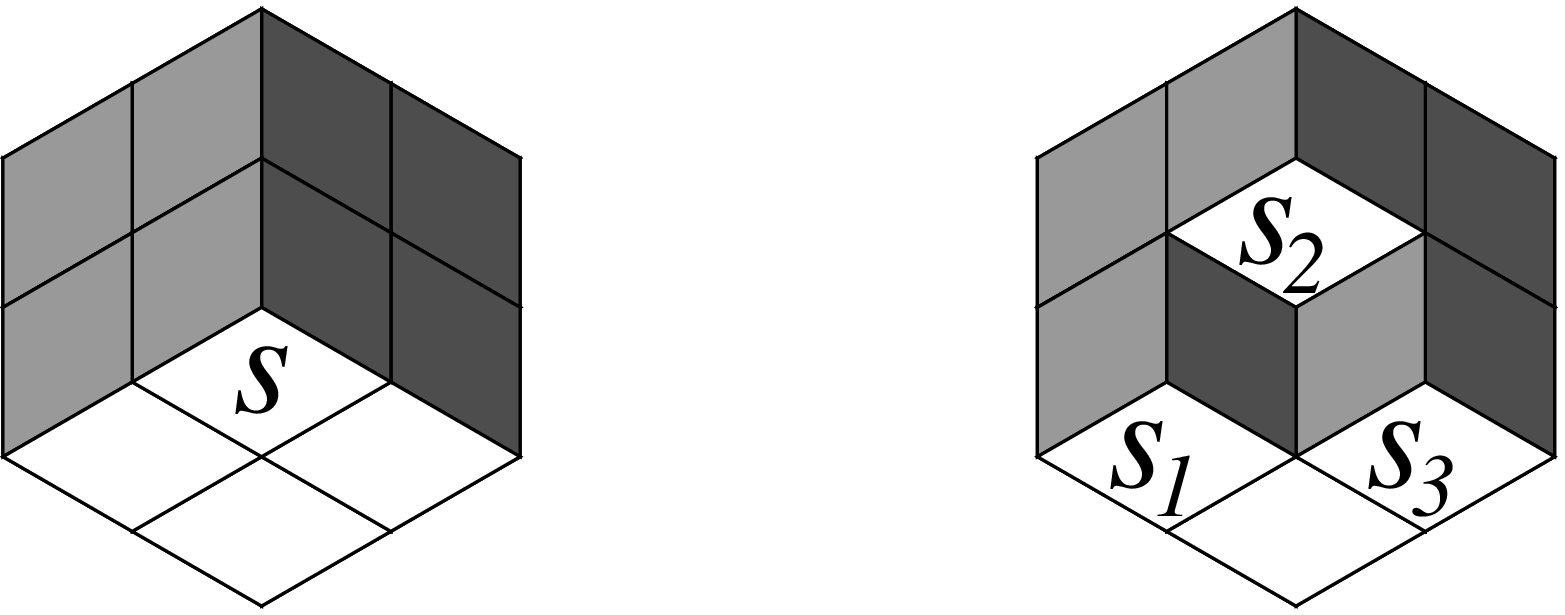, width=4cm, clip=}
	\hfill\mbox{}
	\\[1mm]

	\mbox{}\hfill
	(a) 	\begin{minipage}[t]{4.9cm}
		A range and the maximal range of a seed
		\end{minipage}
	\hfill\mbox{}
	&
	\mbox{}\hfill
	(b) 	\begin{minipage}[t]{4.9cm}
		Covering a seed with a cube potentially creates three new 
		seeds
		\end{minipage}
	\hfill\mbox{}
	\\

	\end{tabular}
	\caption{Seeds and ranges \worklabel{fig:ranges}}
\end{center}
\end{figure}

\begin{definition}
	The \emph{range} $r(s)$ of a seed $s$ in a tiling~$T$ of ${\cal D}$ is
the union of the pseudo-hexagons that contain it. The \emph{maximum range}
$R(s)$ of $s$ is the union of its ranges when~$T$ ranges among all the tilings
of~$\cal D$.
	\worklabel{def:ranges}
	\end{definition}

The set of pseudo-hexagons surrounding a seed $s$ is not empty since it
contains at least $s$ itself. Thus the range of $s$ is well-defined according
to Lemma~\ref{lemma:ph:union}.

\begin{definition}
	A tiling of $r(s)$ $($resp. the maximal, minimal tiling of $R(s)$$)$
will be called a \emph{filling} of $s$ $($resp. the \emph{maximal filling},
\emph{minimal filling} of $s$$)$. We denote by $\text{\rm Max}(s)$ and
$\text{\rm Min}(s)$ the maximal and minimal fillings of $s$; by extension,
$\text{\rm Max}({\cal D})$ and $\text{\rm Min}({\cal D})$ are the maximal and
minimal tilings of ${\cal D}$.

	\worklabel{def:fillings}
	\end{definition}

The intuitive idea behind the term ``filling'' is that one gradually fills
$r(s)$ with (compact piles of) cubes.  Note that a filling of a
seed is a compact pile of cubes and that the maximal filling of $s$ is
obtained by doing all the possible up-flips in~$r(s)$.

Adding a cube to cover a seed $s$ (in other words, flipping $s$) generally
creates three new seeds $s_1$, $s_2$ and~$s_3$ as shown in
Figure~\ref{fig:ranges}~(b).  We thus see that there is a natural partial
order on the ranges of seeds, defined by inclusion. In our example, $r(s_i)
\subset r(s)$ for $i \in \{1,2,3\}$. We will be mostly interested in the
maximal elements for this partial order and now proceed to define this more
precisely.

\begin{definition}
	A seed $s$ is a \emph{child} of a seed $s'$ if $R(s)$ is embedded in
$R(s')$. A seed is a \emph{proper seed} if it is the child of only itself.

	\worklabel{def:proper:seeds}
	\end{definition}

\begin{proposition}
	Every cube in a tiling of ${\cal D}$ belongs to a filling of a proper
seed.

	\worklabel{proposition:cubes:ps}
	\end{proposition}

\noindent \textit{Proof} \quad Every cube is the maximal tiling of a hexagon
of side~1 and it is added when the minimal tiling of this hexagon is flipped.
But this minimal tiling is a seed, whose range is included in the range of a
proper seed. Therefore the cube belongs to a filling of this proper
seed.~\mbox{}\hfill$\square$

\begin{proposition}
	The maximal fillings of two distinct proper seeds are disjoint.

	\worklabel{proposition:disjoint:fillings}
	\end{proposition}

\noindent \textit{Proof} \quad Let $C$ be a cube belonging to the maximal
filling of two proper seeds $s_1$ and $s_2$. Since a filling of a seed is a
compact pile of cubes, there is an uninterrupted line of cubes connecting $C$
to the base planes defined by $s_1$ and $s_2$. Since these planes vary two by
two by a translation, the two basis are in fact the same, so that $s_1 =
s_2$.\hfill$\square$

\begin{corollary}
	The maximal fillings of the proper seeds in a domain ${\cal D}$ form a
partition of the subset of $\mathbb{Z}^3$ defined by the maximal tiling of the
fertile zones of ${\cal D}$.

	\worklabel{coro:partition}
	\end{corollary}

We know that the set of the tilings of ${\cal D}$ is connected by flips;
therefore maximally filling all the seeds that appear in $\text{Min}({\cal
D})$ potentially creates new flippable zones: that is, new seeds and new
proper seeds (see Figure~\ref{fig:seeds}). To clarify the situation, we give
the following definition:

\begin{definition}
	A proper seed is of \emph{order 0} if it appears in the minimal tiling
of ${\cal D}$. It is of \emph{order $n+1$} if it is a proper seed in the
tiling of~$\cal D$ obtained when all the proper seeds of order $k \in
\{1,\ldots,n\}$ are maximally filled with cubes.

	\worklabel{def:seed:order}
	\end{definition}

\fig{
	\begin{tabular}{p{2cm}@{\hskip5mm}
			p{1.8cm}@{\hskip5mm}
			p{3.2cm}@{\hskip5mm}
			p{3cm}
			}

	\mbox{}\hfill
	\epsfig{file=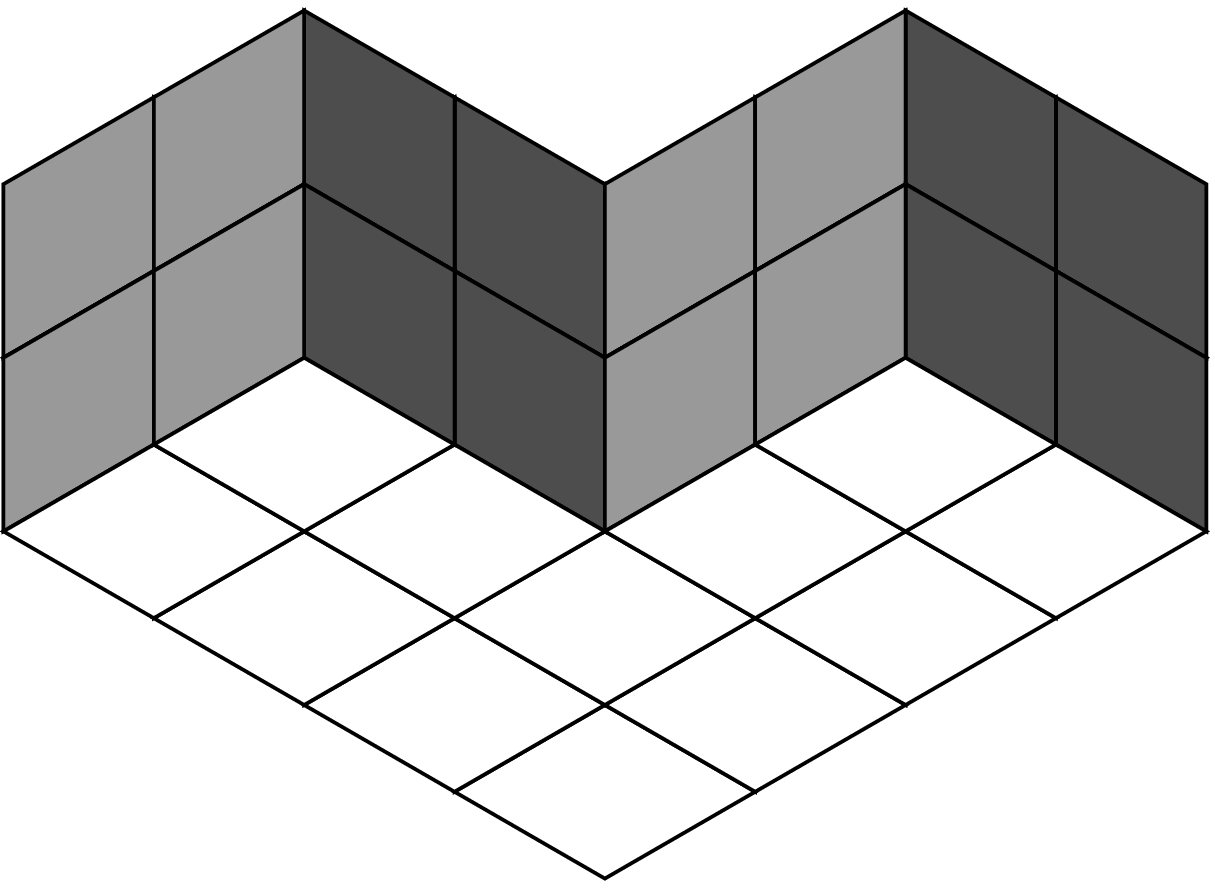, width=1.5cm}
	\hfill\mbox{}
	&
	\mbox{}\hfill
	\epsfig{file=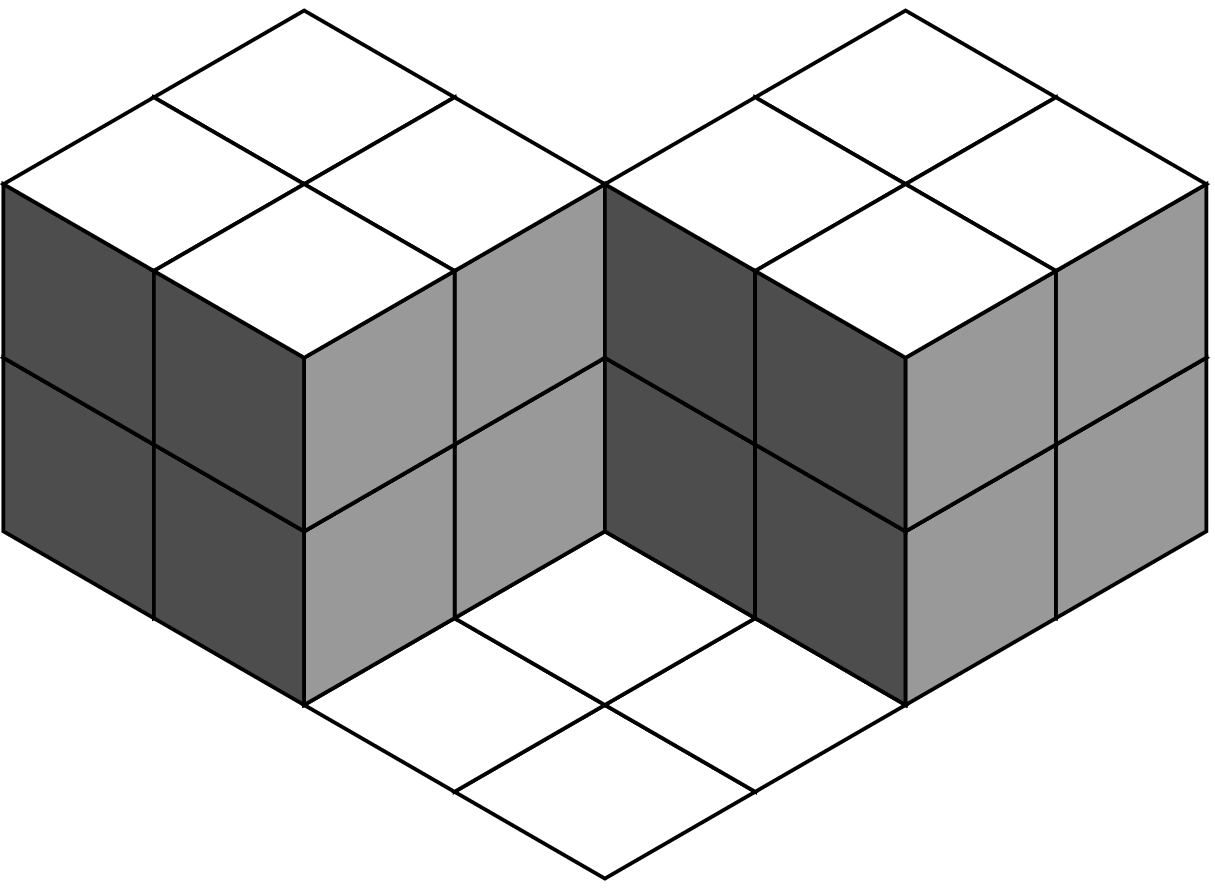, width=1.5cm}
	\hfill\mbox{}
	&
	\mbox{}\hfill
	\epsfig{file=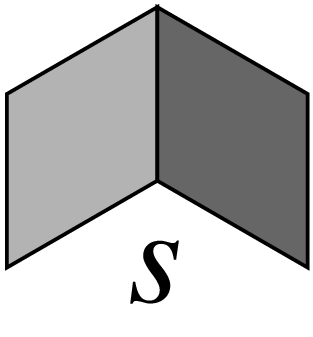, width=1cm}
	\hfill\mbox{}
	&
	\mbox{}\hfill
	\epsfig{file=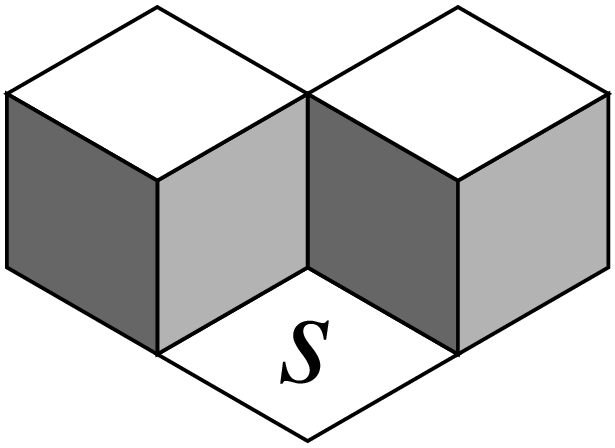, width=2cm}
	\hfill\mbox{}
	\\[1mm]

	Two seeds of order~0
	&
	One seed of order~1
	&
	Two lozenges around a seed of order 1
	&
	The neighbourhood of a seed of order 1
	\\

	\end{tabular}
	\caption{Seeds 
		\worklabel{fig:seeds}
	}
}

%!!!!!!!!!!!!!!!!!!!!!!!!                      !!!!!!!!!!!!!!!!!!!!!!!!
%!!!!!!!!!!!!!!!!!!!!!!!!   Seeds of order 1   !!!!!!!!!!!!!!!!!!!!!!!!
%!!!!!!!!!!!!!!!!!!!!!!!!                      !!!!!!!!!!!!!!!!!!!!!!!!

\subsection{Seeds of order 1}

Let us examine the basic properties of a seed $s$ of order 1. Remember that
this is the minimum tiling of a hexagon and that it does not appear in the
minimal tiling of~$\cal D$. Thus at most two of the three lozenges that make
up $s$ can come from~$\text{Min}({\cal D})$.

Let us first assume that exactly two of the three lozenges in $s$ come from
$\text{Min}({\cal D})$, and exactly one from a proper seed~$t$ of order~0.
Then $s$ appears while filling the range of $t$, and it is a flippable zone:
therefore it is not a proper seed and should be flipped while filling $t$; it
is not a seed of order~1.

So at least two of the lozenges in $s$ come from fillings of proper seeds of
order 0. For reasons of symmetry one can assume that they are as in
Figure~\ref{fig:seeds}. Since they come from fillings they belong to cubes, so
that the immediate neighbourhood of $s$ looks like the right-most picture in
Figure~\ref{fig:seeds}.

%!!!!!!!!!!!!!!!!!                                   !!!!!!!!!!!!!!!!!
%!!!!!!!!!!!!!!!!!   Are there seeds of any order?   !!!!!!!!!!!!!!!!!
%!!!!!!!!!!!!!!!!!                                   !!!!!!!!!!!!!!!!!

\subsection{Seeds of greater order}

For a given domain ${\cal D}$ the order of any seed is bounded since
$\text{Min}({\cal D})$ and $\text{Max}({\cal D})$ vary by a given number of
cubes and each proper seed contributes at least one cube.  Conversely, if $n$
is any integer, there exist domains in which one can find proper seeds of
order~$n$. To illustrate this, consider Figure~\ref{fig:seeds:n}. If one
maximally fills all the seeds of order~0~(b), a sub-domain closely resembling
the first can be outlined~(c); indeed, only the size differs.

Thus one can build a family $({\cal D}_n)$ of domains in which covering all
the seeds of order 0 of ${\cal D}_{n+1}$ exactly yields ${\cal D}_n$ and
therefore all the seeds which are of order $k$ in ${\cal D}_{n}$ are of order
$k+1$ in ${\cal D}_{n+1}$. Using this procedure, one can generate seeds of any
order.

\fig{
	\begin{tabular}{p{3.5cm}@{\hskip3mm}
			p{3.5cm}@{\hskip4mm}
			p{4cm}
			}

	\mbox{}\hfill
	\epsfig{file=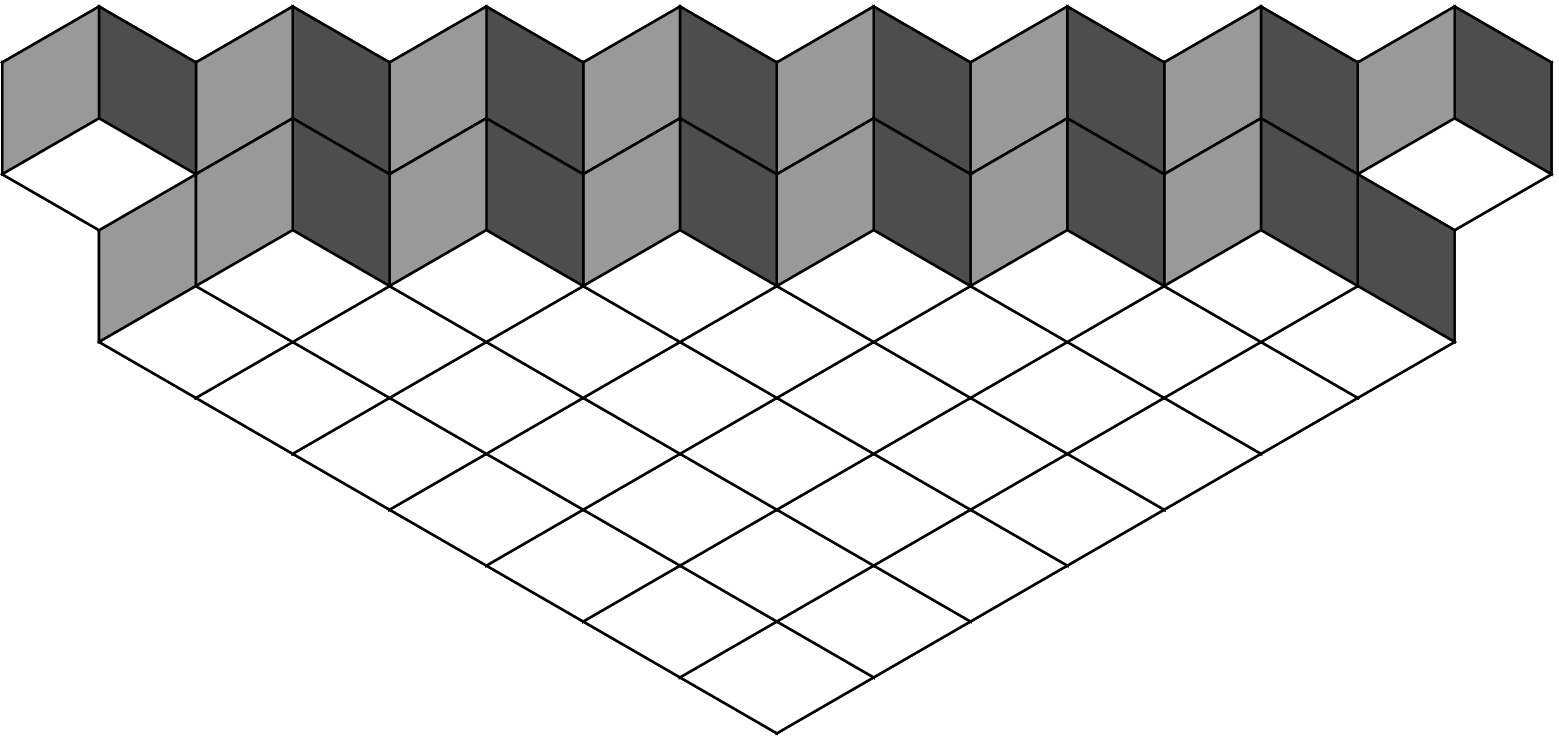, width=3cm}
	\hfill\mbox{}
	&
	\mbox{}\hfill
	\epsfig{file=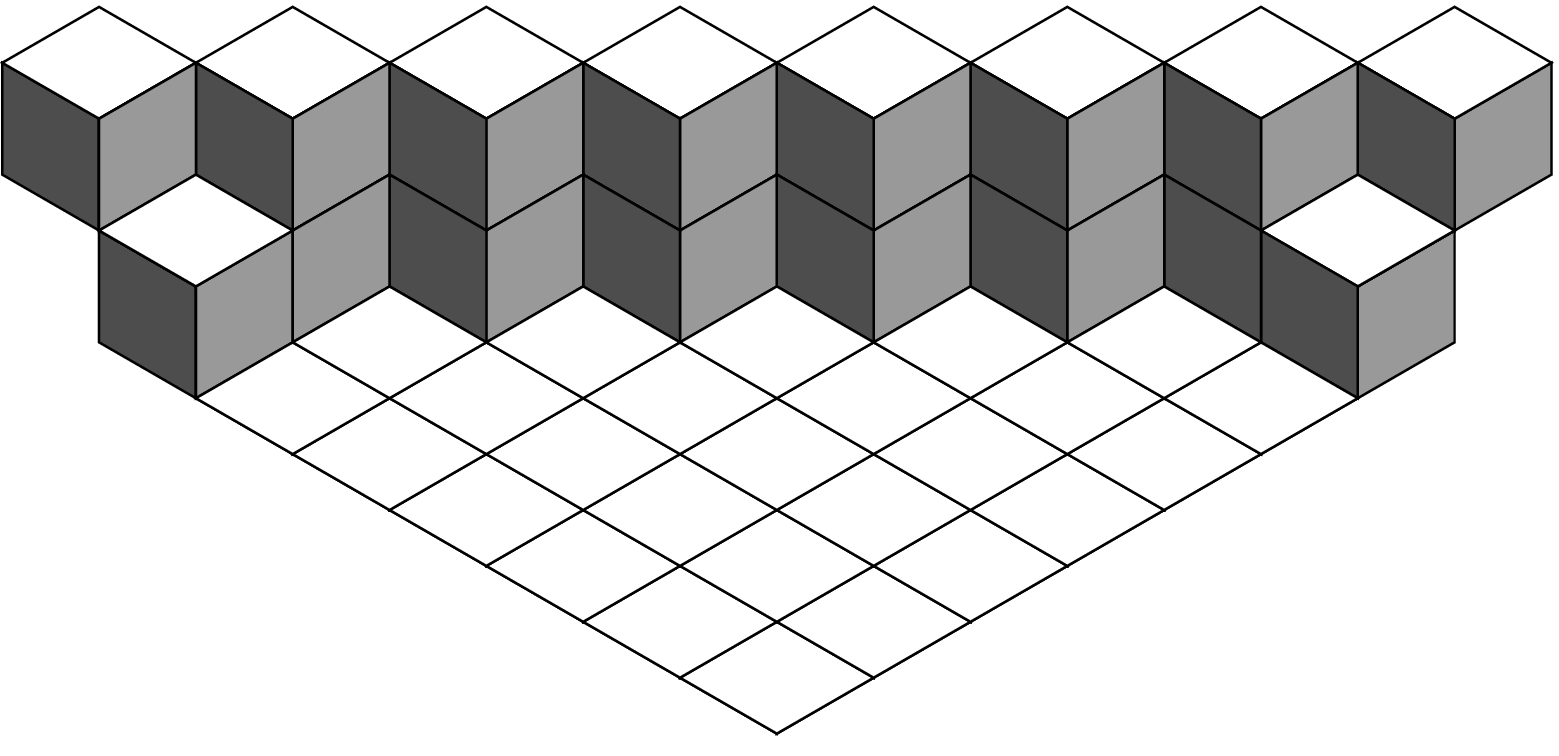, width=3cm}
	\hfill\mbox{}
	&
	\mbox{}\hfill
	\epsfig{file=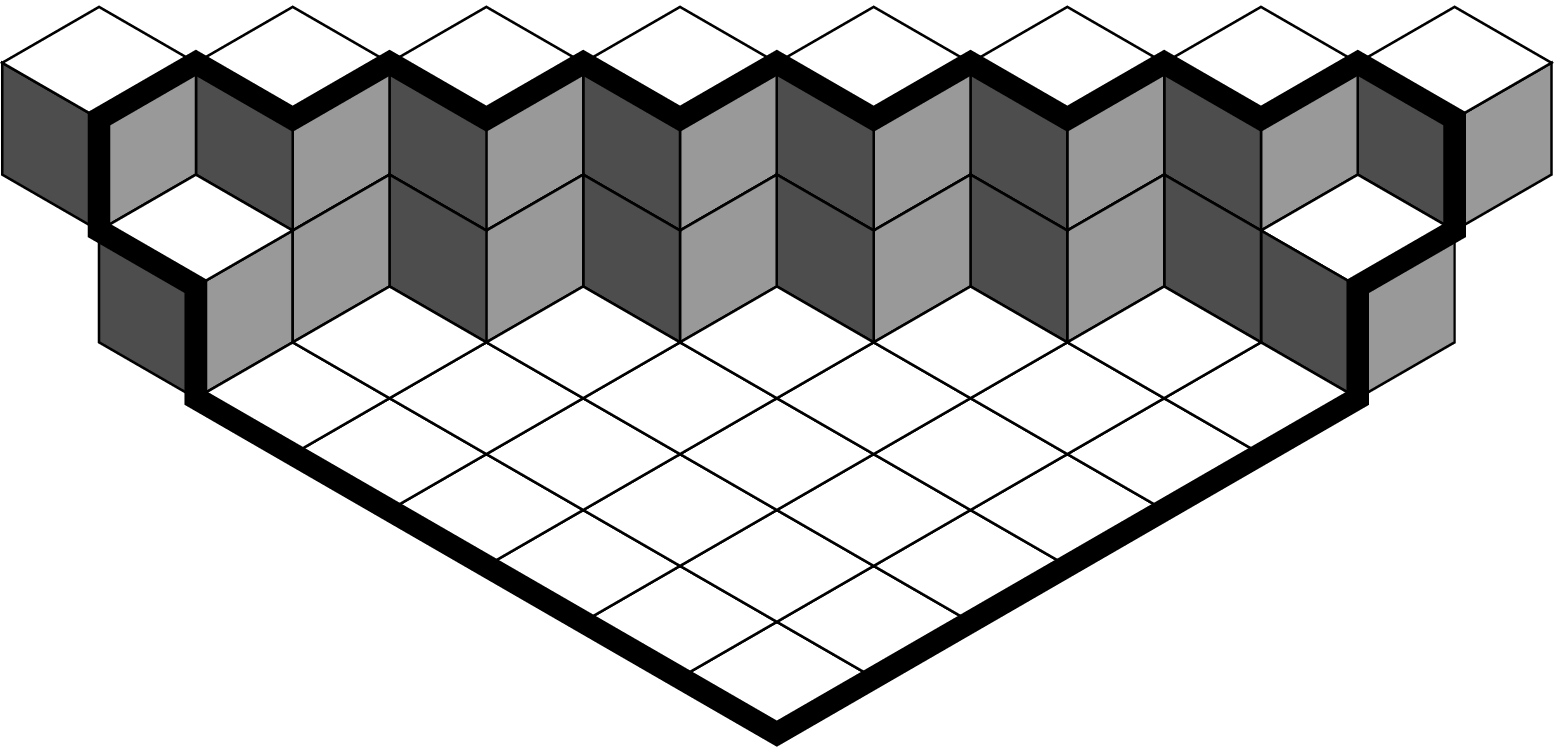, width=3cm}
	\hfill\mbox{}
	\\[1mm]

	\mbox{}\hfill
	(a) 	\begin{minipage}[t]{2cm}
		Seeds of order 0 in ${\cal D}_n$
		\end{minipage}
	\hfill\mbox{}
	&
	\mbox{}\hfill
	(b) 	\begin{minipage}[t]{2cm}
		Seeds of order 1 in ${\cal D}_{n}$
		\end{minipage}
	\hfill\mbox{}
	&
	\mbox{}\hfill
	(c) 	\begin{minipage}[t]{3.4cm}
		Seeds of order 1 in ${\cal D}_n$ are of order 0 in ${\cal 
		D}_{n-1}$
		\end{minipage}
	\hfill\mbox{}
	\\

	\end{tabular}
	\caption{A family of domains \worklabel{fig:seeds:n}
	}
}

%%%%%%%%%%%%%%%%%%%%%%%%%%%%%%%%%%%%%%%%%%%%%%%%%%%%%%%%%%%%%%%%%%%%%%%%
%%%%%%%%%%%%%%%%%%%%%%%%%%%%%%%%%%%%%%%%%%%%%%%%%%%%%%%%%%%%%%%%%%%%%%%%
%%%%%%%%%%%%%%%%%%%%%%%%%%%                   %%%%%%%%%%%%%%%%%%%%%%%%%%
%%%%%%%%%%%%%%%%%%%%%%%%%%%   Sous-treillis   %%%%%%%%%%%%%%%%%%%%%%%%%%
%%%%%%%%%%%%%%%%%%%%%%%%%%%                   %%%%%%%%%%%%%%%%%%%%%%%%%%
%%%%%%%%%%%%%%%%%%%%%%%%%%%%%%%%%%%%%%%%%%%%%%%%%%%%%%%%%%%%%%%%%%%%%%%%
%%%%%%%%%%%%%%%%%%%%%%%%%%%%%%%%%%%%%%%%%%%%%%%%%%%%%%%%%%%%%%%%%%%%%%%%

\section{$C$-minimal tilings and intervals in $\cal L(D)$}
	\worklabel{lattices:intervals}

A close inspection of the lattice $\cal L(D)$ of the tilings of ${\cal D}$
will reveal that it contains interesting intervals, which we will use to
generate all the tilings of a domain by performing the required flips in an
orderly fashion.

%!!!!!!!!!!!!!!!!!!!                                 !!!!!!!!!!!!!!!!!!!
%!!!!!!!!!!!!!!!!!!!   Treillis dans celui d'un PH   !!!!!!!!!!!!!!!!!!!
%!!!!!!!!!!!!!!!!!!!                                 !!!!!!!!!!!!!!!!!!!

\subsection{$C$-minimal tilings of a pseudo-hexagon}

Let ${\cal D}$ be a pseudo-hexagon and let~$C$ be a cube in a filling of the
proper seed~$s$ of ${\cal D}$; as we will see in
Section~\ref{subsection:general:intervals}, marking such a cube can be very
useful. What can one say about the fillings of~$s$ in which~$C$ appears? Since
fillings are compact piles of cubes, each such filling must also contain all
the cubes in the parallelepiped defined by the diagonal going from~$s$ to~$C$
(see Figure~\ref{fig:min:lozenge}~(a)).

\fig{
	\begin{tabular}{p{6cm}@{\hskip5mm}
			p{5cm}
			}

	\mbox{}\hfill
	\begin{minipage}{1.8cm}
		\mbox{}\hfill
		\epsfig{file=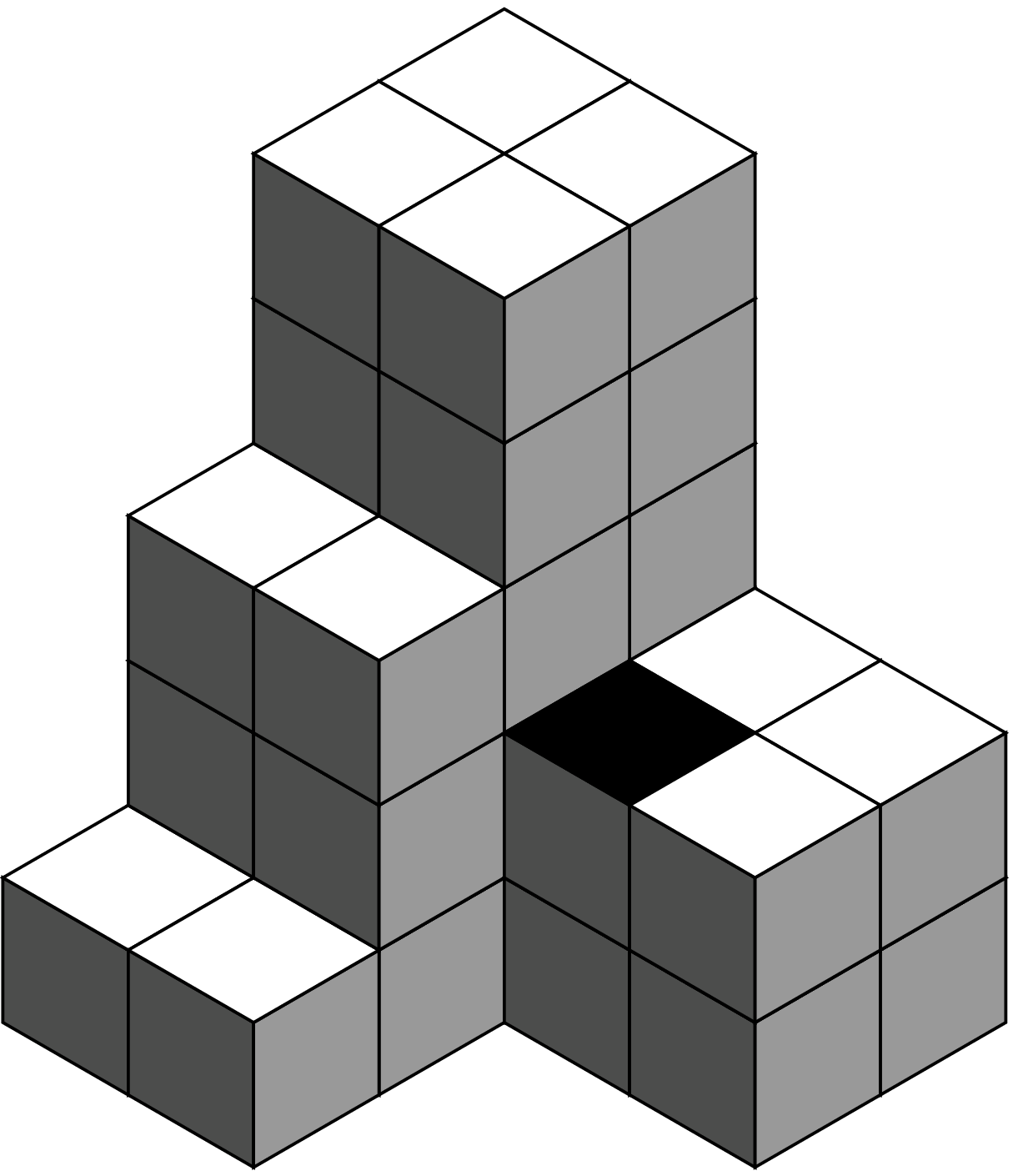, width=1.5cm}
		\hfill\mbox{}
	\end{minipage}
	\hskip5mm
	\begin{minipage}{1.8cm}
		\mbox{}\hfill
		\epsfig{file=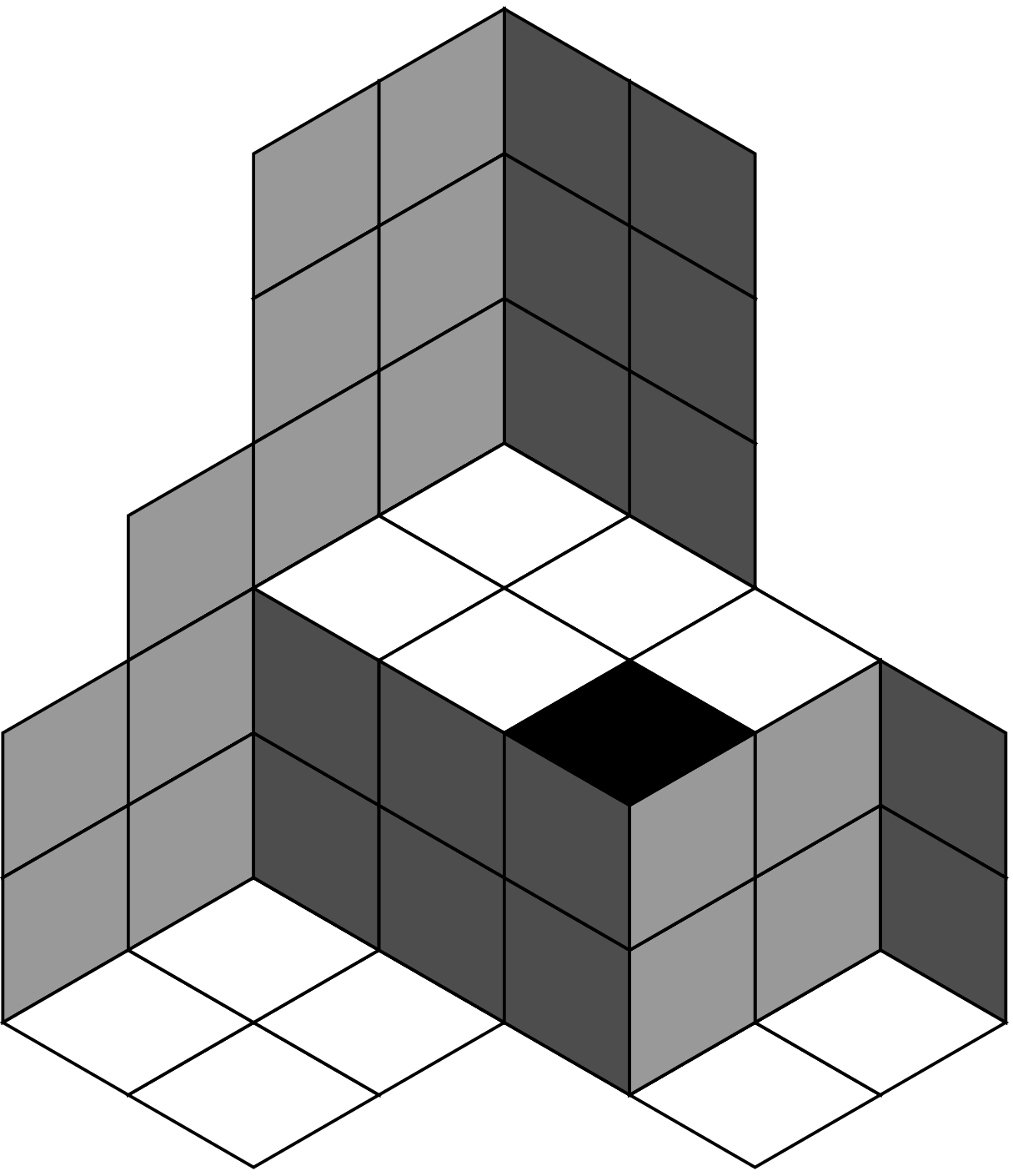, width=1.5cm}
		\hfill\mbox{}
	\end{minipage}
	\hfill\mbox{}
	&
	\mbox{}\hfill
	\begin{minipage}{1.8cm}
		\epsfig{file=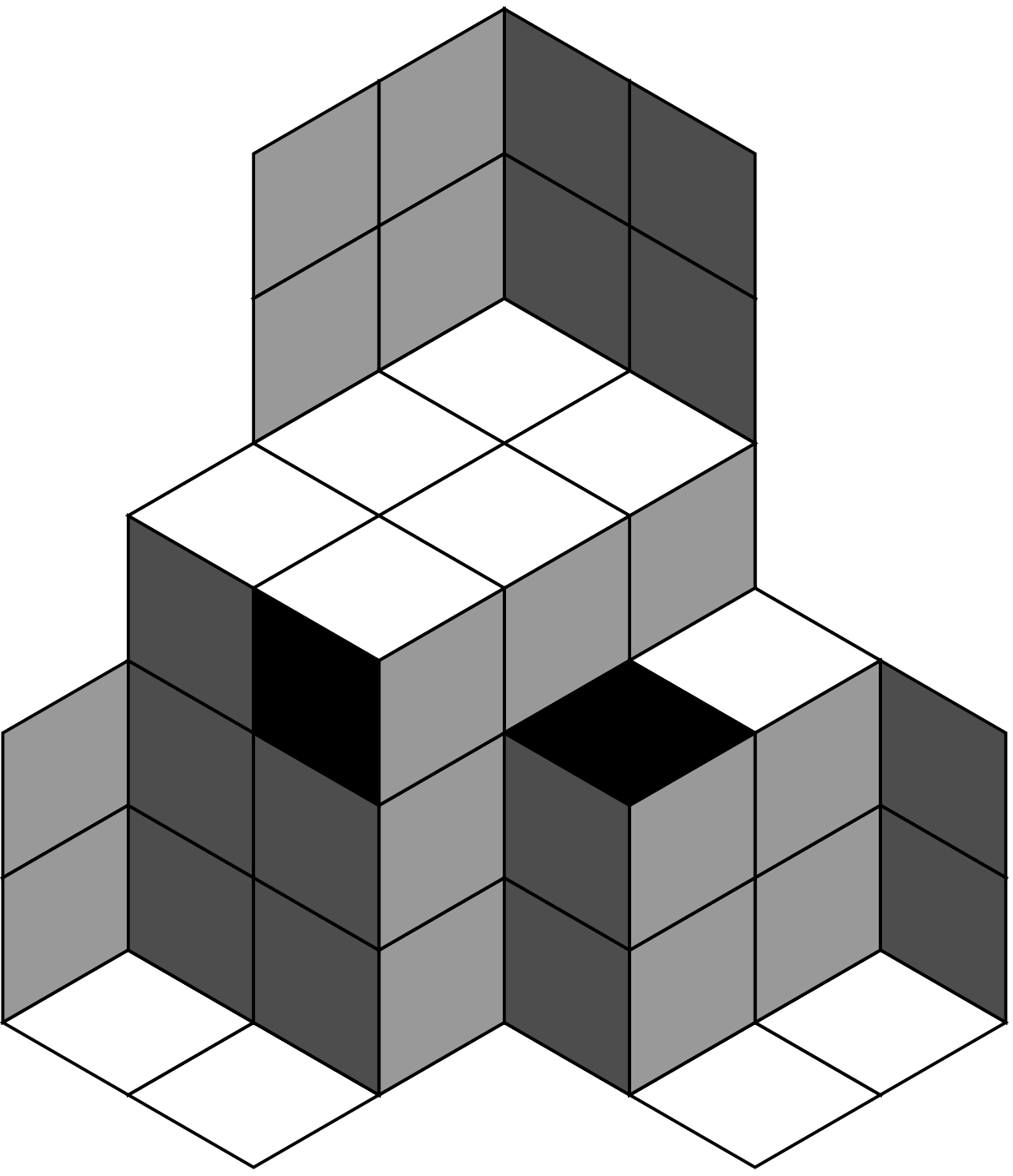, width=1.5cm}
	\end{minipage}
	\hfill\mbox{}
	\\[3mm]

	\mbox{}\hfill
	(a) 	\begin{minipage}[t]{5.4cm}
		A lozenge in $\text{Max}(\cal D)$ and the corresponding 
		$C$-minimal tiling
		\end{minipage}
	\hfill\mbox{}
	&
	\mbox{}\hfill
	(b) 	\begin{minipage}[t]{4.4cm}
		$\text{Min}_{(C_i)}(T)$ with \#$\{C_i\} = 2$
		\end{minipage}
	\hfill\mbox{}
	\\

	\end{tabular}
	\caption{Tilings defined by the presence of a lozenge 
		\worklabel{fig:min:lozenge}
	}
}

\begin{definition}
	Let $P$ be a pseudo-hexagon and $C$ a cube in a filling of its proper
seed~$s$. The $C$-\emph{minimal} tiling of $P$ (which we denote by
$\text{Min}_C(P)$) is the tiling associated with the filling of~$s$ obtained
by starting from $\text{Min}(T)$ and adding all the cubes in the
parallelepiped defined by the diagonal $(s,C)$.

	\worklabel{def:C:minimal:tiling}
	\end{definition}

This definition is consistent. First, a $C$-minimal tiling is a tiling since
the parallelepiped (${\cal P}(C)$) is a compact pile of cubes. Second, if a
filling of $s$ contains ${\cal P}(C)$ then it contains $C$. Finally, if a
filling of $s$ contains $C$ then it must contain ${\cal P}(C)$ since a filling
is a compact pile.

\begin{definition}
	Let $P$ be a pseudo-hexagon and $(C_i)_{1 \leqslant i \leqslant n}$ a
collection of cubes in a filling of the proper seed~$s$ of~$P$. Let $T_k$ be
the $C_k$-minimal tiling of~$P$ for $k\in\{1,\ldots,n\}$.  The
$(C_i)$-\emph{minimal} tiling of $T$ (which we denote by $\text{\rm
Min}_{(C_i)}(P)$) is $\text{\rm Sup}(T_1,\ldots,T_n)$.

	\worklabel{def:Ci:minimal:tiling}
	\end{definition}

An example of such a tiling is given in Figure~\ref{fig:min:lozenge}~(b).

\medskip

All of the above can be reformulated in a dual way. Instead of adding cubes
from the minimal tiling, remove cubes from the maximal one: this can be viewed
as adding ``\emph{anti-cubes}'', the rule being that cubes and anti-cubes
annihilate each other. We could thus define a $C$-\emph{maximal} tiling as a
tiling of ${\cal D}$ which does \emph{not} contain the cube $C$. The notion
can be extended to families, thus defining $(C_i)$-\emph{maximal} tilings.

%!!!!!!!!!!!!!!!!!!!!                              !!!!!!!!!!!!!!!!!!!!
%!!!!!!!!!!!!!!!!!!!!   Intervalles fondamentaux   !!!!!!!!!!!!!!!!!!!!
%!!!!!!!!!!!!!!!!!!!!                              !!!!!!!!!!!!!!!!!!!!

\subsection{The fundamental intervals of $\cal L(D)$}
	\worklabel{subsection:fundamental:intervals}

Let now ${\cal D}$ be any domain (such that the results of
Section~\ref{flips:lattices} hold) and let ${\cal T}({\cal D})$ be the set of
its tilings and $\cal L(D)$ the associated lattice. We have seen that
$\text{Max}({\cal D})$, seen as a pile of cubes, can be partitionned according
to the maximal fillings of its proper seeds.

\begin{definition}
	We denote by ${\cal D}^n$, $n \geqslant 1$, the tiling of ${\cal D}$
obtained by maximally filling the ranges of all the proper seeds of order $k
\in \{1,\ldots,n\}$. By convention, ${\cal D}^0 = \text{Min}({\cal D})$. The
\emph{degree} $d({\cal D})$ of ${\cal D}$ is the minimum of the integers $n$
such that ${\cal D}^n = \text{\rm Max}({\cal D})$. A cube~$C$ is of
\emph{order} $k$ if it belongs to the filling of a seed of order~$k$.

	\worklabel{def:degree}
	\end{definition}

Note that $\text{Max}({\cal D}^k) = \text{Min}({\cal D}^{k+1})$ for $0
\leqslant k < d({\cal D})$.  See Figure~\ref{fig:successive:D:k} for an
example.

\fig{
	\begin{tabular}{p{2.2cm}@{\hskip9mm}
			p{2.2cm}@{\hskip9mm}
			p{2.2cm}@{\hskip9mm}
			p{2.2cm}
			}

	\mbox{}\hfill
	\epsfig{file=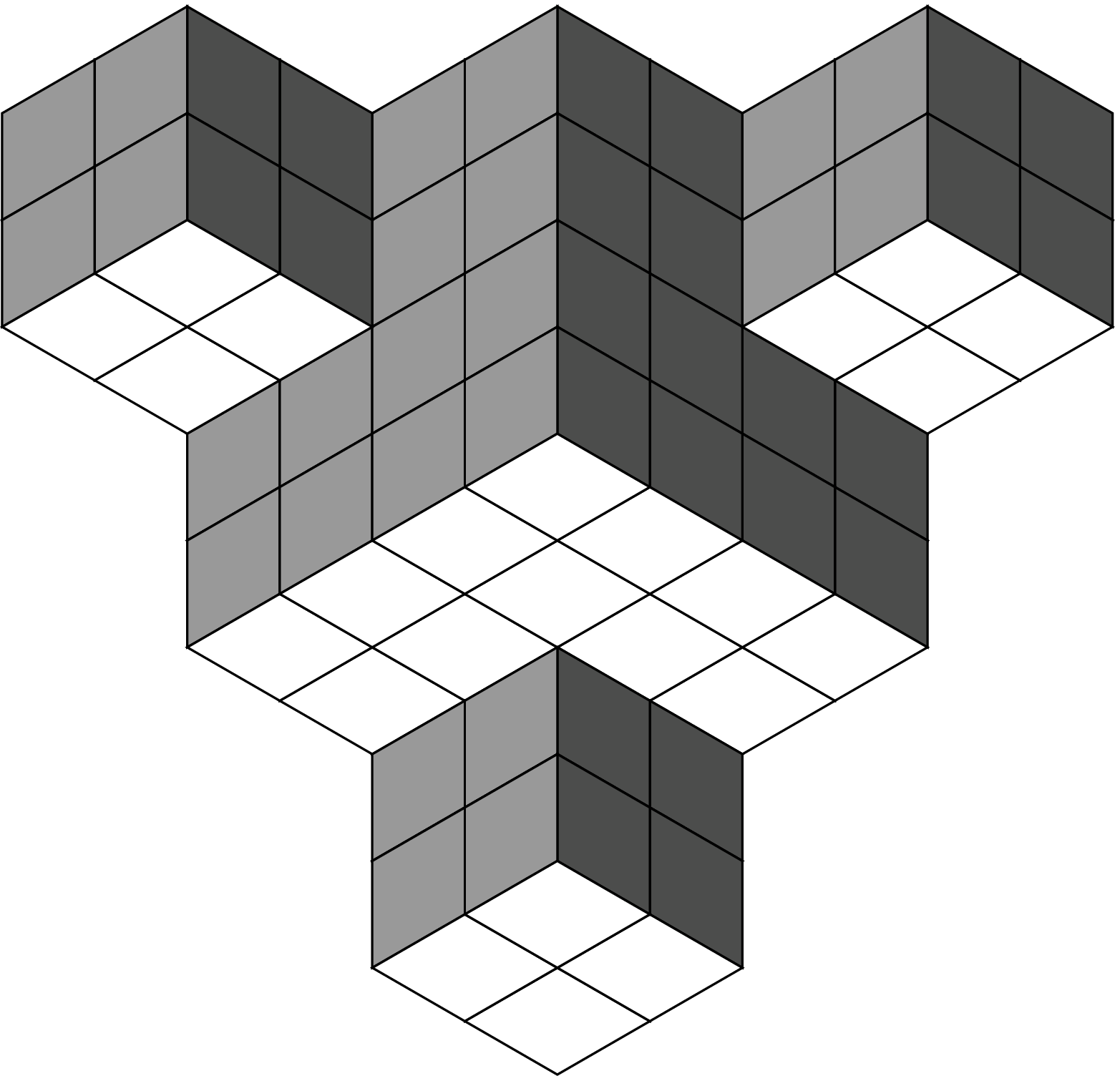, width=2cm}
	\hfill\mbox{}
	&
	\mbox{}\hfill
	\epsfig{file=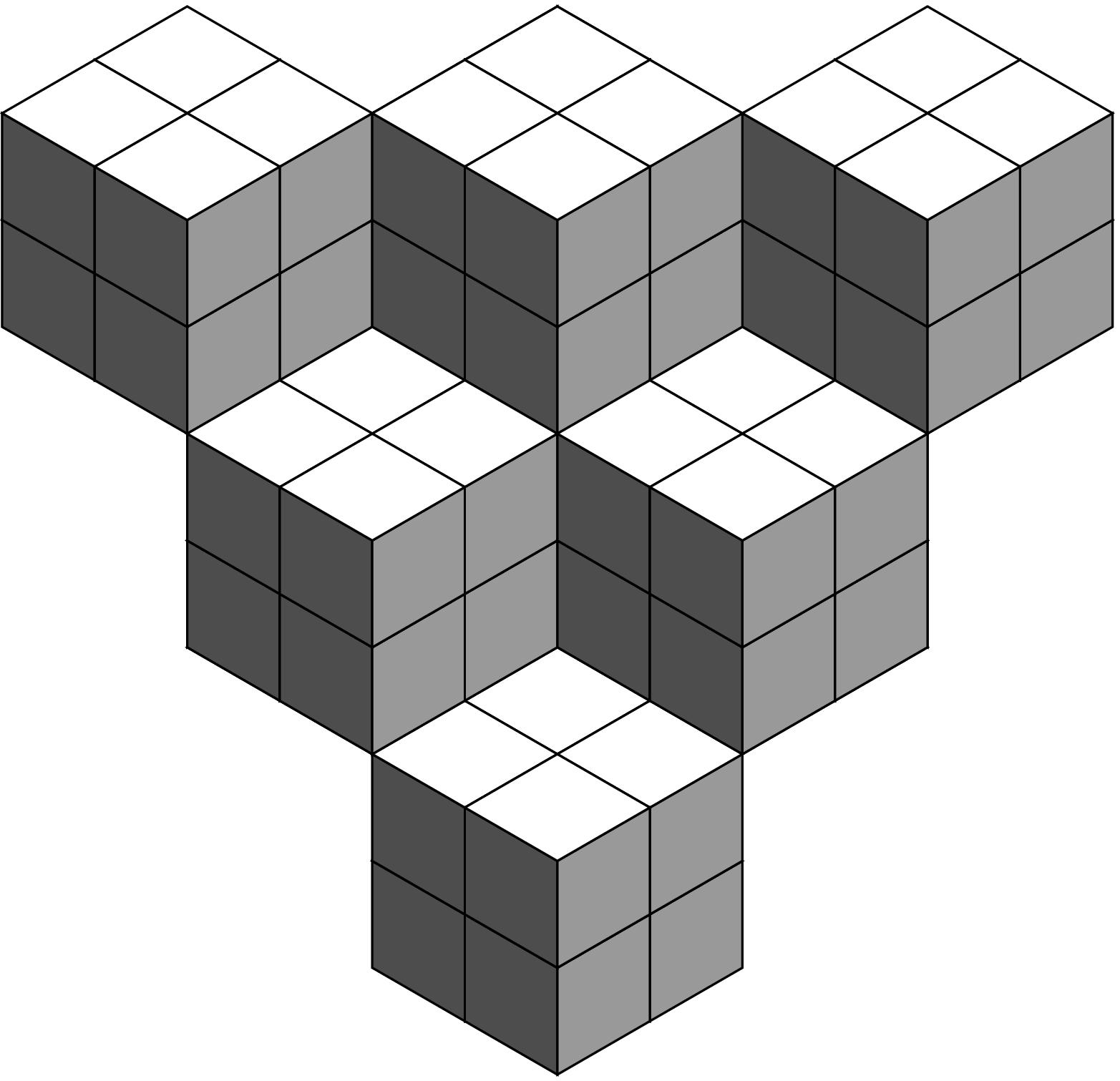, width=2cm}
	\hfill\mbox{}
	&
	\mbox{}\hfill
	\epsfig{file=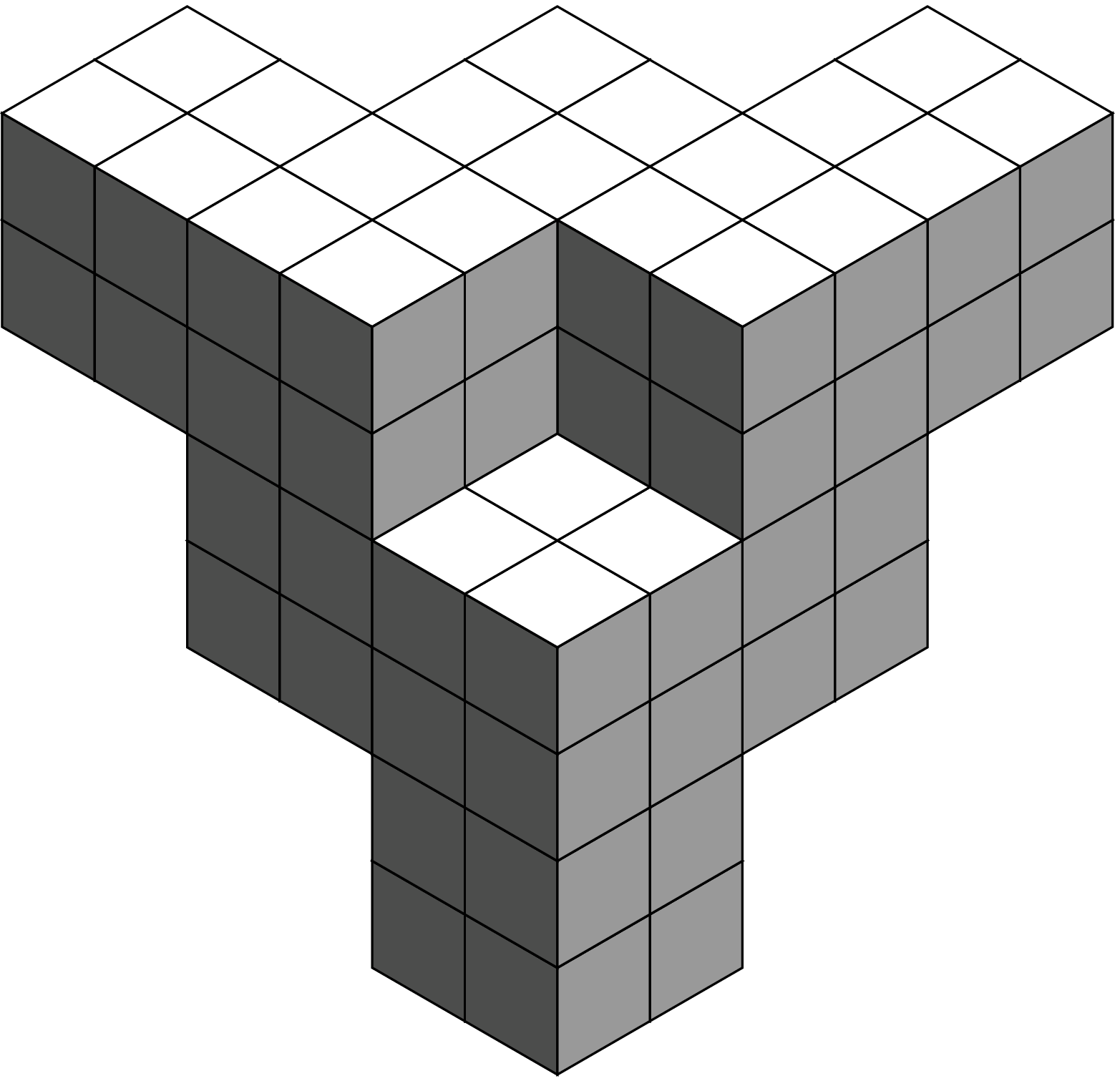, width=2cm}
	\hfill\mbox{}
	&
	\mbox{}\hfill
	\epsfig{file=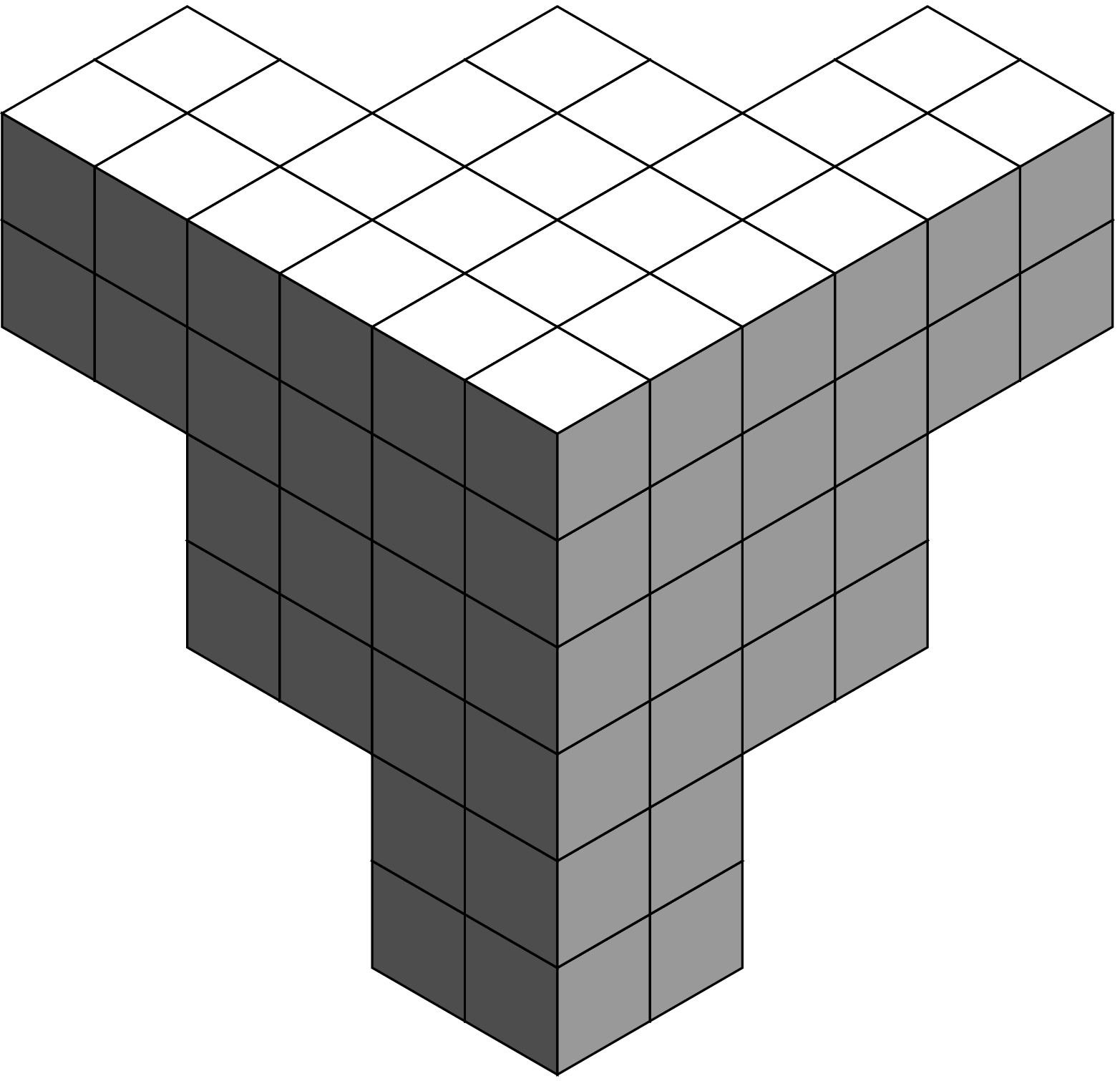, width=2cm}
	\hfill\mbox{}
	\\[1mm]

	\mbox{}\hfill
	${\cal D}^0$
	\hfill\mbox{}
	&
	\mbox{}\hfill
	${\cal D}^1$
	\hfill\mbox{}
	&
	\mbox{}\hfill
	${\cal D}^2$
	\hfill\mbox{}
	&
	\mbox{}\hfill
	${\cal D}^3$
	\hfill\mbox{}
	\\

	\end{tabular}
	\caption{The successive ${\cal D}^k$ for $k=0..d({\cal D})$
		\worklabel{fig:successive:D:k}
	}
}

\begin{notation}
	The set of all the tilings $t$ of ${\cal D}$ for which ${\cal D}^n
\preccurlyeq t \preccurlyeq {\cal D}^{n+1}$ is an interval in the lattice of
the tilings of $\cal D$; we denote this interval by ${\cal T}_n({\cal D})$, $0
\leqslant n < d({\cal D})$). By extension, $\cal T(D)$ denotes the set of all
the tilings of $\cal D$.

	\worklabel{notation:interval:tilings}
	\end{notation}

\begin{definition}
	The ${\cal T}_n({\cal D})$, $0 \leqslant n < d({\cal D})$, are the
\emph{fundamental intervals} of the lattice $\cal L(D)$ of the tilings of
${\cal D}$.

	\worklabel{def:fundamental:intervals}
	\end{definition}

\begin{proposition}
	$\!({\cal T}_0({\cal D}), \ldots, {\cal T}_{d({\cal D})-1}({\cal D}))$
is a maximal chain of intervals in $\cal L(D)$.

	\worklabel{prop:fundamental:chain}
	\end{proposition}

\noindent \textit{Proof} \quad It follows from the definitions that ${\cal
D}^0 = \text{Min}({\cal D})$, ${\cal D}^{d({\cal D})} = \text{Max}({\cal D})$
and $\text{Max}({\cal D}^k) = \text{Min}({\cal D}^{k+1})$ for $0 \leqslant k <
d({\cal D})$.~\mbox{}\hfill$\square$

\medskip

A graphical representation of this chain is given in
Figure~\ref{fig:fundamental:intervals}.

\begin{figure}
\mbox{}\hfill
\epsfig{file=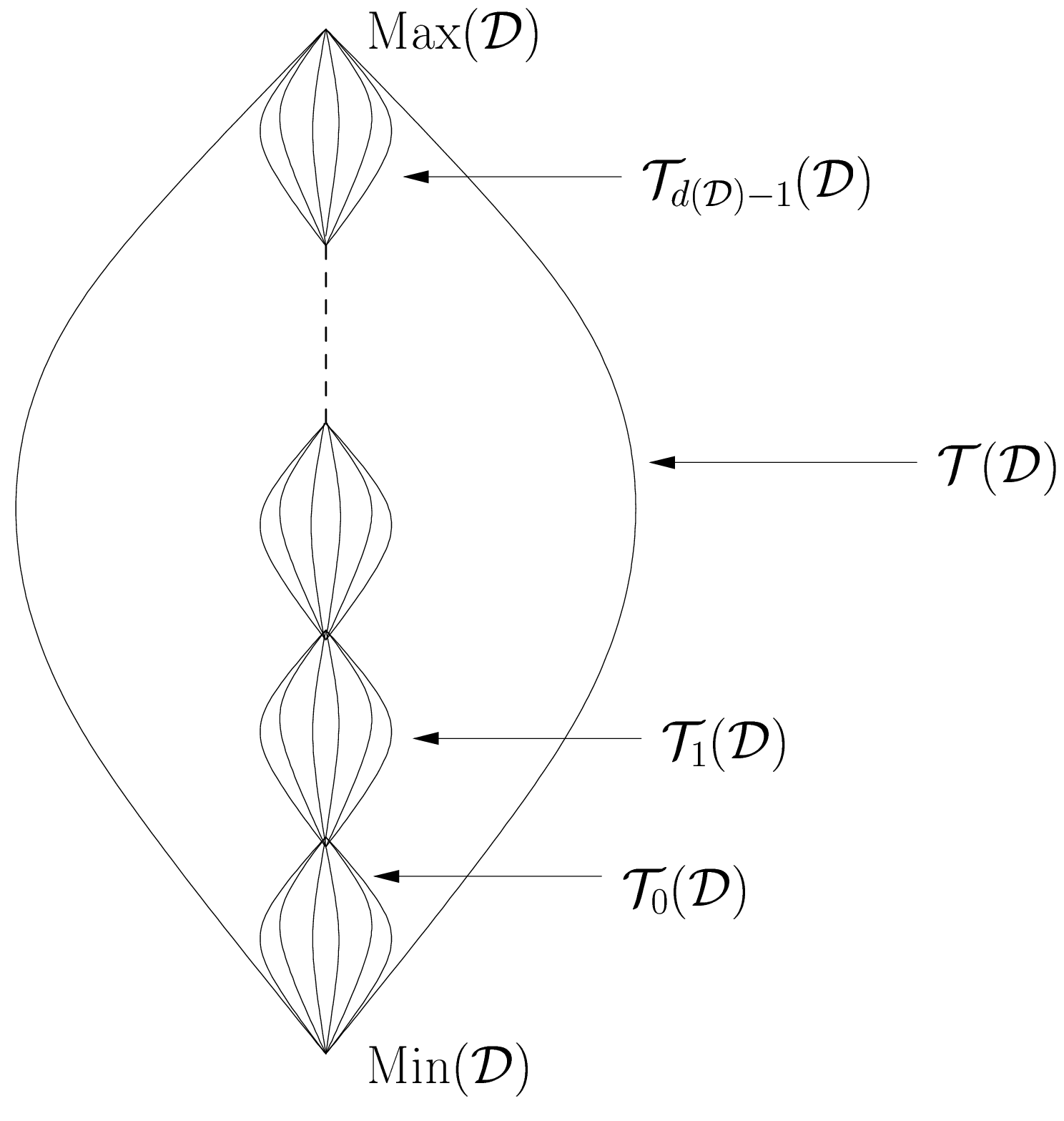, width=4cm}
\hfill\mbox{}
\caption{The fundamental intervals of ${\cal T}({\cal D})$
	\worklabel{fig:fundamental:intervals}
	}
\end{figure}

%!!!!!!!!!!!!!!!!!!!!!!                          !!!!!!!!!!!!!!!!!!!!!!
%!!!!!!!!!!!!!!!!!!!!!!   Intervalles généraux   !!!!!!!!!!!!!!!!!!!!!!
%!!!!!!!!!!!!!!!!!!!!!!                          !!!!!!!!!!!!!!!!!!!!!!

\subsection{$C$-minimal tilings of a fertile zone}
	\worklabel{subsection:general:intervals}

Let~$\cal D$ be a domain such that the results of Section~\ref{flips:lattices}
hold. Using fracture lines (see Section~\ref{fracture:lines}), we can suppose
that~$\cal D$ is a fertile zone. Let $C$ be a cube that appears in a filling
of a proper seed $s$ of order~$n$ of~$\cal D$. It is clear from the
definitions that $C \in {\cal D}^n$. $C$ induces a $C$-minimal tiling of
$R(s)$, let us denote it by ${\cal P}(C)$, whose contour path is a
pseudo-hexagon (whose proper seed is~$s$) and therefore a range of~$s$ (see
Figure~\ref{fig:min:D}~(a)).

\fig{
	\begin{tabular}{p{5.4cm}@{\hskip9mm}p{5.4cm}}

	\mbox{}\hfill
	\begin{minipage}{2.2cm}
		\begin{center}
		\epsfig{file=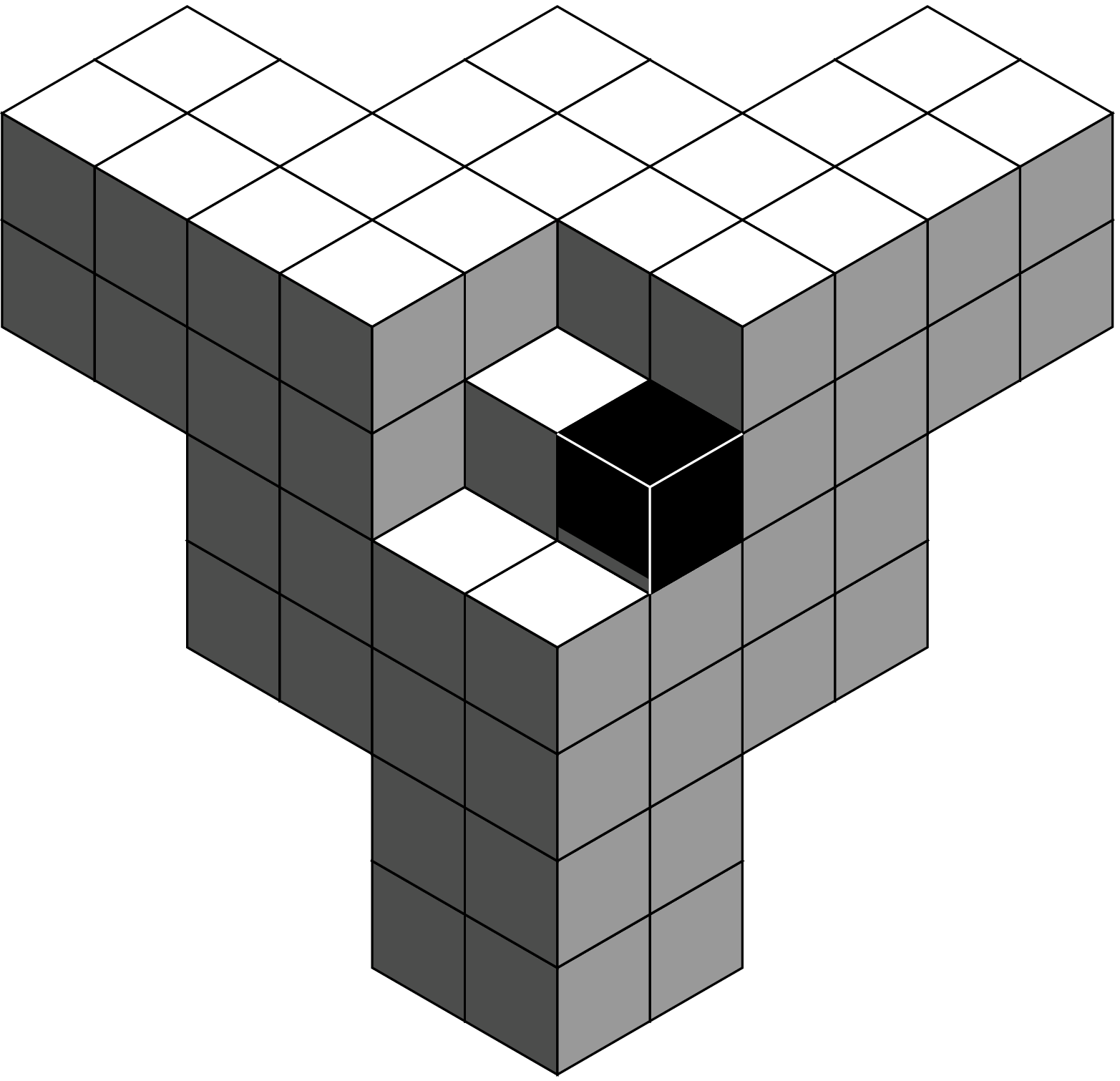, width=2cm}
		\end{center}
	\end{minipage}
		\hskip7mm
	\begin{minipage}{2.2cm}
		\begin{center}
		\epsfig{file=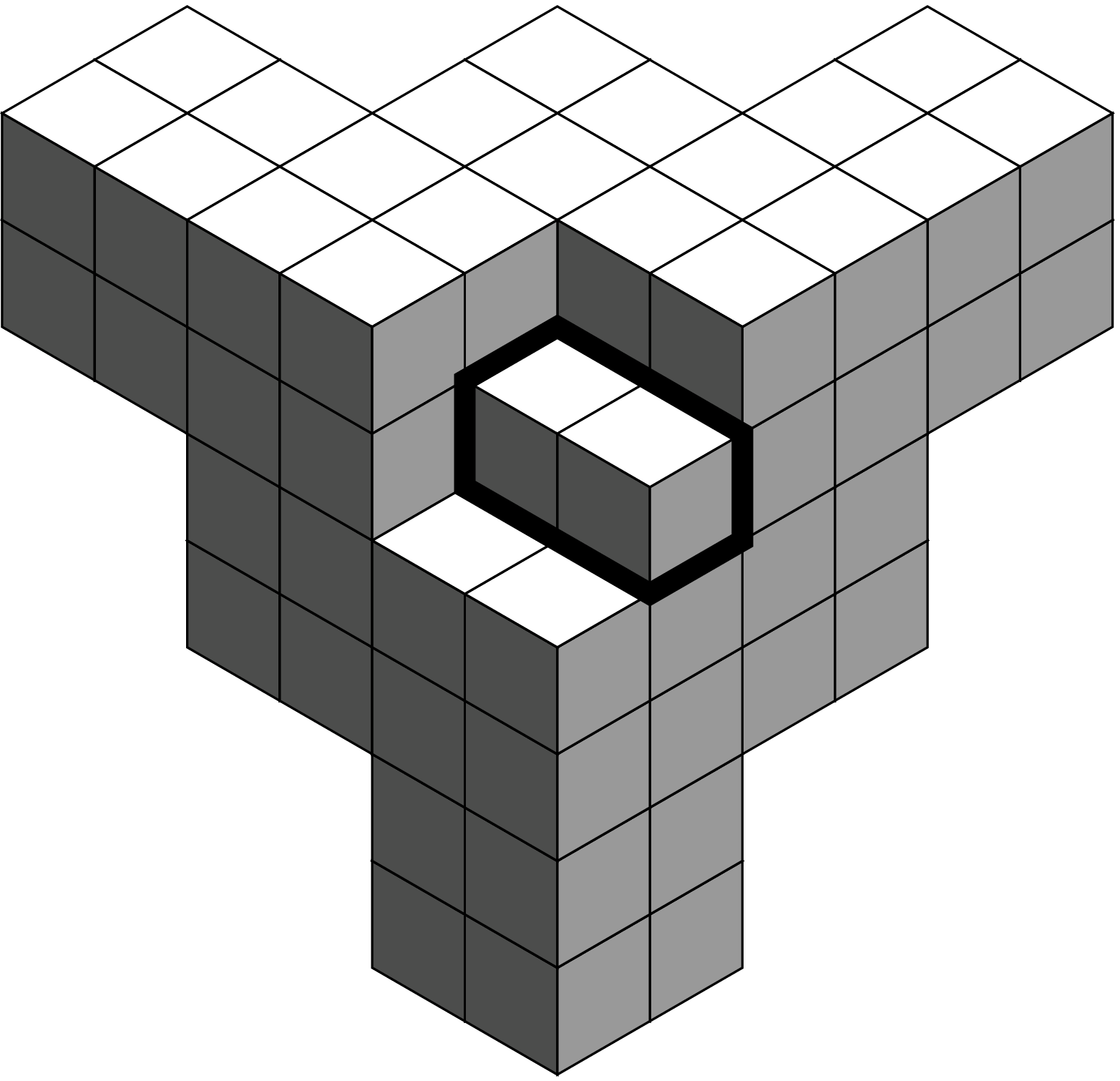, width=2cm}
		\end{center}
	\end{minipage}
	\hfill\mbox{}
	&
	\mbox{}\hfill
	\begin{minipage}{2,2cm}
		\begin{center}
		\epsfig{file=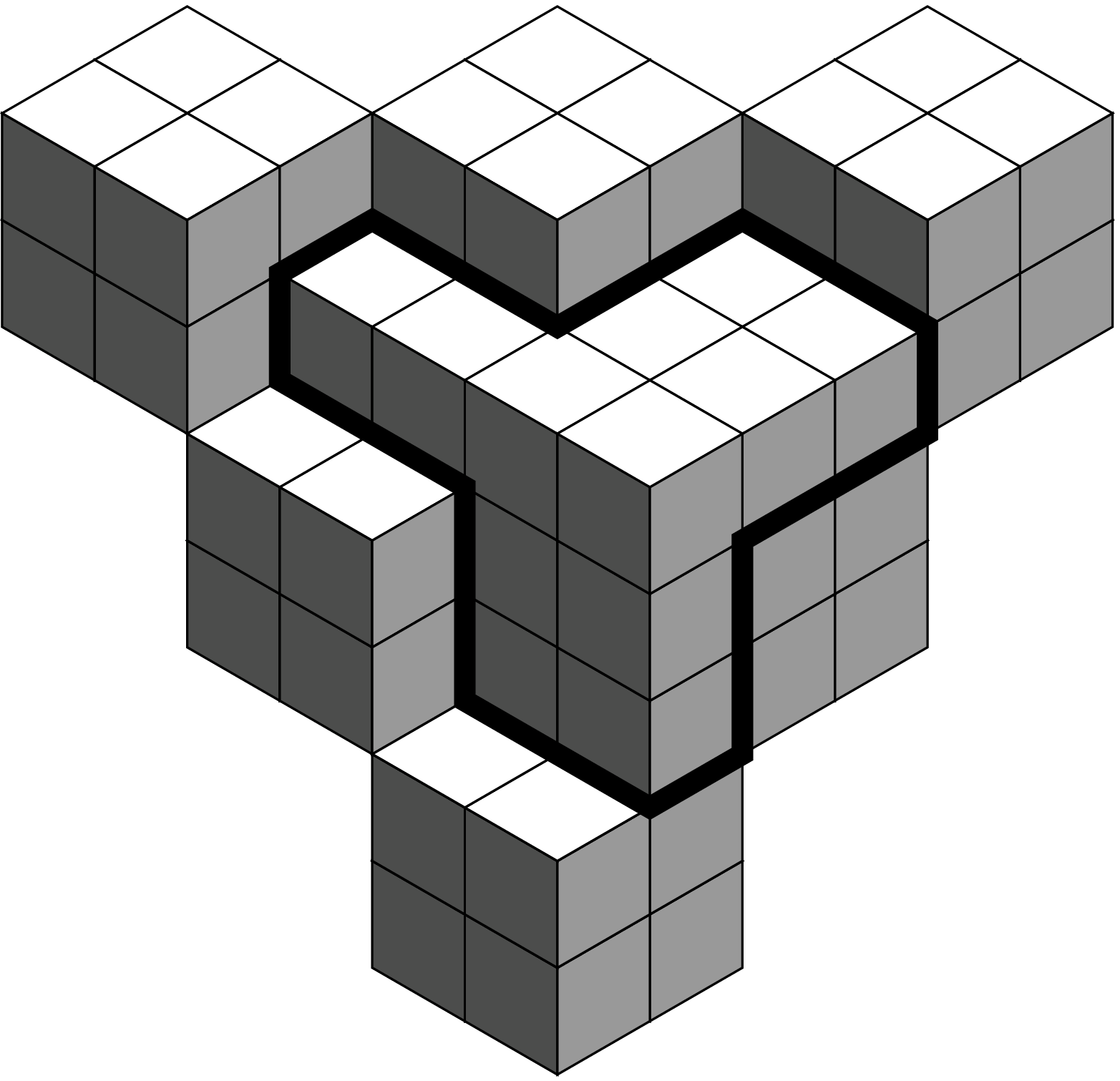, width=2cm}
		\end{center}
	\end{minipage}
	\hskip7mm
	\begin{minipage}{2.2cm}
		\begin{center}
		\epsfig{file=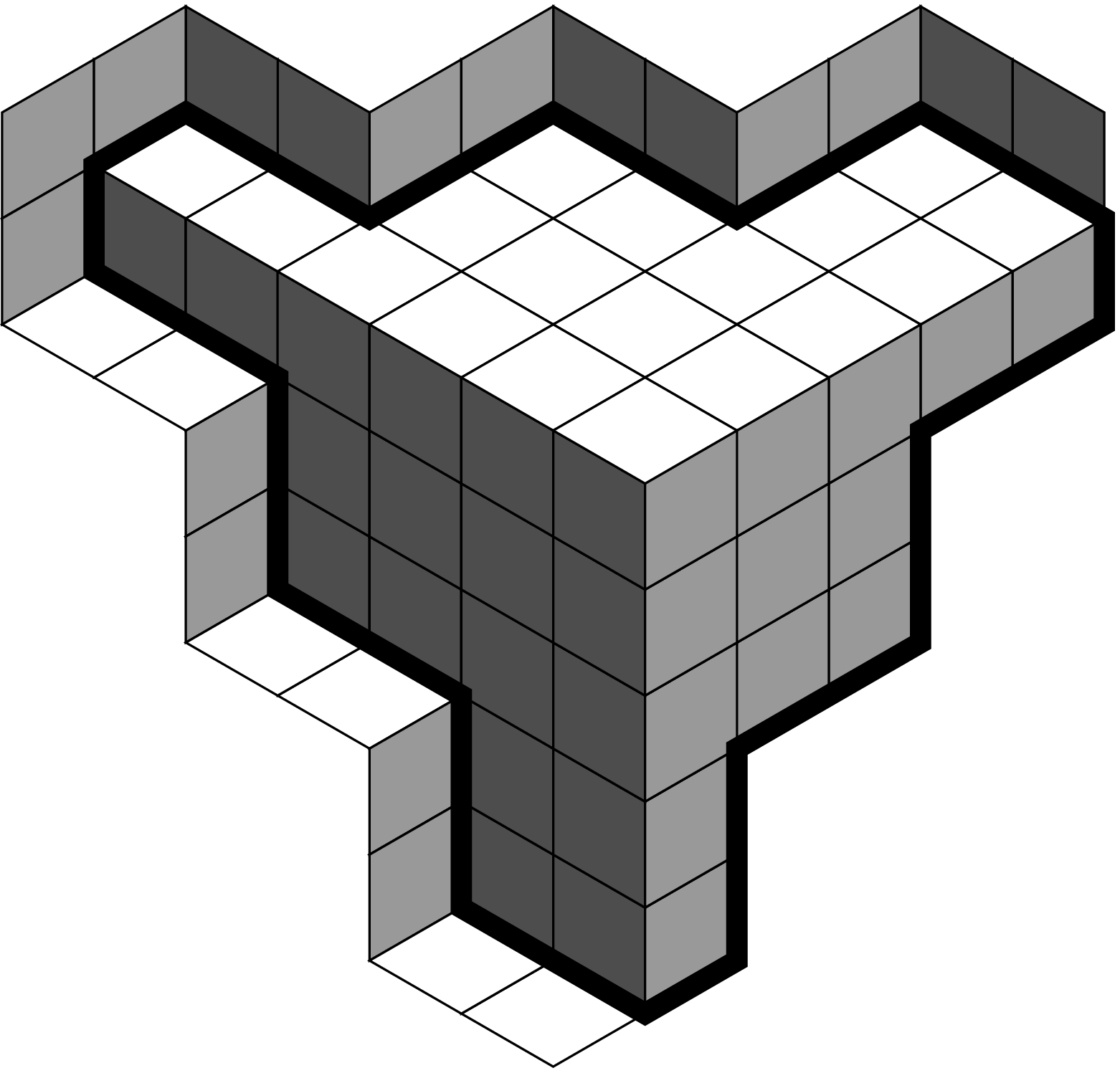, width=2cm}
		\end{center}
	\end{minipage}
	\hfill\mbox{}
	\\[1mm]

	\mbox{}\hfill
	(a) 	\begin{minipage}[t]{4.4cm}
		The contour path of a $C$-minimal tiling is a range of a 
		proper seed
		\end{minipage}
	\hfill\mbox{}
	&
	\mbox{}\hfill
	(b) 	\begin{minipage}[t]{4.6cm}
		Building the smallest of the tilings containing $C$
		\end{minipage}
	\hfill\mbox{}
	\\

	\end{tabular}
	\caption{$C$-minimal tilings in $\cal D$ \worklabel{fig:min:D}}
}

\begin{definition}
	The $C$-minimal tiling of ${\cal D}$, let us denote it by
$\text{Min}_C({\cal D})$, is the infimum of all the tilings of $\cal D$ that
contain~$C$. If $(C_i)_{1 \leqslant i \leqslant n}$ is a family of cubes,
the $(C_i)$-minimal tiling of~$\cal D$ is $\text{Min}_{(C_i)}({\cal D}) =
\text{Sup}(\text{Min}_{C_1}({\cal D}), \ldots, \text{Min}_{C_n}({\cal D}))$.

	\end{definition}

We now sketch a recursive construction of~$\text{Min}_{C}({\cal D})$ (see
Figure~\ref{fig:min:D}~(b)). ${\cal P}(C)$ marks cubes that belong to~${{\cal
D}}^{\:n-1}$ and therefore to the maximal fillings of seeds $t_i$ of order $k
\leqslant n-1$. These cubes give rise to~$(C_j)$-minimal tilings in the~$t_i$,
which in turn mark cubes in~ ${{\cal D}}^{\:n-2}$, and so on. Since we always
choose minimal tilings in the ranges of the proper seeds, the tiling is less
than any tiling containing~$C$, and it contains~$C$, therefore it
is~$\text{Min}_{C}({\cal D})$.

%%%%%%%%%%%%%%%%%%%%%%%%%%%%%%%%%%%%%%%%%%%%%%%%%%%%%%%%%%%%%%%%%%%%%%%%
%%%%%%%%%%%%%%%%%%%%%%%%%%%%%%%%%%%%%%%%%%%%%%%%%%%%%%%%%%%%%%%%%%%%%%%%
%%%%%%%%%%%%%%%%%%%%%%%%%%%%                 %%%%%%%%%%%%%%%%%%%%%%%%%%%
%%%%%%%%%%%%%%%%%%%%%%%%%%%%   Cas general   %%%%%%%%%%%%%%%%%%%%%%%%%%%
%%%%%%%%%%%%%%%%%%%%%%%%%%%%                 %%%%%%%%%%%%%%%%%%%%%%%%%%%
%%%%%%%%%%%%%%%%%%%%%%%%%%%%%%%%%%%%%%%%%%%%%%%%%%%%%%%%%%%%%%%%%%%%%%%%
%%%%%%%%%%%%%%%%%%%%%%%%%%%%%%%%%%%%%%%%%%%%%%%%%%%%%%%%%%%%%%%%%%%%%%%%

\section{An algorithm to generate the lattice of the tilings of a general
domain}
\worklabel{section:algo:main}

In Section~\ref{section:algo:hexagons} we have proposed an algorithm to
generate efficiently the tilings of a pseudo-hexagon. We now proceed to the
general case, where no hypothesis is made on the shape of the domain~$\cal D$,
except that the results of Section~\ref{flips:lattices} hold.  

In order to make use of already found tilings and reduce the computation time,
we need the following important fact: if $(C_i)$ is a collection a cubes, all
of them of order~$n$, the interval of~$\cal L(D)$ between
$\text{Min}_{(C_i)}(\cal D)$ and ${\cal D}^{\:n}$ plus the cubes~$C_i$ is
isomorphic to the interval between $\text{Min}_{(C_i)}($\cal D$)$ minus the
cubes~$C_i$ and ${\cal D}^{\:n}$. This is obvious because cubes of order~$k <
n$ can be added to $\text{Min}_{(C_i)}(\cal D)$ whether the cubes of order~$n$
are present or not.

Formally stated, there is an isomorphism between the intervals
$[\text{Min}_{(C_i)}({\cal D}) \;;\;$ $\text{Sup}(\text{Min}_{(C_i)}({\cal
D}), {\cal D}^n)]$ and $[\text{Inf}(\text{Min}_{(C_i)}({\cal D}),{\cal D}^n)
\;;\; {\cal D}^n]$ of $\cal L(D)$.

\begin{algorithm}
\worklabel{algo:main}
\end{algorithm}

\begin{itemize}

\item \textbf{Input}: A finite-length closed path $\cal P$, enclosing a domain
$\cal D$, in the triangular grid.

\item \textbf{Output}: The lattice $\cal L(D)$ of the tilings of $\cal D$.

\item \textbf{Step 1}: Use Algorithm~\ref{algo:fracture:lines} to break the
domain into mutually independant fracture zones (see
Definition~\ref{def:fracture:zone}). The lattice of the tilings of a single
lozenge is trivial. For each of the fertile zones, follow steps~2 to~4.

\item \textbf{Step 2} ($\cal D$ is now a fertile zone): Recursively build the
tilings~${\cal D}^k$, $0 \leqslant k \leqslant d(\cal D)$. Set the list
$L_{\cal D}$ of the tilings of~$\cal D$ to~$\emptyset$.

\item \textbf{Step 3}: For $k=0 \text{ to } \text{deg}({\cal D})- 1$ do:

\begin{itemize}

\item \textbf{Step 3.1}: For each proper seed of order~$k$, use
Algorithm~\ref{algo:tilings} to generate all the fillings of the seed.

\item \textbf{Step 3.2}: Compute the cartesian product of these fillings.

\item \textbf{Step 3.3}: Each element of this product is a family~$(C_i)$ of
cubes; for each of these families, compute $\text{Min}_{(C_i)}(\cal D)$.

\item \textbf{Step 3.4}: Removing the cubes~$C_i$ from
$\text{Min}_{(C_i)}(\cal D)$ yields a collection of cubes of order at
most~$k-1$, which we denote by $\text{Min}_{(C_i)}^{k-1}(\cal D)$, or the
empty set if~$k=0$. Look up in $L_{\cal D}$ the (already found) tilings
between $\text{Min}_{(C_i)}^{k-1}(\cal D)$ and ${\cal D}^{k}$; for each of
these tilings~$t$, add $t + (C_i)$ to~$L_{\cal D}$.

\end{itemize}

\item \textbf{Step 4}: For each tiling~$T$ in~$L_{\cal D}$, compute the height
function and find the local minima: they correspond to cubes that can be
added. Connect~$T$ to each of the tilings obtained from it by adding exactly
one of these cubes. This generates the lattice of the tilings of the fertile
zone~$\cal D$.

\item \textbf{Step 5}: Compute the product of the lattices found for
each fracture zone.

\end{itemize}

\bigskip

An practical implementation of this algorithm should update the lattice of the
fertile zone at each newly found tiling.

\medskip

\noindent \textit{Proof of the algorithm}: Let $T$ be any tiling of a fertile
zone. It can be associated with a unique pile of cubes. Let $M$ be the maximum
of the orders of these cubes. The cubes of order $M$ belong to the fillings of
the seeds of order $M$, and all of them have been taken care of by
construction. The cubes of order at most $M-1$ in $T$ form a pile of cubes
which has been generated in step 3.4. We conclude that $T$ has been
encountered at least once by the algorithm. Moreover, all the tilings produced
by the algorithm are distinct. Therefore we have generated exactly once each
tiling of each fertile zone of~$\cal D$, and hence exactly once each tiling of
${\cal D}$.~\mbox{}\hfill$\square$

\medskip

\noindent \textit{Space complexity}: Let $|\cal{T(D)}|$ denote the number of
tilings of the domain~$D$ (we believe no closed formula is known). Each tiling
of~$\cal D$ is generated only once; the number of links in $\cal L(D)$
starting from a particular vertex is at most the width of the lattice $\cal
L(D)$, of which we believe nothing is known except that it is trivially
bounded by $|\cal{T(D)}|$; and there is a small overhead to account for
the~${\cal D}^k$. The execution space of the algorithm is thus
$\text{O}(|{\cal T(D)}|^2)$.~\mbox{}\hfill$\square$

\medskip

\noindent \textit{Time complexity}:  The execution time of the algorithm is
controlled by steps~3 and~4, but since the average (or worst case) number of
seeds and size of their maximal ranges is (totally unknown and) highly
dependent on the shape of the domain, no non-trivial bound can be given for
the time being. A rough analysis can be conducted as follows. In step~3.3, one
needs to compute a minimal tiling for a familiy of cubes, of which there is at
most the number of cubes in the maximal tiling of~$\cal D$, which in turn is
less than the number ${\cal T(D)}$ of tilings of~$\cal D$ (since each cube is
bijectively related to the associated minimal tiling). The minimal tiling for
each cube can be computed in a number of steps that is again at most the
number of cubes in the maximal tiling of~$\cal D$. Looking up the already
found tilings in step~3.4 can be done in at most~$|\cal T(D)|$ operations.  
Since steps 3.3 and~3.4 must be conducted for each tiling, the overall cost of
step~3 is $O(|{\cal T(D)}|^3)$. In step~4, computing the height function
requires $|\cal D|$ operations for each tiling and connecting a tiling to its
fathers in the lattice requires at most $\cal T(D)$ operations, so that the
global cost of step~4 is~$O(|{\cal T(D)}|^2 \times |\cal D|)$.

Thus an upper bound for the time complexity is $O(|{\cal T(D)}|^3)$, but the
execution time is probably far less in practice.~\mbox{}\hfill$\square$

%%%%%%%%%%%%%%%%%%%%%%%%%%%%%%%%%%%%%%%%%%%%%%%%%%%%%%%%%%%%%%%%%%%%%%%
%%%%%%%%%%%%%%%%%%%%%%%%%%%%%%%%%%%%%%%%%%%%%%%%%%%%%%%%%%%%%%%%%%%%%%%
%%%%%%%%%%%%%%%%%%%%%%%%%%%%                %%%%%%%%%%%%%%%%%%%%%%%%%%%
%%%%%%%%%%%%%%%%%%%%%%%%%%%%   Conclusion   %%%%%%%%%%%%%%%%%%%%%%%%%%%
%%%%%%%%%%%%%%%%%%%%%%%%%%%%                %%%%%%%%%%%%%%%%%%%%%%%%%%%
%%%%%%%%%%%%%%%%%%%%%%%%%%%%%%%%%%%%%%%%%%%%%%%%%%%%%%%%%%%%%%%%%%%%%%%
%%%%%%%%%%%%%%%%%%%%%%%%%%%%%%%%%%%%%%%%%%%%%%%%%%%%%%%%%%%%%%%%%%%%%%%

\section{Conclusion and perspectives}
	\worklabel{conclusion:perspectives}

We have shown that the isomorphism between Conway and Lagarias' lozenge group
and $\mathbb{Z}^3$ is quite fruitful; it justifies the intuitive geometric
interpretation of lozenge tilings and allows us notably to define the proper
seeds and the important intermediate tilings ${\cal D}^k$, of which Thurston's
minimal and maximal tilings of ${\cal D}$ are particular cases.

Also, we have made much use of the lattice structure; this seems to indicate
that lattice theory is a tool well adapted to the study of tilings.

We believe that the techniques developped in this paper may be extended to
study some related questions:

\begin{itemize}

\item Enumerating the tilings of $\cal D$, without generating them, is a
difficult question; Algorithm~\ref{algo:main} can be used to exhibit an exact
(but not very useful) formula, which is turn could provide good bounds if one
were able to count the number of plane partitions that are both greater and
smaller than two given plane partitions. Cruder bounds can be obtained using
the lattice structure.

\item In the case of tilings with dominoes there is a natural definition of
flips, and therefore of seeds. The range of a seed~$s$ in a tiling~$T$ can be
defined as follows: mark all the seeds in~$T$ except~$s$ as forbidden (they
should not be flipped) and perform all the possible up-flips. The range of~$s$
is the collection of squares around the vertices that have undergone a flip in
the process. It seems therefore likely that our algorithms could be adapted
straightforwardly to the case of dominoes. The geometric interpretation might
be preserved using Levitov tiles (see~\cite{Randall}).

\item Our definition of (the maximal tiling of) a pseudo-hexagon is a natural
generalization of a Ferrers diagram since it is merely a geometrical
representation of a plane partition.  Let us now attribute an integer to each
cube in such a way that sequences are non-increasing along each axis: this is
a generalization of a Ferrers diagram in dimension 4 (this is called a solid
partition in \cite{MacMahon}). Thus one can easily define a generalization of
pseudo-hexagons in dimension $p \geqslant 2$ and may therefore be able to
study tilings in this dimension. In particular, the recursive algorithms of
Section~\ref{section:algo:hexagons} should be readily upgradable.

\item If more counting results were known, one could use
Algorithms~\ref{algo:main} and~\ref{algo:tilings} to generate tilings
uniformly at random.

\item The $C$-minimal tilings seem to correspond to the sup-irreducible
elements of the lattice of the tilings. It would therefore be interesting to
study the link (with Birkhoff's theorem) between these tilings and the order
associated with the lattice.

\end{itemize}

%%%%%%%%%%%%%%%%%%%%%%%%%%%%%%%%%%%%%%%%%%%%%%%%%%%%%%%%%%%%%%%%%%%%%%%%
%%%%%%%%%%%%%%%%%%%%%%%%%%%%%%%%%%%%%%%%%%%%%%%%%%%%%%%%%%%%%%%%%%%%%%%%
%%%%%%%%%%%%%%%%%%%%%%%%%%%%%               %%%%%%%%%%%%%%%%%%%%%%%%%%%%
%%%%%%%%%%%%%%%%%%%%%%%%%%%%%   La Biblio   %%%%%%%%%%%%%%%%%%%%%%%%%%%%
%%%%%%%%%%%%%%%%%%%%%%%%%%%%%               %%%%%%%%%%%%%%%%%%%%%%%%%%%%
%%%%%%%%%%%%%%%%%%%%%%%%%%%%%%%%%%%%%%%%%%%%%%%%%%%%%%%%%%%%%%%%%%%%%%%%
%%%%%%%%%%%%%%%%%%%%%%%%%%%%%%%%%%%%%%%%%%%%%%%%%%%%%%%%%%%%%%%%%%%%%%%%


\begin{thebibliography}{main}

\small

\bibitem[And78]{Andrews} George E. Andrews, \textit{Percy Alexander MacMahon,
Collected Papers}, volume~1, MIT Press, 1978, ISBN 0-262-13121-8.

\bibitem[CL90]{Conway} J. H. Conway and J. C. Lagarias, \textit{Tiling with
Polyominoes and Combinatorial Group Theory}, Journal of Combinatorial Theory,
A~53 (1990), p.~183-208.

\bibitem[Els84]{Elser} Veit Elser, \textit{Solution of the dimer problem on a
hexagonal lattice with boundary}, J.~Phys.~A: Math.~Gen. \textbf{17} (1984)
1509-1513.

\bibitem[Fou96]{Fournier3} J. C. Fournier, \textit{Pavage des figures planes
sans trous par des dominos: fondement graphique de l'algorithme de Thurston,
parall\'elisation, unicit\'e de d\'ecomposition}, Theoretical Computer Science
159 (1996), 105-128.

\bibitem[Fou97]{Fournier1} J. C. Fournier, \textit{Tiling pictures of the
plane with dominoes}, Discrete mathematics 165$/$166 (1997) 313-320.

\bibitem[Fou01]{Fournier2} J. C. Fournier, \textit{Combinatorial of perfect
matchings in plane bipartite graphs and applications to tilings} in the
present volume.

\bibitem[Fra94]{Fraleigh} John B. Fraleigh, \textit{A First Course in Abstract 
Algebra}, 5${}^{\text{th}}$ edition, Addison-Wesley, 1994, ISBN 0-201-59291-6

\bibitem[FS96]{Flajolet} Philippe Flajolet, Robert Sedgewick, \textit{An
introduction to the Analysis of Algorithms}, Addison-Wesley, 1996, ISBN
0-201-40009-X.

\bibitem[FS]{Flaj} Philippe Flajolet, Robert Sedgewick, \textit{Analytic 
Combinatorics}, to be published; chapters 1-3 and 4-5 can be found at\newline
\texttt{http://algo.inria.fr/flajolet/Publications/anacombi1.ps.gz}\newline
\texttt{http://algo.inria.fr/flajolet/Publications/anacombi2.ps.gz}.

\bibitem[Mac16]{MacMahon} Percy Alexander MacMahon, \textit{Combinatory
Analysis}, volume~2, 1916.

\bibitem[R\'em99]{Remila} E.~R\'emila, \textit{The lattice structure of the
set of domino tilings of a polygon}, \'Ecole Normale Sup\'erieure de Lyon,
Research Report N${}^{\circ}$ 1999-25.

\bibitem[RY00]{Randall} Dana Randall and Gary Yngve, \textit{Random
three-dimensional tilings of Aztec octahedra and tetrahedra: an extension of
domino tilings}, Proceedings of the Eleventh Annual ACM-SIAM Symposium on
Discrete Algorithms (San Francisco, CA, 2000), p.~636--645, ACM, New York,
2000.

\bibitem[Thu90]{Thurston} W. P. Thurston, \textit{Conway's Tiling Groups}, 
American Mathematical Monthly, oct.~1990, p.~757-773.

\end{thebibliography}
\end{document}